\documentclass[times]{cnmauth}

\def\volumeyear{2017}

\usepackage{amsmath}
\usepackage{amsfonts}
\usepackage{bm}
\usepackage{cite}
\usepackage{color}
\usepackage{colortbl}
\usepackage{epsfig}
\usepackage{multirow}
\usepackage{psfrag}
\usepackage{url}

\renewcommand \d [2]{\frac{{\rm d}    #1}{{\rm d}  #2}}
\newcommand   \D [2]{\frac{\partial   #1}{\partial #2}}

\renewcommand{\vec}[1]{\bm{\mathrm{#1}}}

\def \Div{\nabla_h \cdot \mbox{}}
\def \Grad{\nabla_h}
\def \Lap{\nabla_h^2}
\def \dt{\Delta t}
\def \dx{\Delta x}

\def \CC{\mathbb{C}}
\def \EE{\mathbb{E}}
\def \FF{\mathbb{F}}
\def \II{\mathbb{I}}
\def \PP{\mathbb{P}}
\def \PPe{\PP^{\text{e}}}
\def \RR{\mathbb{R}}

\def \E{\vec{E}}
\def \F{\vec{F}}
\def \G{\vec{G}}
\def \N{\vec{N}}
\def \T{\vec{T}}
\def \U{\vec{U}}
\def \V{\vec{V}}
\def \X{\vec{\chi}}

\def \f{\vec{f}}
\def \g{\vec{g}}
\def \n{\vec{n}}
\def \s{\vec{X}}
\def \t{\vec{t}}
\def \u{\vec{u}}
\def \v{\vec{v}}
\def \x{\vec{x}}
\def \grad{\nabla}
\def \half{\frac{1}{2}}
\def \p{\partial}

\def \Dx{\mbox{d} \x}
\def \Ds{\mbox{d} \s}
\def \DA{\mbox{d} A(\s)}

\def \cJ{\mathcal J}
\def \cM{\mathcal M}
\def \cS{\mathcal S}
\def \cT{\mathcal T}

\def \sigmaf{\vec{\sigma}^{\text{f}}}
\def \sigmae{\vec{\sigma}^{\text{e}}}
\def \fe{\f^{\text{e}}}
\def \We{W^{\text{e}}}

\def \Mfac{M_{\text{fac}}}

\begin{document}

\runningheads{B.~E.~Griffith and X.~Y.~Luo}{Finite difference/finite element immersed boundary method}

\title{Hybrid finite difference/finite element immersed boundary
  method}

\author{Boyce E.~Griffith\affil{1}\corrauth, Xiaoyu Luo\affil{2}}
\address{\affilnum{1}Departments of Mathematics and Biomedical Engineering, Carolina Center for Interdisciplinary Applied Mathematics, and McAllister Heart Institute, University of North Carolina, Chapel Hill, NC, USA\break
\affilnum{2}School of Mathematics and Statistics, University of Glasgow, Glasgow, UK}

\corraddr{Phillips Hall, Campus Box 3250, University of North Carolina, Chapel Hill, NC 27599-3250, USA. E-mail: \url{boyceg@email.unc.edu}}

\begin{abstract}

The immersed boundary method is an approach to fluid-structure interaction that uses a Lagrangian description of the structural deformations, stresses, and forces along with an Eulerian description of the momentum, viscosity, and incompressibility of the fluid-structure system.
The original immersed boundary methods 
described immersed elastic structures using systems of flexible fibers, 
and even now, most immersed boundary methods still require Lagrangian meshes that are finer than the Eulerian grid.
This work introduces a coupling scheme for the immersed boundary method to link the Lagrangian and Eulerian variables that facilitates independent spatial discretizations for the structure and background grid.
This approach employs a finite element discretization of the structure while retaining a finite difference scheme for the Eulerian variables.
We apply this method to benchmark problems involving elastic, rigid, and actively contracting structures, including an idealized model of the left ventricle of the heart.
Our tests include cases in which, for a fixed Eulerian grid spacing, coarser Lagrangian structural meshes yield discretization errors that are as much as several orders of magnitude smaller than errors obtained using finer structural meshes.
The Lagrangian-Eulerian coupling approach developed in this work enables the effective use of these coarse structural meshes with the immersed boundary method.
This work also contrasts two different weak forms of the equations, one of which is demonstrated to be more effective for the coarse structural discretizations facilitated by our coupling approach.
\end{abstract}

\keywords{immersed boundary method; fluid-structure interaction; finite element method; finite difference method; incompressible elasticity; incompressible flow}

\maketitle

\section{Introduction}
\label{s:introduction}

Since its introduction \cite{Peskin72b,Peskin77}, the immersed boundary (IB) method has been widely used to simulate biological fluid dynamics and other problems in which a structure is immersed in a fluid flow \cite{Peskin02}.
The IB formulation of such problems uses a Lagrangian description of the deformations, stresses, and forces of the structure and an Eulerian description of the momentum, viscosity, and incompressibility of the fluid-structure system.
Lagrangian and Eulerian variables are coupled by integral transforms with delta function kernels.
When the continuous equations are discretized, the Lagrangian equations are approximated on a curvilinear mesh, the Eulerian equations are approximated on a Cartesian grid, and the Lagrangian-Eulerian interaction equations are approximated by replacing the singular kernels by regularized delta functions.
A major advantage of this approach is that it permits nonconforming discretizations of the fluid and the immersed structure.
Specifically, the IB method does not require dynamically generated body-fitted meshes, a property that is especially useful for problems involving large structural deformations or displacements, or contact between structures.

In many applications of the IB method, the elasticity of the immersed structure is described by systems of fibers that resist extension, compression, or bending \cite{Peskin02}.
Such descriptions can be well suited for highly anisotropic materials commonly encountered in biological applications, and have facilitated significant work in biofluid dynamics \cite{Miller04,Miller05,FogelsonGuy08,XYang08,CYHsu09,LAMiller09,EDTytell10,TSkorczewski14-platelet_adhesion,AHoover2015-benefits,SKJones15-lift_vs_drag,WKou15-esophageal_transport,WKou15-muscles_in_esophageal_transport,EDTytell16-insights}, including three-dimensional simulations of cardiac fluid dynamics \cite{PeskinMcQueen96,McQueenPeskin97,Lemmon00a,Lemmon00b,McQueenPeskin00,McQueenPeskin01,BEGriffith07-ibamr_paper,BEGriffith09-heart_valves,BEGriffith09-ibamr_chapter,BEGriffith12-aortic_valve,XYLuo12,XSMa13-human_mv}.
Fiber models are also convenient to use in practice because they permit an especially simple discretization as collections of points that are connected by springs or beams.
The fiber-based approach to elasticity modeling also presents certain challenges.
For instance, it can be difficult to incorporate realistic shear properties into spring network models \cite{PEHammer11-mass_spring}.
Further, fiber models often must use extremely dense collections of Lagrangian points to avoid leaks if the models are subjected to very large deformations \cite{WKou15-esophageal_transport}.

The fiber-based elasticity models often used with the conventional IB method are a special case of finite-deformation structural models in which the structural response depends only on strains in a single material direction (i.e., strains in the \emph{fiber direction}).
The mathematical framework of the IB method is not restricted to fiber-based material models, however, and several distinct extensions of this methodology have sought to use more general finite-deformation structural mechanics models.
For instance, Liu, Wang, Zhang, and co-workers developed the immersed finite element (IFE) method \cite{LZhang04,WKLiu06,LTZhang07,XSWang10,XSWang12-semi_implicit_IFE}, which is a version of the IB method in which finite element (FE) approximations are used for both the Lagrangian and the Eulerian equations.
Like the IB method, the IFE method couples Lagrangian and Eulerian variables by discretized integral transforms with regularized delta function kernels, although because the IFE method is designed to use unstructured discretizations of the Eulerian momentum equation, the IFE method uses different families of smoothed kernel functions than those typically used with the IB method.
Devendran and Peskin proposed an energy functional-based version of the conventional IB method that obtains a nodal approximation to the elastic forces generated by an immersed hyperelastic material via an FE-type approximation to the deformation of the material \cite{Peskin-energy-functions,DDevendran12}, and Boffi, Costanzo, Gastaldi, Heltai, and co-workers developed a fully variational IB method that avoids regularized delta functions altogether \cite{DBoffi08,Heltai12,SRoy15-benchmarking_ibfe}.
Other work led to the development of the immersed structural potential method, which uses a meshless method to describe the mechanics of 
hyperelastic structures immersed in fluid \cite{AJGil10,AJGil13}.

In this paper, we describe an alternative approach to using FE structural discretizations with the IB method that combines a Cartesian grid finite difference method for the incompressible Navier-Stokes equations with a nodal FE method for the structural mechanics.
We consider both flexible structures with a hyperelastic material response, and also rigid immersed structures in which rigidity constraints are approximately imposed via a simple penalty method.
The primary contribution of this work is its treatment of the equations of Lagrangian-Eulerian interaction.
Conventionally, structural forces are \emph{spread} directly from the nodes of the Lagrangian mesh using the regularized delta function, and velocities are \emph{interpolated} directly from the Eulerian grid to the Lagrangian mesh nodes.
This discrete Lagrangian-Eulerian coupling approach is adopted by the classical IB method as well as the IFE method \cite{LZhang04,WKLiu06,LTZhang07,XSWang10,XSWang12-semi_implicit_IFE}, the energy-based method of Devendran and Peskin \cite{Peskin-energy-functions,DDevendran12}, and the immersed structural potential method \cite{AJGil10,AJGil13}.
A significant limitation of this approach is that if the physical spacing of the Lagrangian nodes is too large in comparison to the spacing of the background Eulerian grid, severe leaks will develop at fluid-structure interfaces.
A widely used empirical rule that generally prevents such leaks is to require the Lagrangian mesh to be approximately twice as fine as the Eulerian grid \cite{Peskin02}, potentially requiring very dense structural meshes, especially for applications that involve large structural deformations.

The Lagrangian-Eulerian coupling scheme introduced in this paper overcomes this longstanding limitation of the classical IB method.
Specifically, rather than spreading forces from the nodes of the Lagrangian mesh and interpolating velocities to those mesh nodes, we instead spread forces from and interpolate velocities to dynamically selected \emph{interaction points} located within the Lagrangian structural elements.
Here, the interaction points are constructed by quadrature rules defined on the structural elements.
In contrast to both the conventional IB method and the IFE method, in our method, the nodes of the structural mesh act as \emph{control points} that determine the positions of the interaction points, but the mesh nodes are not required to be locations of direct Lagrangian-Eulerian interaction.\footnote{
	In fact, the structural mesh nodes can be both control points and interaction points if the quadrature rules that generate the interaction points include the mesh nodes as quadrature points, 
	as in Gauss-Lobatto rules.
}
Our approach thereby takes full advantage of the kinematic information provided by the FE description of the structural deformation, which yields approximations not only to the positions of the nodes of the Lagrangian mesh, but also to the positions of all material points of the structure.

A related IB method was introduced by Shankar et al.~\cite{VShankar15-IB-RBF}, who use a radial basis function (RBF) approach to describe the deforming geometry of two-dimensional models of circulating platelets, which are described as elliptical shells with linear elasticity models.
In that approach, velocities are interpolated to \emph{data sites}, which are analogous to the control points of the present FE scheme, but forces are spread from a fixed collection of \emph{sample sites} that are distributed along the surface of the platelet.
Unlike the RBF scheme, the interaction operators of the present FE approach satisfy a discrete adjoint property that implies that the method satisfies energy conservation during Lagrangian-Eulerian interaction.
As demonstrated by Shankar et al., satisfying the discrete adjoint property is not required to obtain a convergent scheme, but this property is crucial for some applications, such as the design of efficient IB methods for rigid bodies \cite{BKallemov16-RigidIBAMR, FBalboaUsabiaga16-rigid_multiblob}.
The importance of satisfying the adjoint property also seems to depend on the choice of regularized kernel functions.
Lagrangian-Eulerian coupling schemes that do not satisfy the adjoint property seem to be prone to induce numerical instabilities when used with kernel functions developed by Peskin that satisfy moment conditions along with a ``sum-of-squares'' condition \cite{Peskin02}, like those used in this work.
Shankar et al.~consider a kernel function that may be less sensitive to these aspects of the discretization.

We apply the present IB method to benchmark fluid-structure interaction (FSI) problems involving elastic, rigid, and actively contracting structures, including an idealized model of the left ventricle of the heart \cite{SLand15-cardiac_verification}.
For elastic structures, we consider two weak formulations of the equations of motion suitable for standard nodal ($C^0$) FE methods for structural mechanics.
One of these formulations, referred to herein as the \emph{unified weak form}, is similar to those used by earlier IB-like methods \cite{LZhang04,WKLiu06,LTZhang07,XSWang10,XSWang12-semi_implicit_IFE,DBoffi08,Heltai12,SRoy15-benchmarking_ibfe}.
This formulation uses a single volumetric force density to describe the mechanical response of the immersed structure.
The other formulation, referred to herein as the \emph{partitioned weak form}, uses both an internal volumetric force density, which is supported throughout the immersed elastic structure, and also a transmission surface force density, which is supported only on the surface of the immersed structure.
In a numerical scheme, the unified formulation effectively \emph{regularizes} the transmission surface force density by projecting it onto the volumetric FE shape functions.
The partitioned formulation described herein, which does not appear to have been used in previous versions of the IB or IFE methods, avoids this additional regularization step by treating the surface and volume force densities separately.
We present numerical tests that demonstrate that this partitioned scheme can yield accurate results even for Lagrangian meshes that are significantly \emph{coarser} than the background Eulerian grid. 

An important feature of our discretization approach is that obtaining a ``watertight'' structure simply requires using a dense enough collection of interaction points so as to prevent leaks. 
Moreover, because it is straightforward to use dynamic quadrature schemes that account for highly deformed elements, 
this approach can ensure that the structure remains watertight even in the presence of very large structural deformations.
One benefit of our approach is that the Lagrangian resolution may be determined primarily by accuracy requirements for the structural model rather than by the requirements of the Lagrangian-Eulerian coupling scheme.
Further, for problems involving immersed rigid structures, we demonstrate that using coarser Lagrangian meshes can reduce discretization errors by an order of magnitude compared to finer Lagrangian meshes for a given Eulerian grid spacing.
The present method also yields improved volume conservation in comparison to the IFE method \cite{XSWang10} for both coarse and fine structural meshes.
To our knowledge, the present IB method is the first IFE-type method to explicitly enable the effective use of such relatively coarse Lagrangian meshes.


\section{Continuous formulations}
\label{s:continuous_formulations}

\subsection{Immersed elastic bodies}
\label{s:immersed_elastic_structures}

\begin{figure}
\centering
\input{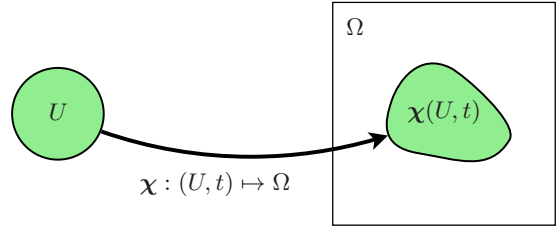}
\caption{
	Lagrangian and Eulerian coordinate systems.
	The Lagrangian material coordinate domain is $U$, and the Eulerian physical coordinate domain is $\Omega$.
	The physical position of material point $\s$ at time $t$ is $\X(\s,t)$, the physical region occupied by the structure is $\X(U,t)$, and the physical region occupied by the fluid is $\Omega \setminus \X(U,t)$.
}
\label{f:coordinates}
\end{figure}

In the IB formulation of problems involving an immersed elastic body, the momentum, velocity, and incompressibility of the coupled fluid-structure system are described in Eulerian form, whereas the deformation and elastic response of the immersed structure are described in Lagrangian form.
A similar formulation is used in the case of an immersed \emph{rigid} structure; see sec.~\ref{s:immersed_rigid_structures}.
Let \mbox{$\x = (x_1,x_2,\ldots) \in \Omega \subset \RR^d$}, $d=2$ or $3$, denote Cartesian physical coordinates, with $\Omega$ denoting the physical region that is occupied by the coupled fluid-structure system; let $\s = (X_1,X_2,\ldots) \in U \subset \RR^d$ denote Lagrangian material coordinates that are attached to the structure, with $U$ denoting the Lagrangian coordinate domain; and let $\X(\s,t) \in \Omega$ denote the physical position of material point $\s$ at time $t$.
The physical region occupied by the structure at time $t$ is $\X(U,t) \subseteq \Omega$, and the physical region occupied by the fluid at time $t$ is $\Omega \setminus \X(U,t)$.
See fig.~\ref{f:coordinates}.
We do not assume that the Lagrangian coordinates are the initial coordinates of the elastic structure, nor, more generally, do we require that $U \subseteq \Omega$.

To use an Eulerian description of the fluid and a Lagrangian description of the elasticity of the immersed structure, it is necessary to describe the stress of the fluid-structure system in both Eulerian and Lagrangian forms.
Let $\vec{\sigma}(\x,t)$ denote the Cauchy stress tensor of the \emph{coupled fluid-structure system}.
In the present formulation, we assume that
\begin{equation}
  \vec{\sigma}(\x,t) = \sigmaf(\x,t) +
  \begin{cases}
    \sigmae(\x,t) & \mbox{for $\x \in \X(U,t)$,} \\
    \vec{0} & \mbox{otherwise,}
  \end{cases}
  \label{e:cauchy_stress_tensor}
\end{equation}
in which $\sigmaf(\x,t)$ is the stress tensor of a viscous incompressible fluid and $\sigmae(\x,t)$ is the stress tensor that describes the elastic response of the immersed structure.
The fluid stress tensor is the usual one for a viscous incompressible fluid,
\begin{equation}
  \sigmaf(\x,t) = -p(\x,t) \, \II + \mu \left[ \grad\vec{u}(\x,t) + \left(\grad\vec{u}(\x,t)\right)^{\mathrm T} \right], \label{e:fluid_stress_tensor}
\end{equation}
in which $p(\x,t)$ is the pressure, $\mu$ is the dynamic viscosity, and $\u(\x,t)$ is the Eulerian velocity field.

To describe the elasticity of the structure with respect to the Lagrangian material coordinate system, we use the first Piola-Kirchhoff elastic stress tensor $\PPe(\s,t)$, which 
is defined in terms of $\sigmae$ via
\begin{equation}
  \PPe(\s,t) = J(\s,t) \, \sigmae(\X(\s,t),t) \, \FF^{-{\mathrm T}}(\s,t),  \label{e:PK1_stress_alt}
\end{equation}
in which the deformation gradient tensor associated with the deformation $\X:(U,t)\mapsto\Omega$ is
\begin{equation}
	\FF(\s,t) = \D{\X}{\s}(\s,t),
\end{equation}
and $J(\s,t) = \det(\FF(\s,t))$ is the Jacobian determinant of the deformation gradient.
Although $\PPe$ is only defined within the solid region, it is convenient to extend $\sigmae(\x,t)$ by zero in the fluid region.

For simplicity, we primarily restrict our attention to hyperelastic constitutive models, which may be characterized by a strain-energy functional $\We(\FF)$.
For such constitutive laws,
\begin{equation}
	\PPe = \D{\We}{\FF}.
\end{equation}
Our formulation does not rely on the existence of such an energy functional, however, and it permits a material description defined only in terms of a stress response.
For instance, separate work using the present methodology relies on this feature 
to treat active tension generation in dynamic models of muscle contraction \cite{HGao14-iblv_mi,SLand15-cardiac_verification,WWChen16-coupled_lv,HGaoXX-coupled_mv_lv}.

\subsubsection{Strong formulation}

As shown by Boffi et al.~\cite{DBoffi08}, the strong form of the equations of motion is
\begin{align}
  \rho\frac{\mathrm{D}\u}{\mathrm{D}t}(\x,t) &= - \grad p(\x,t) + \mu \grad^2 \u(\x,t) + \fe(\x,t),         \label{e:momentum}          \\
  \grad \cdot \u(\x,t) &= 0,                                                                                \label{e:incompressibility} \\
  \fe(\x,t)            &= \int_U \grad_{\s} \cdot \PPe(\s,t) \, \delta(\x - \X(\s,t)) \, \Ds                \label{e:force_density}     \\
                       & \ \ \ \mbox{} - \int_{\p U} \PPe(\s,t) \, \N(\s) \, \delta(\x - \X(\s,t)) \, \DA,  \nonumber                   \\
  \D{\X}{t}(\s,t)      &= \int_\Omega \u(\x,t) \, \delta(\x - \X(\s,t)) \, \Dx,            \label{e:interp}
\end{align}
in which $\rho$ is the mass density of the coupled fluid-structure system, $\frac{\mathrm{D}\u}{\mathrm{D}t}(\x,t) = \D{\u}{t}(\x,t) + \u(\x,t) \cdot \grad \u(\x,t)$ is the material derivative, $\fe(\x,t)$ is the Eulerian elastic force density, and $\delta(\x) = \prod_{i=1}^{d} \delta(x_i)$ is the $d$-dimensional delta function.

Two different types of Lagrangian elastic force densities appear in these equations.
The Lagrangian \emph{internal force density}, $\grad_{\s} \cdot \PPe$, is a volumetric force density that is distributed throughout the elastic body, whereas the Lagrangian \emph{transmission force density}, $-\PPe \, \N$, is a surface force density that is distributed along the fluid-solid interface.
Eq.~\eqref{e:force_density} generates a corresponding Eulerian description of the elastic forces from the volumetric and surface forces.
Notice that $\fe(\x,t)$ will generally have a $\delta$-layer of force along the fluid-solid interface.
It is possible to show that $\fe$ is variationally equivalent to $\grad \cdot \sigmae$ (for instance, by integrating against a test function).
Doing so, it is clear that $\fe$ is singular along the fluid-solid interface wherever $\sigmae \n$ is discontinuous.
Indeed, because elastic stresses are present only within the structure, $\sigmae \n$ is generally discontinuous at the fluid-structure interface.
These discontinuities must be exactly balanced by discontinuities in $\sigmaf \n$ to ensure that the \emph{total} stress, $\vec{\sigma} = \sigmaf + \sigmae$, has a continuous traction vector.
If $\PPe(\s,t)$ is sufficiently smooth, the internal force acts as a regular (i.e., nonsingular) body force and may be treated with higher-order accuracy by the IB method \cite{BEGriffith05-ib_accuracy,DBoffi08}.
However, the transmission force always acts as a singular force layer, and although this force will induce jumps in the pressure and shear stress along $\p \X(U,t)$, such force layers may also be readily treated by the IB method.

An integral transform is also used in eq.~\eqref{e:interp} to determine the velocity of the immersed elastic structure from the material velocity field $\u(\x,t)$.
The defining property of $\delta(\x)$ implies that eq.~\eqref{e:interp} is equivalent to
\begin{equation}
	\D{\X}{t}(\s,t) = \u(\X(\s,t),t),
\end{equation}
which may be interpreted as the no-slip and no-penetration conditions of a viscous incompressible fluid.
Notice, however, that the no-slip and no-penetration conditions do not appear as constraints on the fluid motion.
Instead, these conditions determine the motion of the immersed structure.

\subsubsection{Weak formulations}

To use standard nodal ($C^0$) FE methods for nonlinear elasticity with the IB formulation, it is necessary to introduce a weak formulation of the equations of motion.
Here we consider two different formulations that each employ a weak form of the Lagrangian equations.
Because we use a finite difference scheme to approximate the incompressible Navier-Stokes equations, we do not develop a weak formulation for the Eulerian equations or the equations of Lagrangian-Eulerian interaction.

We begin with the \emph{partitioned weak formulation} of the equations.
To do so, we first define the internal elastic force density, $\F(\s,t)$, that is variationally equivalent to $\grad_{\s} \cdot \PPe$ by requiring
\begin{align}
  \int_U \F(\s,t) \cdot \V(\s) \, \Ds =& -\int_U \PPe(\s,t) : \grad_{\s} \V(\s) \, \Ds \label{e:w2_virtual_work} \\
  & \ \ \ \ \ \mbox{} + \int_{\p U} \left(\PPe(\s,t) \, \N(\s)\right) \cdot \V(\s) \, \DA  \nonumber
\end{align}
to hold for arbitrary smooth $\V(\s)$.  
Integrating by parts, it is clear that eq.~\eqref{e:w2_virtual_work} implies that
\begin{equation}
   \int_U \F(\s,t) \cdot \V(\s) \, \Ds = \int_U \left(\grad_{\s} \cdot \PPe(\s,t)\right) \cdot \V(\s) \, \Ds,
\end{equation}
for all smooth $\V(\s)$.
Assuming sufficient regularity, $\F = \grad_{\s} \cdot \PPe$ pointwise.
It is also convenient to define the transmission force density, $\T(\s,t)$, pointwise along $\p U$ via
\begin{equation}
  \T(\s,t) = -\PPe(\s,t) \, \N(\s). \label{e:w2_transmission_force} 
\end{equation}
With $\F$ and $\T$ so defined, the equations of motion become
\begin{align}
  \rho\frac{\mathrm{D}\u}{\mathrm{D}t}(\x,t) &= -\grad p(\x,t) + \mu \grad^2 \u(\x,t) + \f(\x,t) + \t(\x,t),  \label{e:w2_momentum}                   \\
  \grad \cdot \u(\x,t)                       &= 0,                                                            \label{e:w2_incompressibility}          \\
  \f(\x,t)                                   &= \int_U \F(\s,t) \, \delta(\x - \X(\s,t)) \, \Ds,              \label{e:w2_internal_force_density}     \\
  \t(\x,t)                                   &= \int_{\p U} \T(\s,t) \, \delta(\x - \X(\s,t)) \, \DA,         \label{e:w2_transmission_force_density} \\
  \D{\X}{t}(\s,t)                            &= \int_\Omega \u(\x,t) \, \delta(\x - \X(\s,t)) \, \Dx,         \label{e:w2_velocity_interp}
\end{align}
in which $\f(\x,t)$ is the Eulerian internal elastic force density and $\t(\x,t)$ is the Eulerian transmission elastic force density.
It is important to notice that, under relatively mild regularity requirements, $\F$ and $\T$ are both smooth functions on their domains of definition, and $\f$ is piecewise smooth.
Under these conditions, the only singular function in this formulation is $\t$.

Alternative weak definitions of the Lagrangian elastic force density are possible.
The formulation typically used with the IFE method defines a \emph{total} elastic force per unit volume, $\G(\s,t)$, by requiring
\begin{equation}
  \int_U \G(\s,t) \cdot \V(\s) \, \Ds = -\int_U \PPe(\s,t) : \grad_{\s} \V(\s) \, \Ds  \label{e:w1_virtual_work}
\end{equation}
for arbitrary $\V(\s)$.
It is possible to see that $\G(\s,t)$ accounts for the effects of \emph{both} the internal and transmission force densities of the strong form of the equations.
To do so, we again integrate by parts to find that
\begin{align}
  \int_U \G(\s,t) \cdot \V(\s) \, \Ds =& \int_U \left(\grad \cdot \PPe(\s,t) \right) \cdot \V(\s) \, \Ds \\
  & \ \ \ \ \ \mbox{} - \int_{\p U} \left(\PPe(\s,t) \, \N(\s)\right) \cdot \V(\s) \, \DA, \nonumber
\end{align}
for all $\V(\s)$.
Thus, $\G(\s,t)$ can be a continuous function only if $\PPe \N \equiv \vec{0}$.
In general, $\G$ is in fact a distribution that accounts for both the internal force per unit volume and the transmission force per unit area, in which the transmission force gives rise to a singular force layer concentrated along $\p U$.
It is possible to use eq.~\eqref{e:w1_virtual_work} with standard finite element methods; however, when doing so, the transmission surface force density is effectively projected (in an $L^2$ sense) onto the volumetric basis functions.
Specifically, the FE basis functions serve to \emph{regularize} the transmission force.

Using this definition of $\G$, we may state a \emph{unified weak formulation} that includes only a single, unified body forcing term:
\begin{align}
  \rho\frac{\mathrm{D}\u}{\mathrm{D}t}(\x,t) &= -\grad p(\x,t) + \mu \grad^2 \u(\x,t) + \g(\x,t),      \label{e:w1_momentum}          \\
  \grad \cdot \u(\x,t)                       &= 0,                                                     \label{e:w1_incompressibility} \\
  \g(\x,t)                                   &= \int_U \G(\s,t) \, \delta(\x - \X(\s,t)) \, \Ds,       \label{e:w1_force_density}     \\
  \D{\X}{t}(\s,t)                            &= \int_\Omega \u(\x,t) \, \delta(\x - \X(\s,t)) \, \Dx,  \label{e:w1_velocity_interp}
\end{align}
in which $\g(\x,t)$ is the Eulerian total elastic force density.
This weak form of the equations of motion is essentially the formulation employed in the IFE method \cite{LZhang04,WKLiu06,LTZhang07,XSWang10,XSWang12-semi_implicit_IFE}, the energy-based formulation of Devendran and Peskin \cite{Peskin-energy-functions,DDevendran12}, and the fully variational IB method \cite{DBoffi08,Heltai12,SRoy15-benchmarking_ibfe}.
To our knowledge, the partitioned formulation described here has not been widely used in previous work.

\subsection{Immersed rigid structures}
\label{s:immersed_rigid_structures}

The equations of motion for a fixed, rigid structure are similar to those used in the case of an immersed elastic structure:
\begin{align}
  \rho\frac{\mathrm{D}\u}{\mathrm{D}t}(\x,t)
                       &= - \grad p(\x,t) + \mu \grad^2 \u(\x,t) + \f(\x,t),               \label{e:momentum_rigid}          \\
  \grad \cdot \u(\x,t) &= 0,                                                               \label{e:incompressibility_rigid} \\
  \f(\x,t)             &= \int_U \F(\s,t) \, \delta(\x - \X(\s,t)) \, \Ds,                 \label{e:force_density_rigid}     \\
  \D{\X}{t}(\s,t)      &= \int_\Omega \u(\x,t) \, \delta(\x - \X(\s,t)) \, \Dx = \vec{0},  \label{e:interp_rigid}
\end{align}
in which here, $\F(\s,t)$ is a Lagrange multiplier for the constraint $\D{\X}{t} \equiv \vec{0}$.
Thus, this fully constrained formulation takes the form of an extended saddle-point problem with two Lagrange multipliers, $p(\x,t)$ for the incompressibility constraint and $\F(\s,t)$ for the rigidity constraint.
Solving this system effectively requires specialized techniques that 
are the subject of active research \cite{BKallemov16-RigidIBAMR, FBalboaUsabiaga16-rigid_multiblob}.
In this work, we consider instead a penalty formulation for an immersed stationary structure, in which the Lagrange multiplier force is approximated by
\begin{equation}
  \F(\s,t) = \kappa\left(\X(\s,0) - \X(\s,t)\right) - \eta \D{\X}{t}(\s,t), \label{e:rigidity}
\end{equation}
in which $\kappa \ge 0$ is a stiffness penalty parameter and $\eta \ge 0$ is a damping penalty parameter.
As $\kappa \rightarrow \infty$, $\X(\s,t) \rightarrow \X(\s,0)$ and $\D{\X}{t} \rightarrow \vec{0}$, so that this formulation is equivalent to the constrained formulation.
In principle, it is not necessary to include the damping parameter $\eta$; however, we have found that using damping reduces numerical oscillations, especially at moderate-to-high Reynolds numbers.

\section{Spatial discretization}
\label{s:spatial_discretization}

For simplicity, we describe the numerical scheme in two spatial dimensions.
The extension of the numerical scheme to the case $d=3$ is straightforward. 

\subsection{Eulerian discretization}

\begin{figure}
	\centering
	\small
	\input{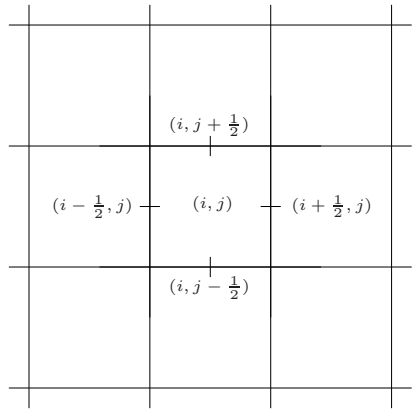}
	\caption{Schematic of the staggered-grid layout of Eulerian degrees of freedom in two spatial dimensions.
	The pressures are approximated at cell centers, indicated by $(i,j)$, the $x_1$-components of the velocity and force are approximated on the $x_1$-edges, $(i-\half,j)$ or $(i+\half,j)$, and the $x_2$-components of the velocity and force are approximated on the $x_2$-edges, $(i,j-\half)$ or $(i,j+\half)$.}	
	\label{f:cartesian_grid}
\end{figure}

We use a staggered-grid finite difference scheme to discretize the incompressible Navier-Stokes equations in space.
Compared to collocated discretizations (i.e., purely cell-~or node-centered schemes), such Eulerian discretization approaches yield superior accuracy when used with the conventional IB method \cite{BEGriffith12-ib_volume_conservation}.
To simplify the exposition, assume that $\Omega$ is the unit square and is discretized using a regular $N \times N$ Cartesian grid with grid spacings $\dx_1 = \dx_2 = h = \frac{1}{N}$.
Let $\mbox{$(i,j)$}$ label the individual Cartesian grid cells for integer values of $i$ and $j$, $0 \le i,j < N$.
The components of the Eulerian velocity field $\u = (u_1,u_2)$ are approximated at the centers of the $x_1$-~and $x_2$-edges of the Cartesian grid cells, i.e., at positions $\x_{i-\half,j} = \mbox{$\left( ih, \left(j+\half\right)h \right)$}$ and $\x_{i,j-\half} = \mbox{$\left( \left(i+\half\right)h, jh \right)$}$, respectively.
A staggered scheme is also used for the Eulerian body force $\f = (f_1,f_2)$.
We use the notation $(u_1)_{i-\half,j}$, $(u_2)_{i,j-\half}$, $(f_1)_{i-\half,j}$, and $(f_2)_{i,j-\half}$ to denote the discrete values of $\u$ and $\f$.
The pressure $p$ is approximated at the centers of the Cartesian grid cells, and its discrete values are denoted $p_{i,j}$.
See fig.~\ref{f:cartesian_grid}.

Let $\Div \u \approx \grad \cdot \u$, $\Grad p \approx \grad p$, and $\Lap \u \approx \grad^2 \u$ respectively denote standard second-order accurate finite difference approximations to the divergence, gradient, and Laplace operators \cite{BEGriffith09-efficient}.
In this approach, $\Div \u$ is defined at cell centers, whereas both $\Grad p$ and $\Lap \u$ are defined at cell edges.
We use a staggered-grid version \cite{BEGriffith09-efficient,BEGriffith12-ib_volume_conservation} of the xsPPM7 variant \cite{RiderEtAl07} of the piecewise parabolic method (PPM) \cite{ColellaWoodward84} to discretize the nonlinear advection terms.
Where needed, physical boundary conditions are discretized and imposed along $\p \Omega$ as described by Griffith \cite{BEGriffith09-efficient}.
In some numerical examples, we use a locally refined Eulerian discretization that employs Cartesian grid adaptive mesh refinement (AMR) following the discretization approach described by Griffith \cite{BEGriffith12-aortic_valve}.

\subsubsection{Eulerian inner products}

If $\u$ and $\v$ are discrete staggered-grid vector fields, we denote by $[\u]$ and $[\v]$ the corresponding vectors of grid values.
If $\Omega$ has periodic boundaries, we define the discrete $L^2$ inner product on $\Omega$ by
\begin{equation}
  \left( \u , \v \right)_{\x} = [\u]^{\text{T}} [\v] \, h^2.
\end{equation}
Minor adjustments to this definition are required when non-periodic physical boundary conditions are used \cite{BEGriffith09-efficient} or near coarse-fine interfaces in locally refined grids \cite{BEGriffith12-aortic_valve}.

\subsection{Lagrangian discretization}

Let $\cT_h = \cup_e U^{e}$ be a triangulation of $U$ composed of elements $U^e$, in which $e$ indexes the elements of the mesh.
We denote by $\{\s_l\}_{l=1}^{M}$ the nodes of the mesh, and by $\{\phi_l(\s)\}_{l=1}^{M}$ nodal (Lagrangian) FE basis functions.
The time-dependent physical positions of the nodes of the Lagrangian mesh are $\{\X_l(t)\}_{l=1}^{M}$, and, using the FE basis functions, we define an approximation to $\X(\s,t)$ by
\begin{equation}
  \X_h(\s,t) = \sum_{l=1}^{M} \X_l(t) \, \phi_l(\s).
\end{equation}
An approximation to the deformation gradient is given by
\begin{equation}
  \FF_h(\s,t) = \D{{\X}_h}{\s}(\s,t) = \sum_{l=1}^{M} \X_l(t) \, \D{\phi_l}{\s}(\s).
\end{equation}

\subsubsection{Immersed elastic structures}

For an immersed elastic structure, we use the FE approximation to the deformation gradient tensor $\FF_h(\s,t)$ to compute $\PPe_h(\s,t)$ and $\T_h(\s,t)$, which are approximations to the first Piola-Kirchhoff stress tensor and the Lagrangian transmission force density, respectively. 
We approximate the Lagrangian force densities $\F(\s,t)$ and $\G(\s,t)$ by
\begin{align}
  \F_h(\s,t) &= \sum_{l=1}^{M} \F_l(t) \, \phi_l(\s),\text{ and} \\
  \G_h(\s,t) &= \sum_{l=1}^{M} \G_l(t) \, \phi_l(\s).
\end{align}
The nodal values $\{\F_l\}_{l=1}^{M}$ and $\{\G_l\}_{l=1}^{M}$ must be determined from $\PPe_h(\s,t)$ via discretizations of eq.~\eqref{e:w2_virtual_work} and eq.~\eqref{e:w1_virtual_work}.
We use a standard Galerkin approach, so that after rearranging terms, eq.~\eqref{e:w2_virtual_work}
becomes
\begin{align}
  \sum_{l=1}^{M} \left(\int_U \phi_l(\s) \, \phi_m(\s) \, \Ds \right) \F_l(t)
  =& -\int_U \PPe_h(\s,t) \, \grad_{\s} \phi_m(\s) \, \Ds \nonumber \\
  & \ \ \ \ \ \mbox{} + \int_{\p U} \PPe_h(\s,t) \, \N(\s) \, \phi_m(\s) \, \DA, \label{e:discrete_F_l}
\end{align}
for each $m=1,\ldots,M$.  Similarly, eq.~\eqref{e:w1_virtual_work} becomes
\begin{equation}
  \sum_{l=1}^{M} \left(\int_U \phi_l(\s) \, \phi_m(\s) \, \Ds \right) \G_l(t) = -\int_U \PPe_h(\s,t) \, \grad_{\s} \phi_m(\s) \, \Ds, \label{e:discrete_G_l}
\end{equation}
for each $m=1,\ldots,M$.
In practice, these integrals are approximated via Gaussian quadrature.

\subsubsection{Immersed rigid structures}

For a fixed, rigid immersed structure, we directly evaluate the discretized Lagrangian penalty force $\F_h(\s,t)$ via
\begin{equation}
  \F_h(\s,t) = \kappa\left(\X_h(\s,0) - \X_h(\s,t)\right) - \eta \D{\X_h}{t}(\s,t).
\end{equation}
Because $\X_h(\s,t)$ and $\F_h(\s,t)$ are defined in terms of the same basis functions, $\F_h(\s,t)$ is given by
\begin{equation}
  \F_h(\s,t) = \sum_l \F_l(t) \, \phi_l(\s)
\end{equation}
in which
\begin{equation}
  \F_l(t) = \kappa\left(\X_l(0) - \X_l(t)\right) - \eta \d{\X_l}{t}(t).
\end{equation}

\subsubsection{Lagrangian inner products}

Letting $[\F]$ denote the vector of nodal coefficients of $\F_h$, we write eq.~\eqref{e:discrete_F_l} as
\begin{equation}
  [\cM] \, [\F] = [\vec{B}], \label{e:discrete_F}
\end{equation}
in which $[\cM]$ is the mass matrix that has entries of the form $\int_U \phi_l(\s) \, \phi_m(\s) \, \Ds$.
Eq.~\eqref{e:discrete_G_l} may be rewritten similarly.
The mass matrix $[\cM]$ can also be used to evaluate the $L^2$ inner product of Lagrangian functions on $U$.
In particular, for any $\U_h(\s,t) = \sum_l \U_l(t) \, \phi_l(\s)$ and $\V_h(\s,t) = \sum_l \V_l(t) \, \phi_l(s)$,
\begin{equation}
  \left( \U_h , \V_h \right)_{\s} = [\U]^{\text{T}} \, [\cM] \, [\V].  \label{e:discrete_L2_U}
\end{equation}
Different choices of mass matrices (e.g., lumped mass matrices) induce different discrete inner products on $U$.

To simplify notation, in the remainder of this paper we drop the subscript ``$h$'' from our numerical approximations to the Lagrangian variables.

\subsection{Lagrangian-Eulerian interaction}

We next describe Lagrangian-Eulerian coupling operators that take advantage of the kinematic information encoded in the FE approximation to the deformation of the immersed structure.
As in the conventional IB method, we approximate the singular delta function kernel appearing in the Lagrangian-Eulerian interaction equations by a smoothed $d$-dimensional Dirac delta function $\delta_h(\x)$ that is of the tensor-product form $\delta_h(\x) = \mbox{$\prod_{i=1}^{d} \delta_h(x_i)$}$.
Except where otherwise noted, in this work, we take the one-dimensional smoothed delta function $\delta_h(x)$ to be the four-point delta function of Peskin~\cite{Peskin02}.

To compute an approximation to $\f = (f_1,f_2)$ on the Cartesian grid, we construct for each element $U^{e} \in \cT_h$ a Gaussian quadrature rule with $N^e$ quadrature points $\s_Q^e \in U^{e}$ and weights $\omega_Q^e$, $Q = 1,\ldots,N^e$.
We then compute $f_1$ and $f_2$ on the edges of the Cartesian grid cells via
\begin{align}
  (f_1)_{i-\half,j}
  &= \sum_{U^{e} \in \cT_h} \sum_{Q=1}^{N^e} F_1(\s_Q^e,t) \, \delta_h(\x_{i-\half,j} - \X(\s_Q^e,t)) \, \omega_Q^e,\text{ and} \label{e:discrete_force_prolongation1} \\
  (f_2)_{i,j-\half}
  &= \sum_{U^{e} \in \cT_h} \sum_{Q=1}^{N^e} F_2(\s_Q^e,t) \, \delta_h(\x_{i,j-\half} - \X(\s_Q^e,t)) \, \omega_Q^e, \label{e:discrete_force_prolongation2}
\end{align}
with $\F(\s,t) = (F_1(\s,t),F_2(\s,t))$.
We use the notation
\begin{equation}
  \f(\x,t) = \cS\left(\X(\cdot,t)\right) \, \F(\s,t),  \label{e:discrete_force_prolongation}
\end{equation}
in which $\cS\left(\X(\cdot,t)\right)$ is the \emph{force-prolongation operator} implicitly defined by eqs.~\eqref{e:discrete_force_prolongation1} and \eqref{e:discrete_force_prolongation2}.
Notice that
\begin{align}
  (f_1)_{i-\half,j}
  &= \sum_{U^{e} \in \cT_h} \int_{U^{e}} F_1(\s,t) \,
  \delta_h(\x_{i-\half,j} - \X(\s,t)) \, \Ds + O(\Delta \s^q),\text{ and} \\
  (f_2)_{i,j-\half}
  &= \sum_{U^{e} \in \cT_h} \int_{U^{e}} F_2(\s,t) \,
  \delta_h(\x_{i,j-\half} - \X(\s,t)) \, \Ds + O(\Delta \s^q),
\end{align}
in which $\Delta \s$ is proportional to the Lagrangian mesh spacing and $O(\Delta \s^q)$ corresponds to quadrature error that may be controlled by the choice of numerical quadrature rules.

A corresponding \emph{velocity-restriction operator} $\cJ\left(\X(\cdot,t)\right)$ determines the motion of the nodes of the Lagrangian mesh from the Cartesian grid velocity field via
\begin{equation}
  \D{\X}{t}(\s,t) = \cJ\left(\X(\cdot,t)\right) \, \u(\x,t).  \label{e:discrete_velocity_restriction}
\end{equation}
There are many possible ways to construct $\cJ$; however, we have found that an effective approach is to require $\D{\X}{t} = \cJ \, \u$ to satisfy the discrete power identity,
\begin{equation}
  \left(\F, \D{\X}{t}\right)_{\s} = \left(\f , \u\right)_{\x}, \label{e:discrete_power_identity}
\end{equation}
which implies that the semi-discrete unified formulation conserves energy during Lagrangian-Eulerian interaction \cite{Peskin02}.
This power identity can be rewritten as
\begin{equation}
  \left(\F, \cJ \, \u\right)_{\s} = \left(\cS \, \F , \u\right)_{\x}, \label{e:adjoint_property}
\end{equation}
i.e., $\cJ = \cS^{*}$.

To construct $\cJ$ explicitly, it is convenient to use matrix notation.
Identifying $[\cS]$ and $[\cJ]$ with the matrix representations of $\cS$ and $\cJ$, we have that
\begin{align}
  [\f] &= [\cS] \, [\F],\text{ and} \\
  \d{[\X]}{t} &= [\cJ] \, [\u].
\end{align}
Eq.~\eqref{e:adjoint_property} then becomes
\begin{equation}
  [\F]^{\text{T}} \, [\cM] \, [\cJ] \, [\u] = \left([\cS] \, [\F]\right)^{\text{T}} \, [\u] \, h^2.  \label{e:discrete_adjoint_property}
\end{equation}
If eq.~\eqref{e:discrete_adjoint_property} is to hold for any $[\F]$ and $[\u]$, then $[\cJ]$ must be defined via
\begin{equation}
  [\cJ] = [\cM]^{-1} \, [\cS]^{\text{T}} \, h^2.
\end{equation}
In our time-stepping scheme, which is stated below in appendix~\ref{s:temporal_discretization}, notice that we need only to apply $[\cJ]$ to discrete velocity fields defined on the Cartesian grid.  Specifically, we do not need to compute $[\cJ]$ explicitly.

It is straightforward to show that this construction of $\cJ$ implies that $\D{\X}{t}(\s,t)$ is an approximation to the $L^2$ projection of the Lagrangian vector field $\U^{\text{IB}}(\s,t) =
(U^{\text{IB}}_1(\s,t),U^{\text{IB}}_2(\s,t))$, with
\begin{align}
  U^{\text{IB}}_1(\s,t) &= \sum_{i,j} \left(u_1\right)_{i-\half,j} \, \delta_h(\x_{i-\half,j} - \X(\s_,t)) \, h^2,\text{ and} \\
  U^{\text{IB}}_2(\s,t) &= \sum_{i,j} \left(u_2\right)_{i,j-\half} \, \delta_h(\x_{i,j-\half} - \X(\s_,t)) \, h^2.
\end{align}
Because the components of $\U^{\text{IB}}(\s,t)$ are not generally linear combinations of the FE basis functions, generally $\D{\X}{t} \neq \U^{\text{IB}}$.

For the semi-discretized partitioned formulation, $\f$ is computed on the Cartesian grid via $\f = \cS \, \F$.
The Eulerian transmission force density $\t$ is computed in a similar manner, but in this case, the numerical integration occurs only over those element boundaries that coincide with $\p U$.
We use the notation $\t = \cS^{\p U} \, \T$ to denote this operation.
We use the same regularized delta function $\delta_h(\x)$ to construct both $\cS$ and $\cS^{\p U}$.
For simplicity, we use the same velocity-restriction operator for both formulations.
This choice ensures that the two formulations coincide whenever $\T \equiv 0$.

\begin{figure}
  \centering
  \begin{tabular}{lclclc}
    \vspace{-10pt} {\bf A.} & & {\bf B.} & & {\bf C.} & \\
    & \includegraphics[width=92pt]{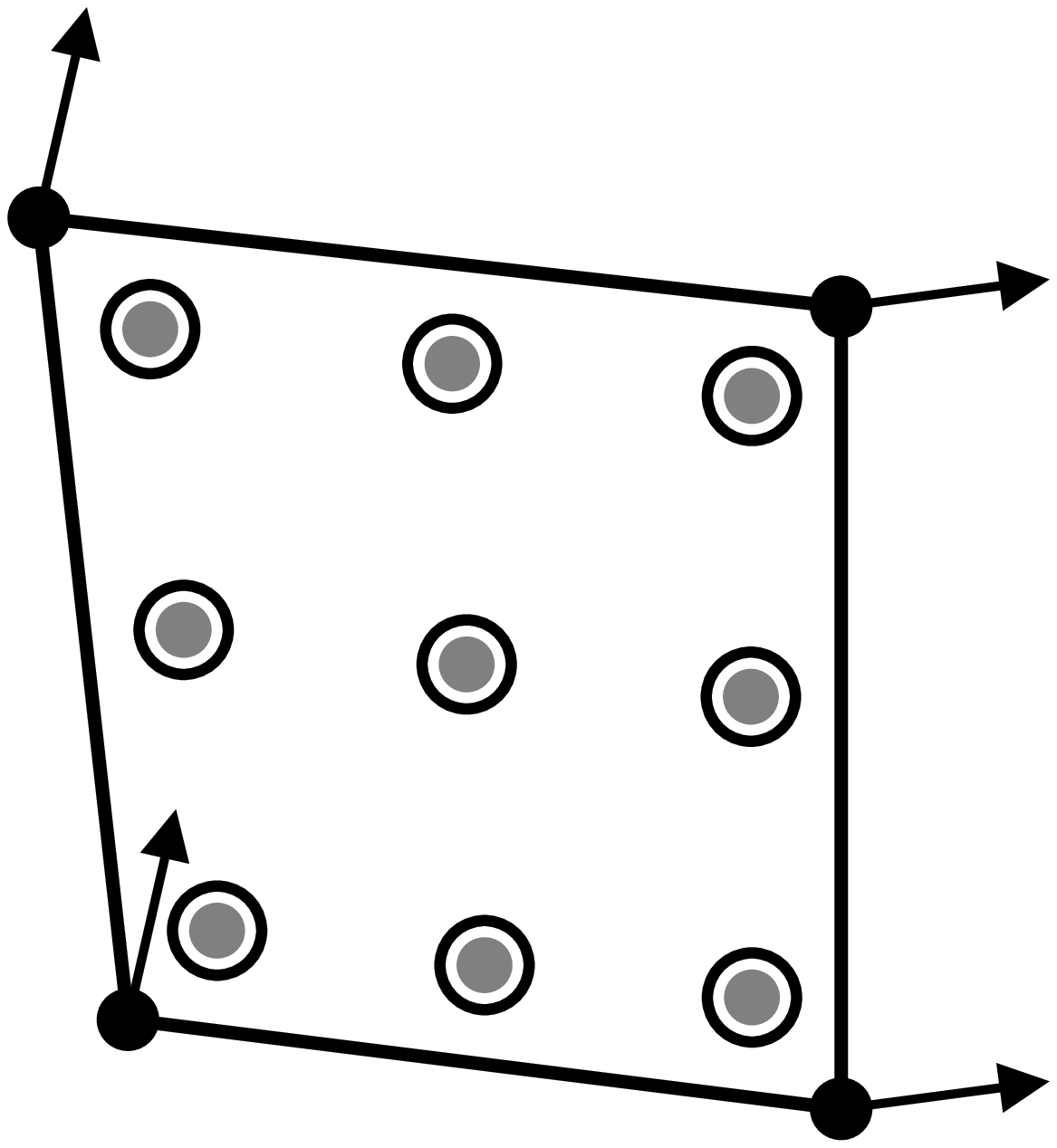} &
    & \includegraphics[width=92pt]{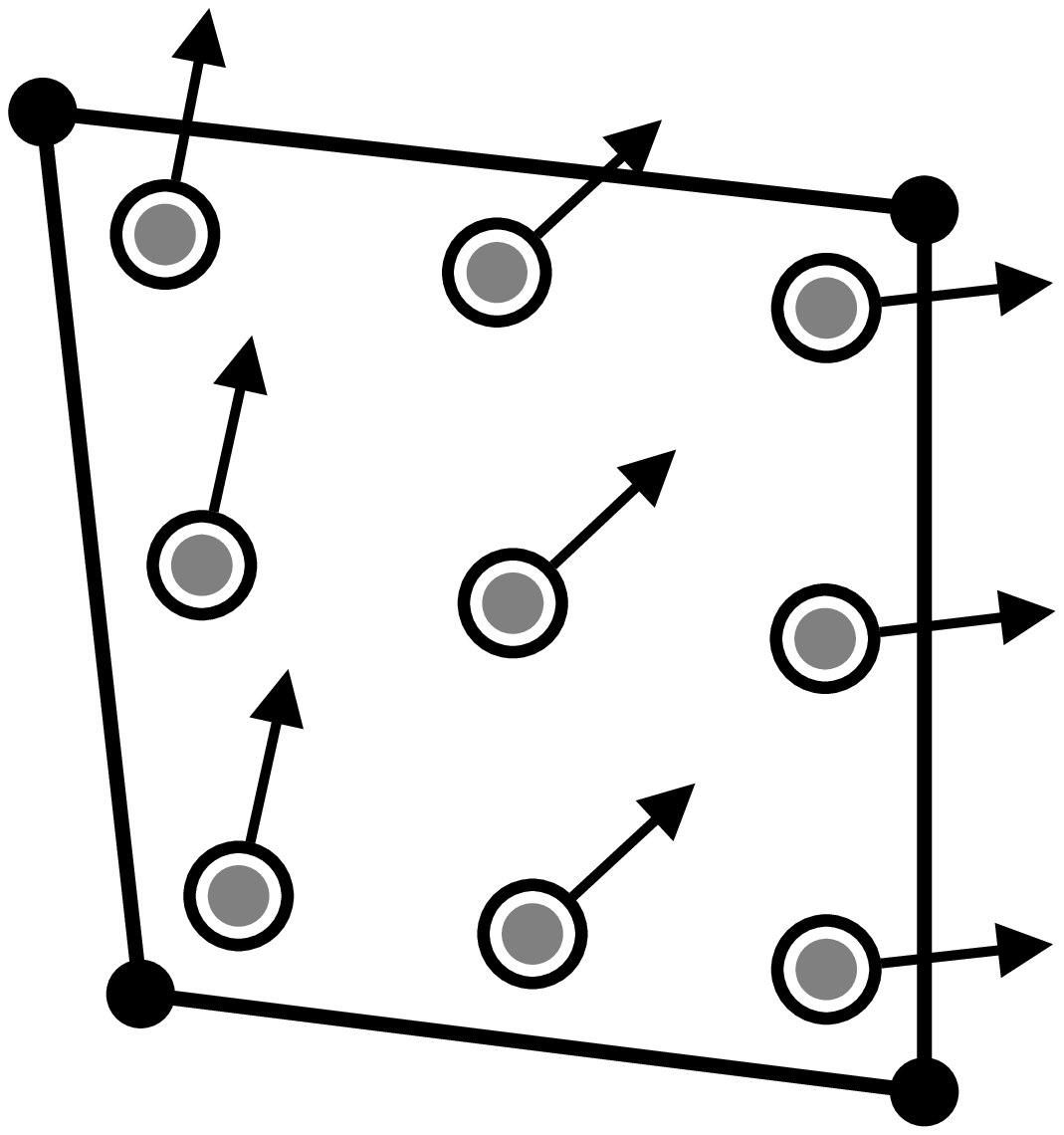} &
    & \includegraphics[width=92pt]{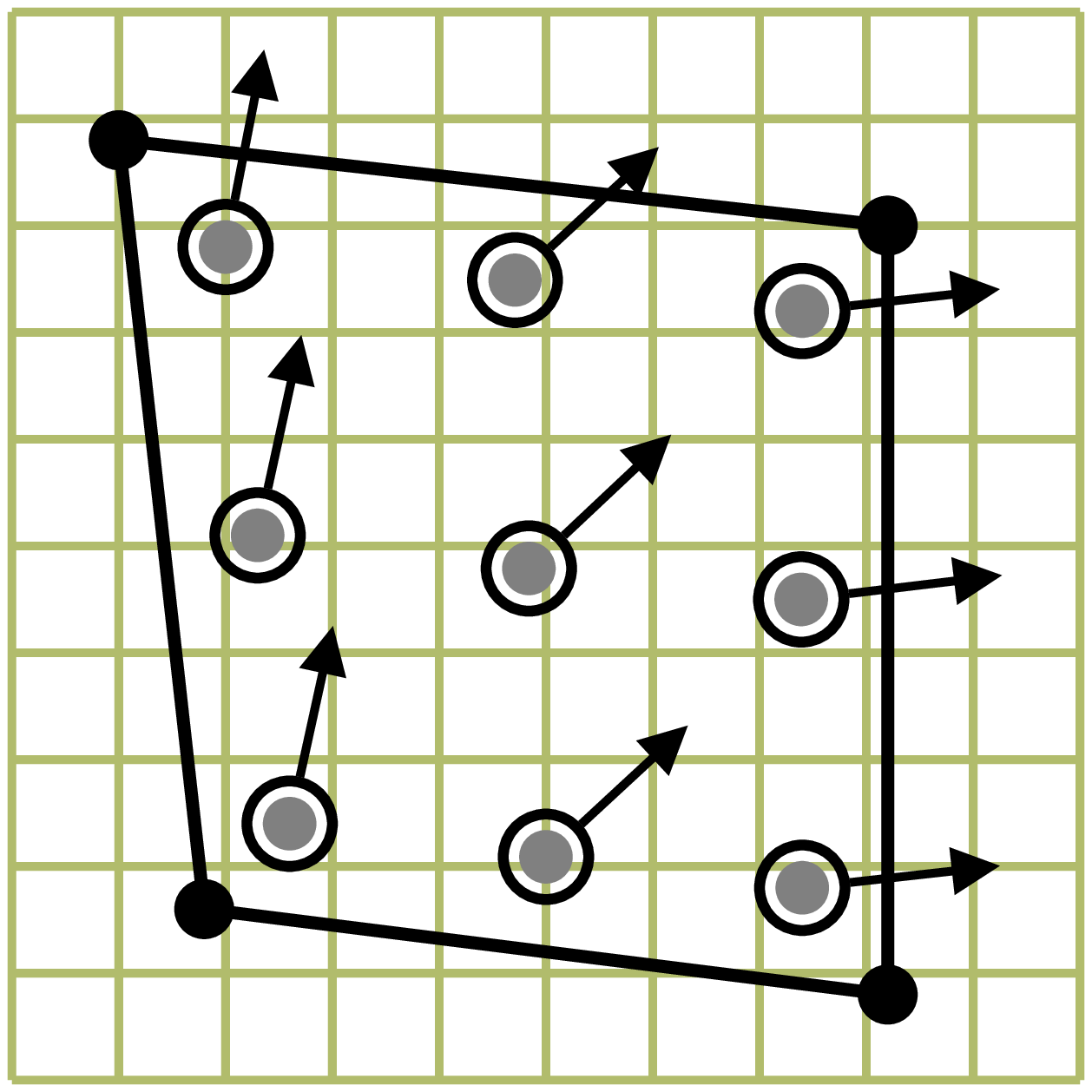}
  \end{tabular}
  \caption{
  	Prolonging the elastic force density from the Lagrangian mesh onto the Eulerian grid.
  	Starting with an approximation to the force density at the nodes of the Lagrangian mesh ({\bf A}), we use the interpolatory FE basis functions to determine the force density at interaction points that are defined by a quadrature rule ({\bf B}), and then spread those forces from the interaction points to the background Eulerian grid using the smoothed delta function $\delta_h(\x)$ ({\bf C}).
  	This approach permits Lagrangian meshes that are significantly coarser than the Eulerian grid so long as the ``net'' of interaction points is sufficiently dense.
  	Denser nets of interaction points can be obtained, for instance, by increasing the order of the numerical quadrature scheme.
  }
  \label{f:gaussian_quad}
\end{figure}

The Lagrangian-Eulerian interaction operators introduced in this work are different from analogous operators generally used in standard IB methods.
Standard IB methods and schemes such as the IFE method use regularized delta functions to apply nodal forces directly to the Cartesian grid and to interpolate Cartesian grid velocities directly to the Lagrangian nodes \cite{Peskin02}.
In such schemes, $\f(\x,t)$ would be approximated on the Eulerian grid by expressions similar to
\begin{align}
  (f_1^{\text{IB}})_{i-\half,j}
  &= \sum_{l=1}^{M} \left(F_1\right)_l(t) \, \delta_h(\x_{i-\half,j} - \X_l(t)) \, \omega_l^{\text{IB}},\text{ and} \label{e:IB_force_prolongation1} \\
  (f_2^{\text{IB}})_{i,j-\half}
  &= \sum_{l=1}^{M} \left(F_2\right)_l(t) \, \delta_h(\x_{i,j-\half} - \X_l(t)) \, \omega_l^{\text{IB}}, \label{e:IB_force_prolongation2}
\end{align}
in which $\f^{\text{IB}}$ denotes the Eulerian force determined by the standard IB force-spreading operator and $\omega_l^{\text{IB}}$ denotes the volume associated with Lagrangian node $l$.
In this approach, each nodal force $\F_l$ is applied only to Cartesian grid cells within the support of $\delta_h(\x-\X_l)$, and the Lagrangian mesh must therefore be finer than the Cartesian grid to avoid leaks \cite{Peskin02}.
The corresponding approach to velocity restriction used by such methods would be to set $\d{\X_l}{t} = \U^{\text{IB}}(\s_l,t)$.

Our Lagrangian-Eulerian interaction operators 
communicate between a collection of Lagrangian control points (the nodes of the structural mesh) and the Cartesian grid via a collection of interaction points (the Lagrangian quadrature points).
The force-prolongation operator can be seen as the composition of two operations: first, the Lagrangian force density is evaluated at the interaction points in terms of data defined at the control points; then, the standard IB delta function $\delta_h(\x)$ spreads volume-~or area-weighted force densities from the interaction points to the Cartesian grid.
See fig.~\ref{f:gaussian_quad}.
Velocity restriction is similar:
First, the Cartesian velocity field is interpolated to the interaction points using $\delta_h(\x)$; then, these velocities are accumulated to form the right-hand-side of a system of equations that determines $\D{\X}{t}$ at the control points.
Our approach is similar to methods used in RBF-based IB methods \cite{VShankar15-IB-RBF}.

In general, it is necessary that the same interaction points are used in the discrete force-spreading and velocity-interpolation operators if those operators are to satisfy a discrete adjoint property.
It is possible to construct an interpolation operator that uses the control points as the interaction points, but then satisfying the discrete adjoint property requires that the control points are also used as the interaction points in the spreading operator.

Our numerical tests indicate that in our scheme, 
the Lagrangian structure appears watertight to the fluid so long as the net of interaction points is sufficiently dense.
Denser nets of interaction points can be obtained by increasing the order of the quadrature rule, and this may be done in an adaptive manner as the simulation progresses.
In our numerical tests, we use \emph{dynamically adapted} Gaussian quadrature rules to construct $\cS$ and $\cJ$ that provide, on average, at least a $3 \times 3$ net of quadrature points per Cartesian grid cell.

\section{Implementation}
\label{s:implementation}

This version of the IB method is implemented in the open-source IBAMR software \cite{IBAMR-web-page}, which is a C++ framework for developing FSI models that use the IB method.
IBAMR provides support for distributed-memory parallelism and adaptive mesh refinement (AMR).
IBAMR relies upon the SAMRAI \cite{samrai-web-page,HornungKohn02,HornungWissinkKohn06}, PETSc \cite{petsc-web-page,petsc-user-ref,petsc-efficient}, \emph{hypre} \cite{hypre-web-page,FalgoutYang02}, and \texttt{libMesh} \cite{libMesh-web-page,libMeshPaper} libraries for much of its functionality.

\section{Numerical results}
\label{s:numerical_results}

\subsection{Thick elliptical shell}
\label{s:thick_elliptical_shell}

\begin{figure}
  \centering
  \tabcolsep 2.5pt
  \begin{tabular}{lrrr}
    \vspace{-10pt} {\bf A.} & & & \\
    & \includegraphics[height=100pt]{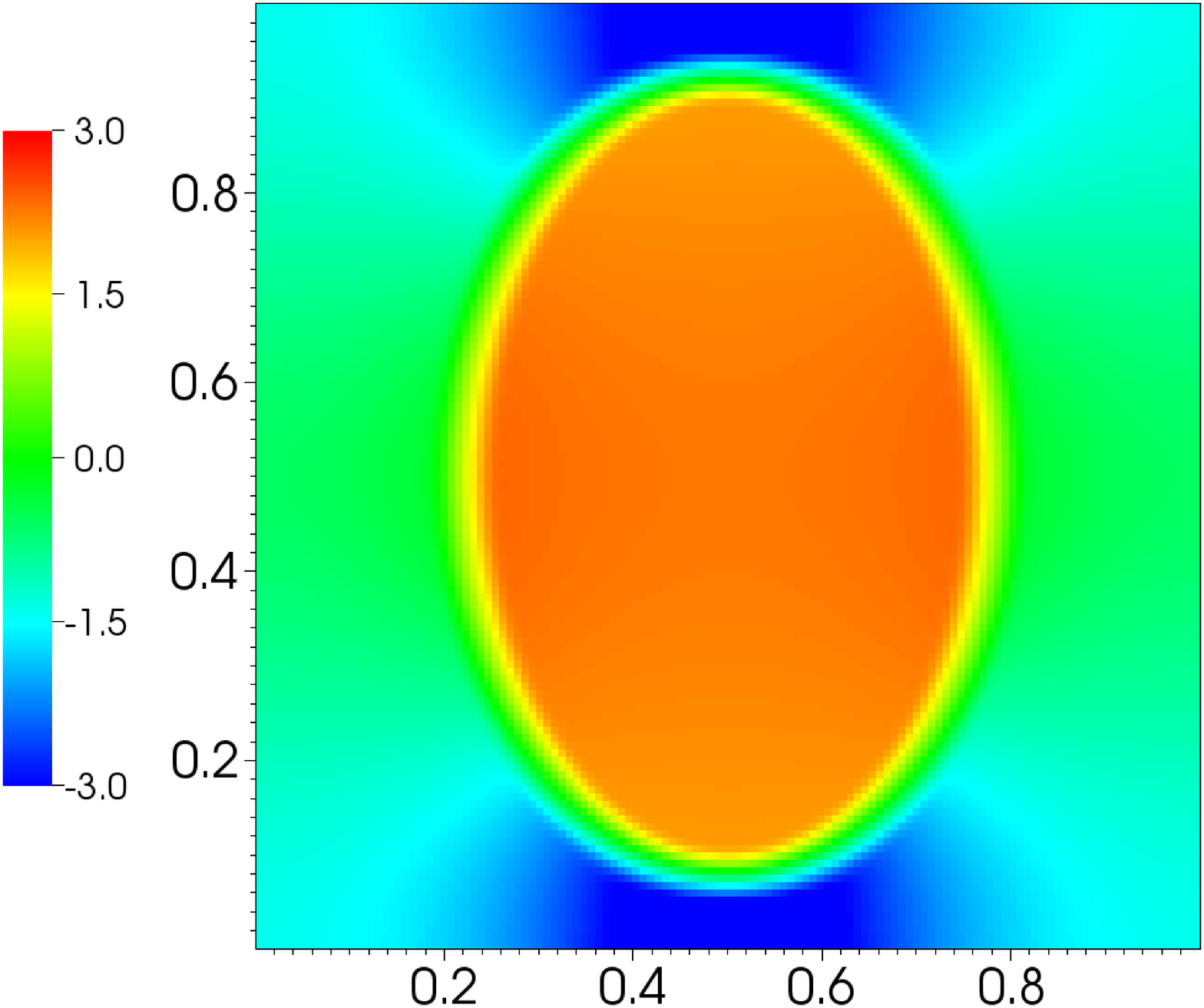}
    & \includegraphics[height=100pt]{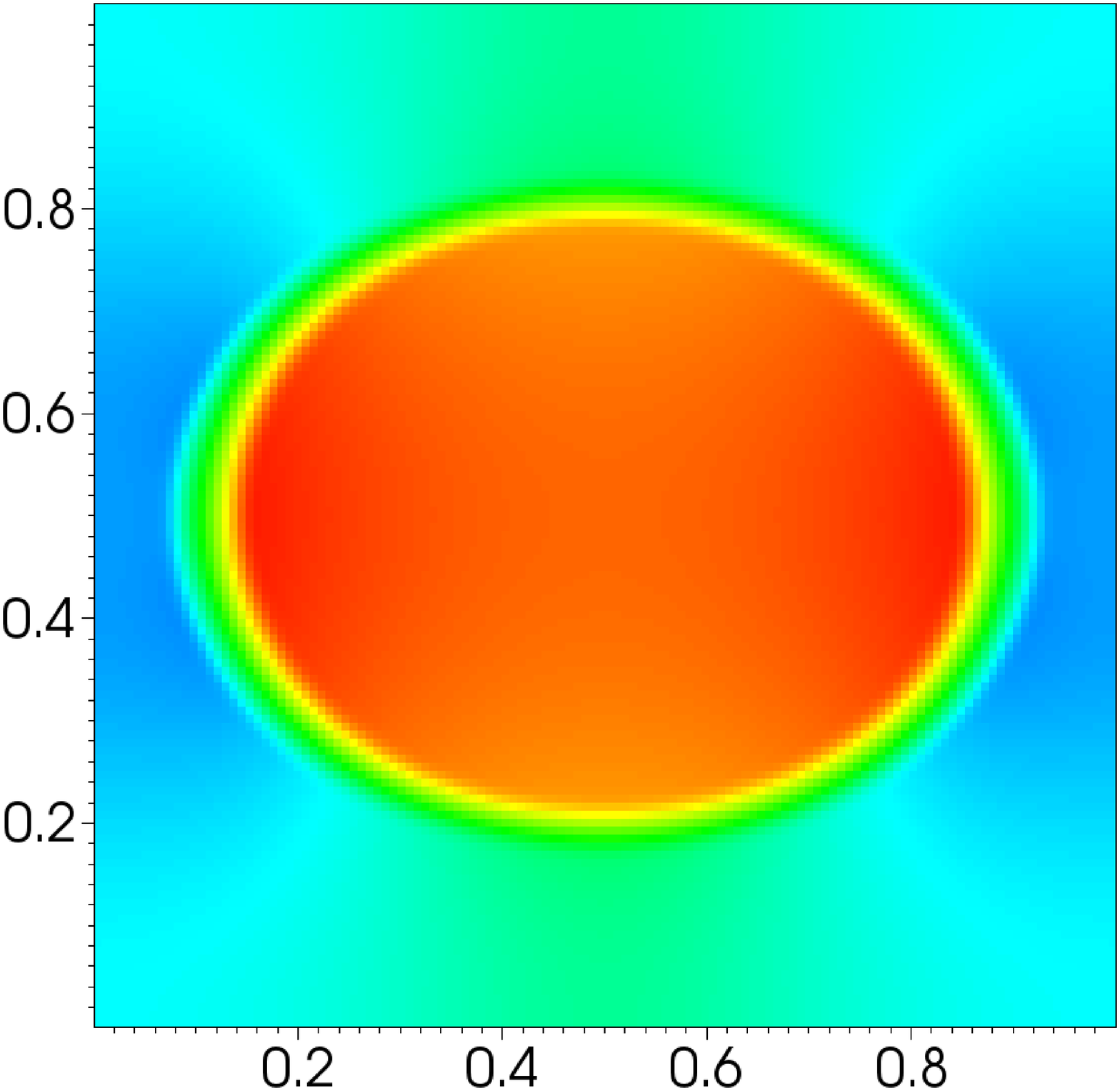}
    & \includegraphics[height=100pt]{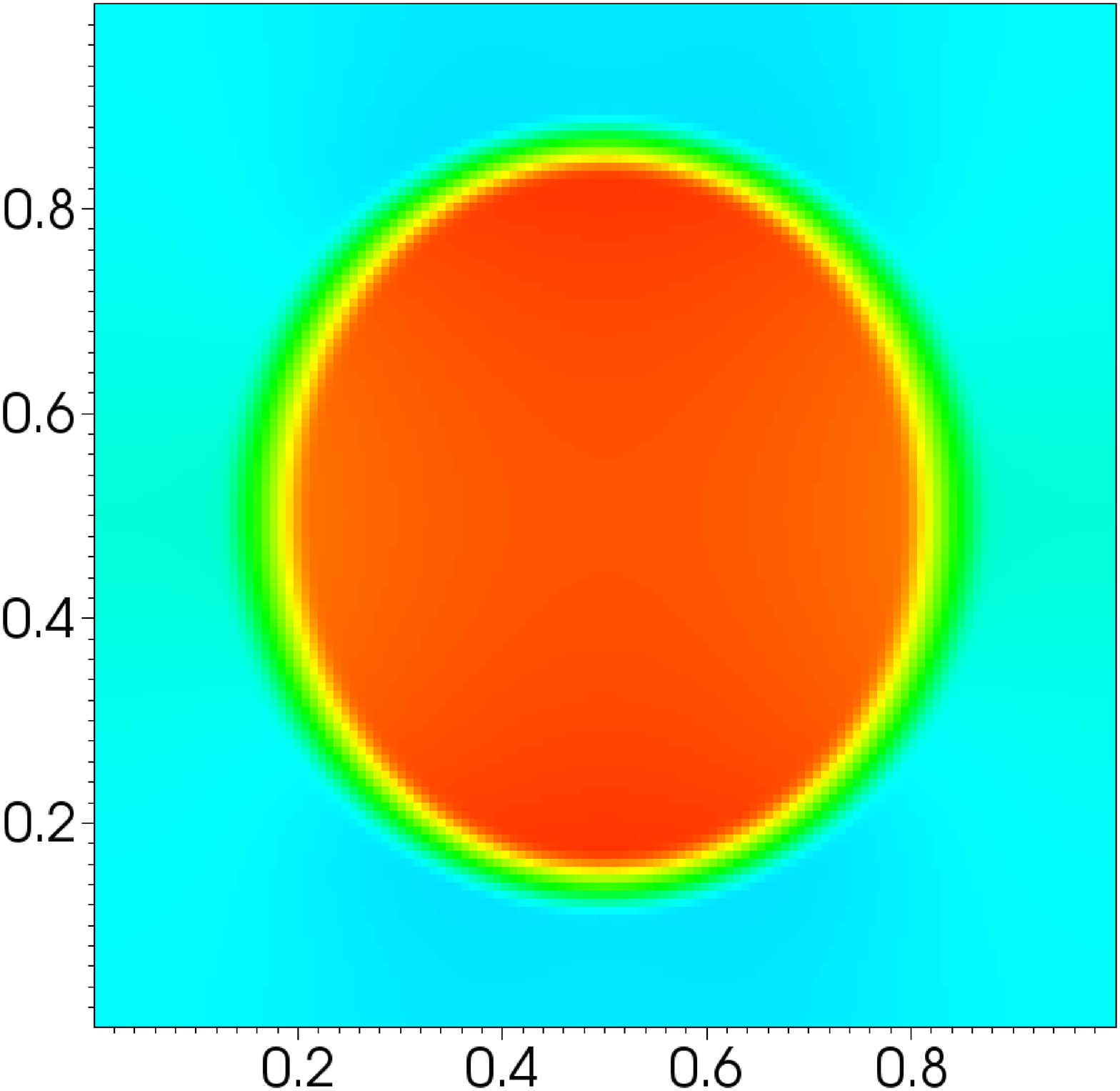} \\
    \vspace{-10pt} {\bf B.} & & & \\
    & \includegraphics[height=100pt]{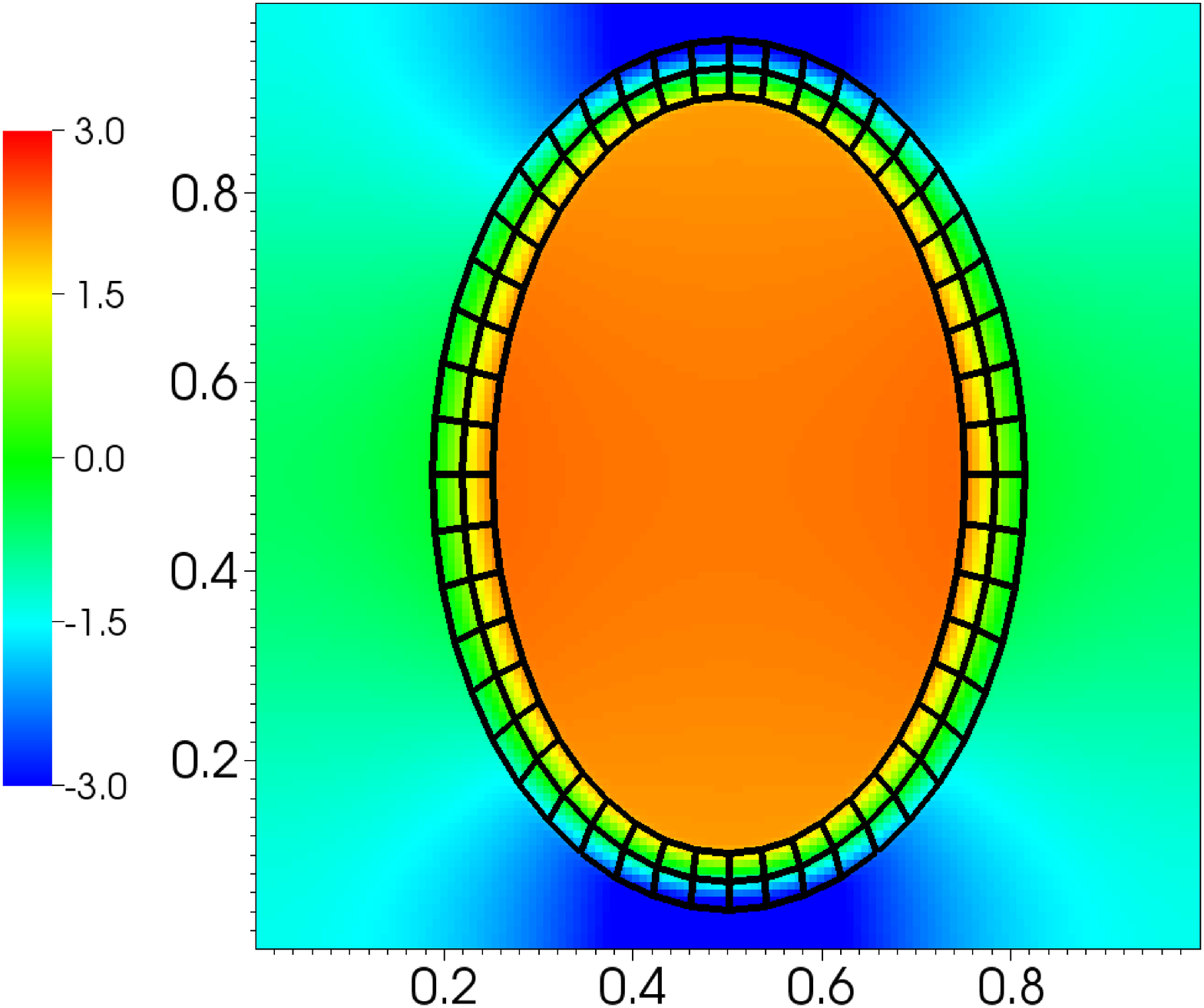}
    & \includegraphics[height=100pt]{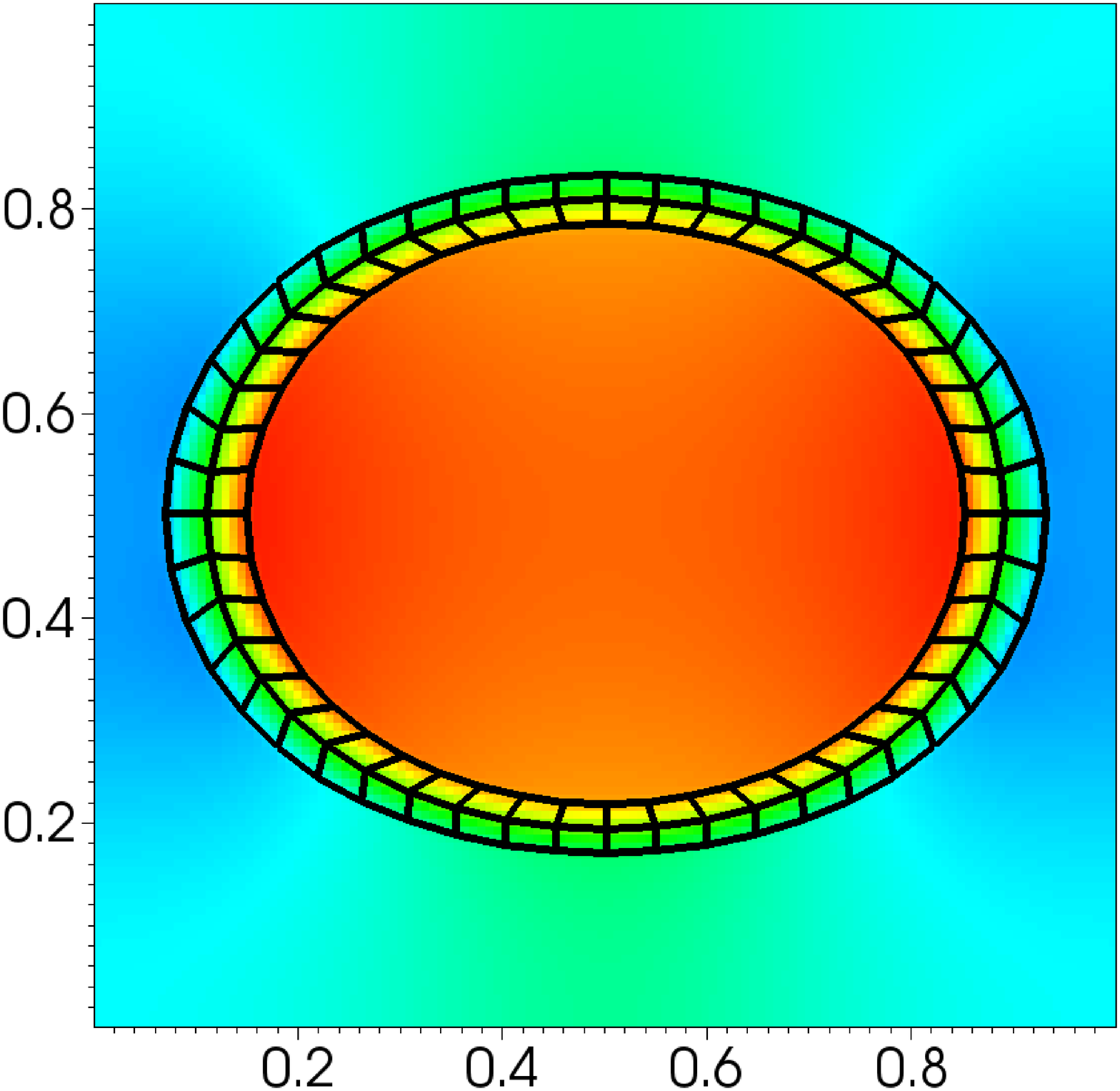}
    & \includegraphics[height=100pt]{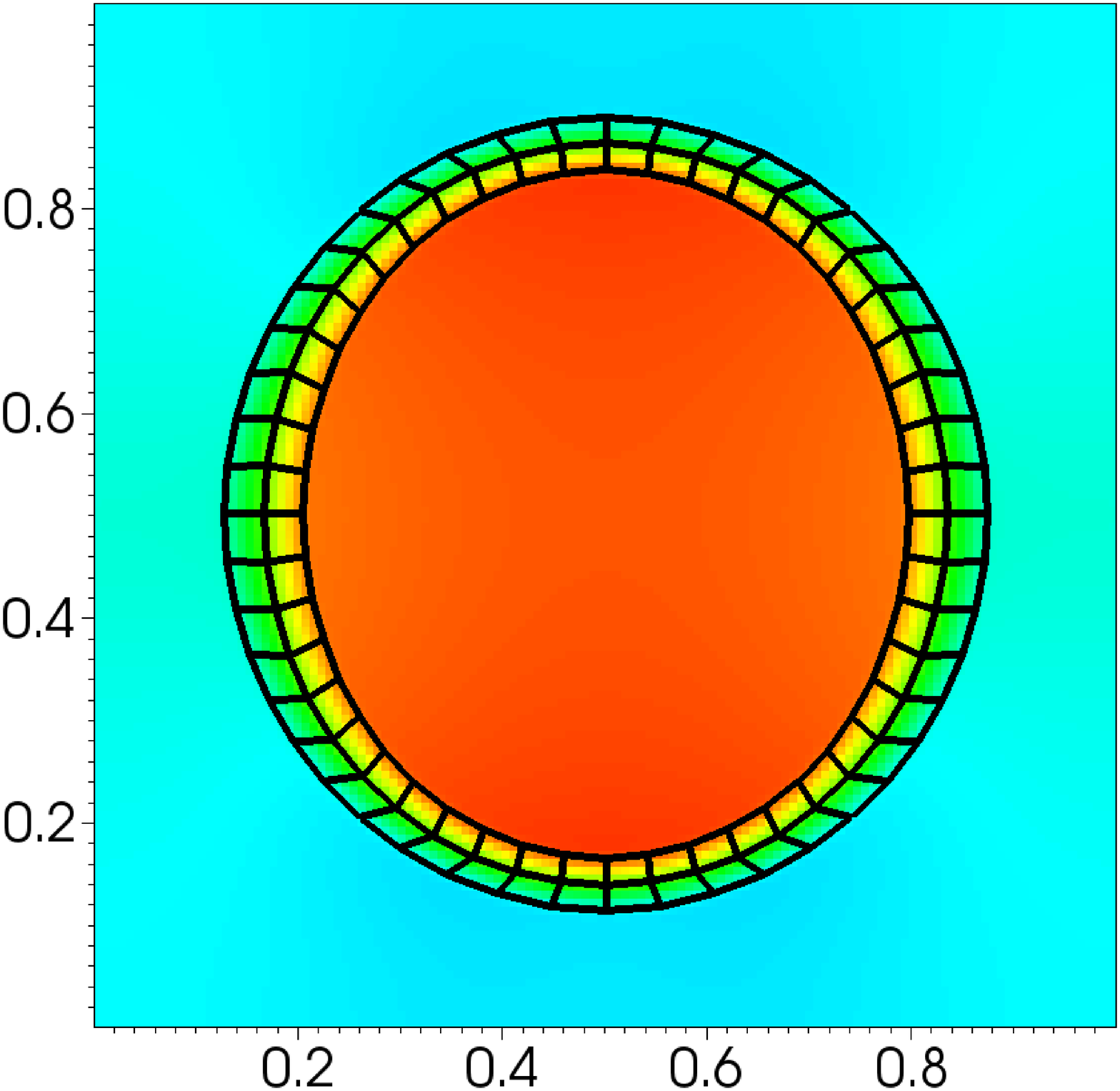} \\
    \vspace{-10pt} {\bf C.} & & & \\
    & \includegraphics[height=100pt]{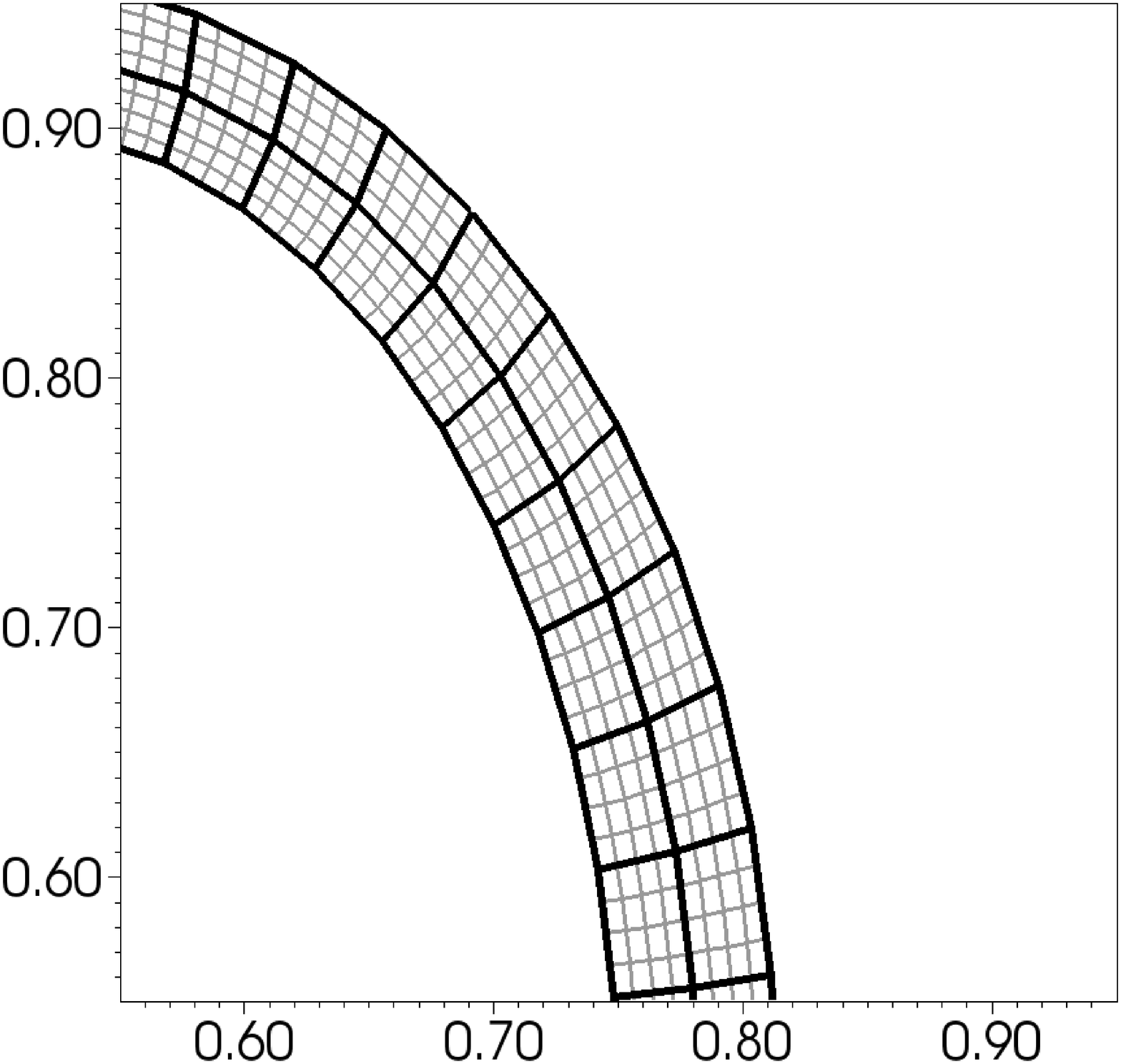}
    & \includegraphics[height=100pt]{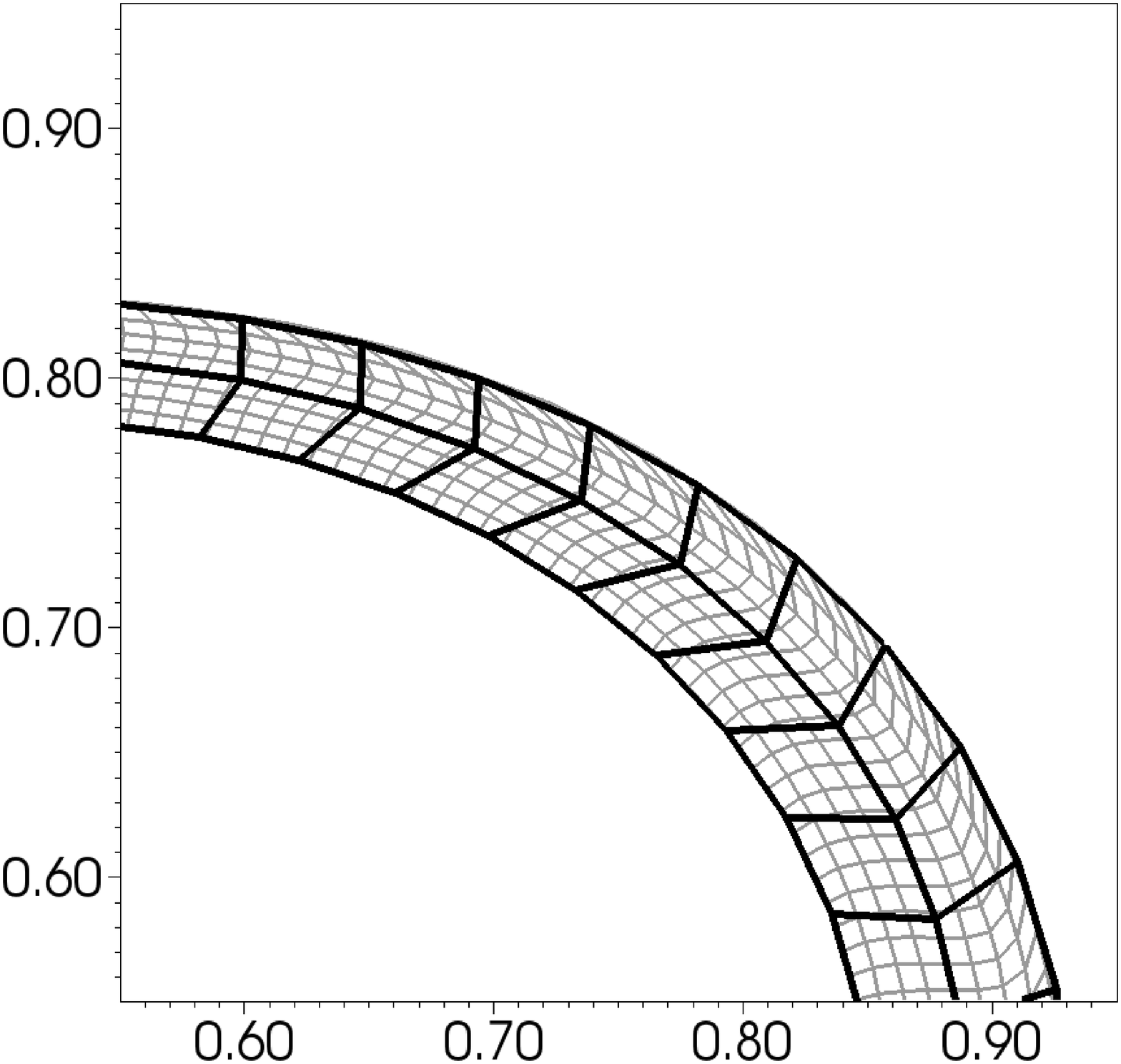}
    & \includegraphics[height=100pt]{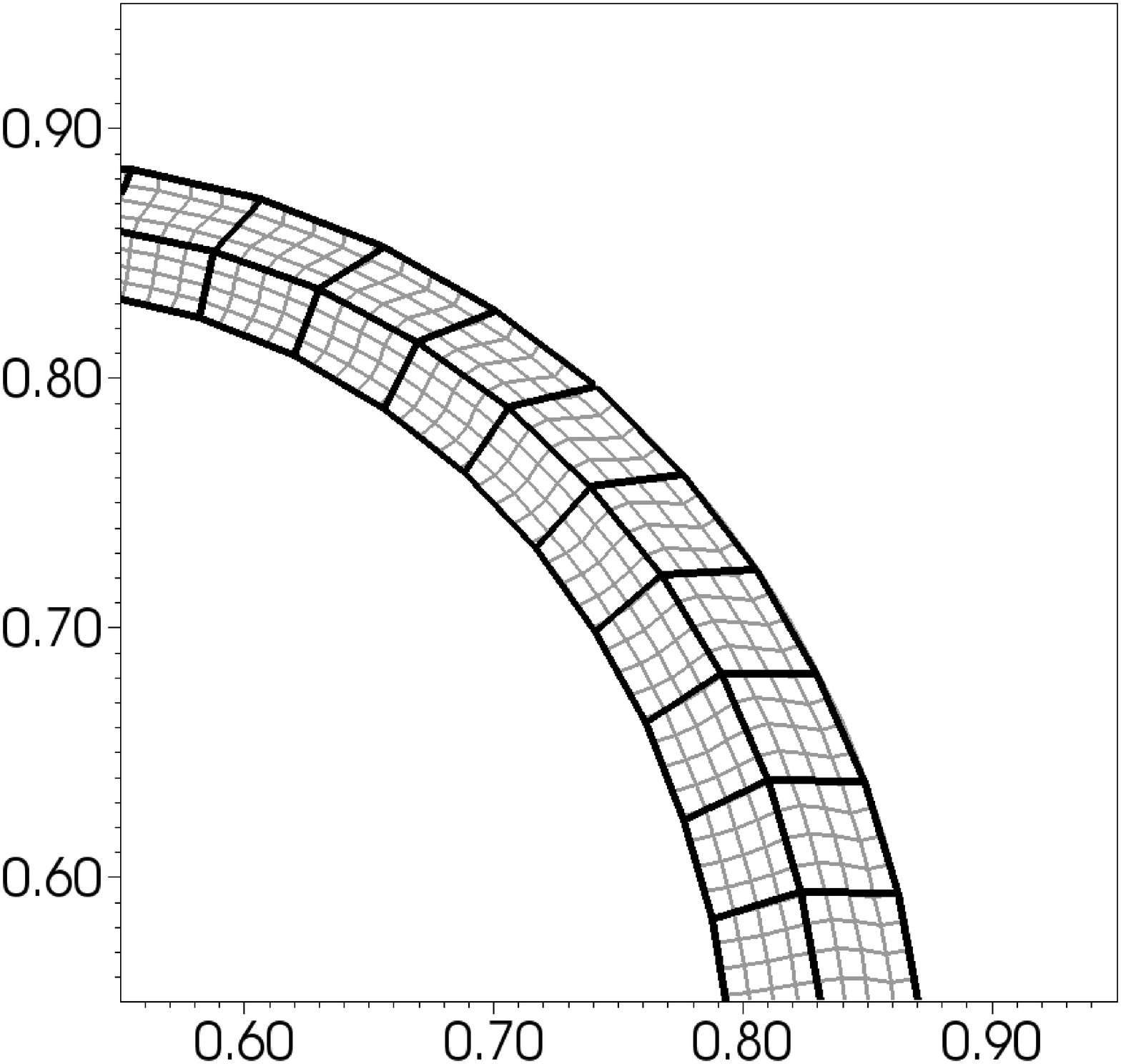} \\
  \end{tabular}
  \caption{
  	Representative results from the dynamic ($\gamma=0.15$) version of the idealized anisotropic shell model of sec.~\ref{s:anisotropic_shell} for $N=128$ over the time interval $0 \leq t \leq 0.75$.
  	The computed pressure and structure deformation for $\Mfac = 4$ are shown in {\bf A} and {\bf B}.
  	The computed deformations obtained with $\Mfac = 1$ and $\Mfac = 4$ are compared in {\bf C}.
  	The coarse and fine structural meshes yield essentially identical kinematics.
  }
  \label{f:anisotropic_dynamic}
\end{figure}

\begin{figure}
  \centering
  \tabcolsep 2.5pt
  \begin{tabular}{lrrr}
    \vspace{-10pt} {\bf A.} & & & \\
    & \includegraphics[height=100pt]{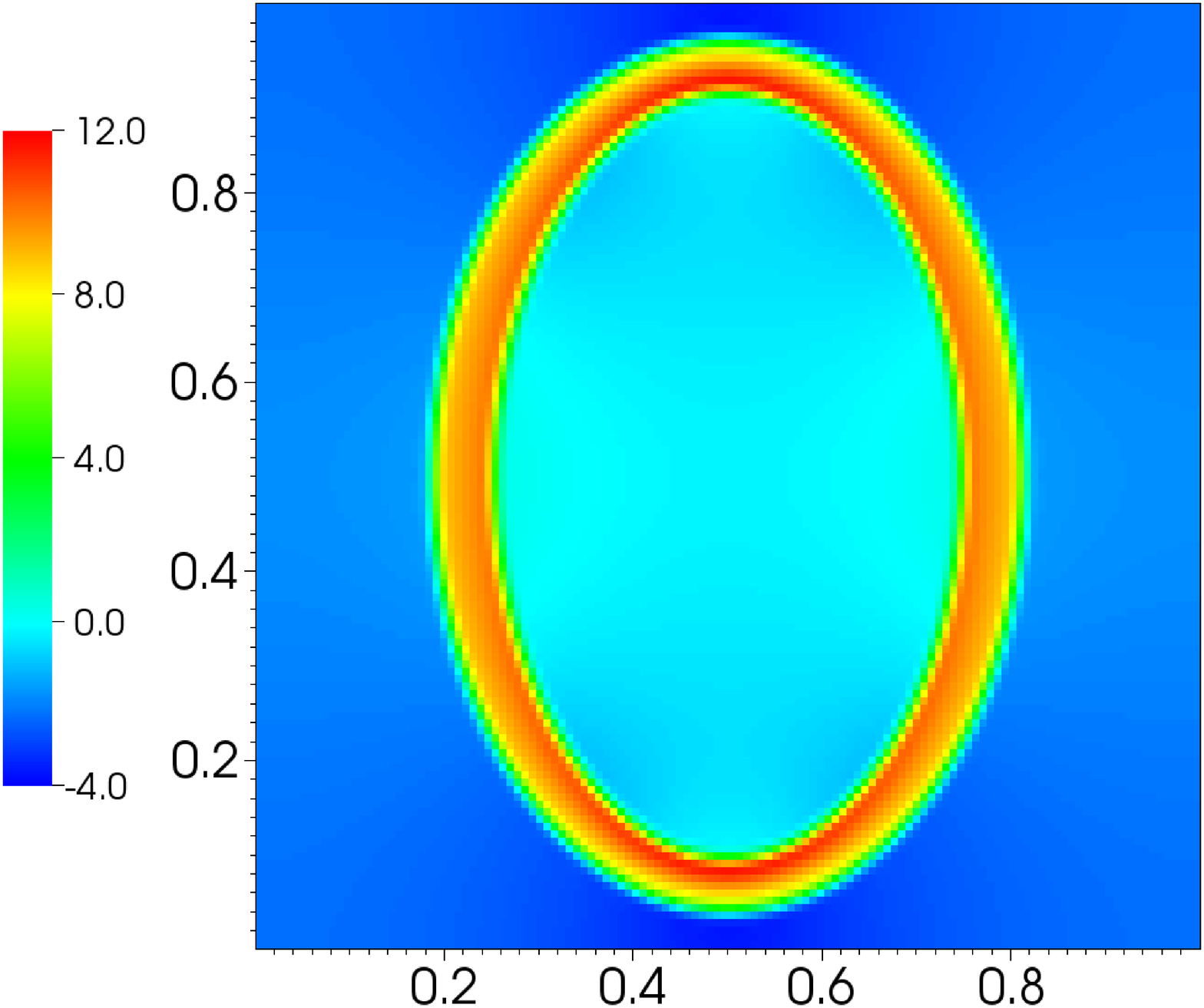}
    & \includegraphics[height=100pt]{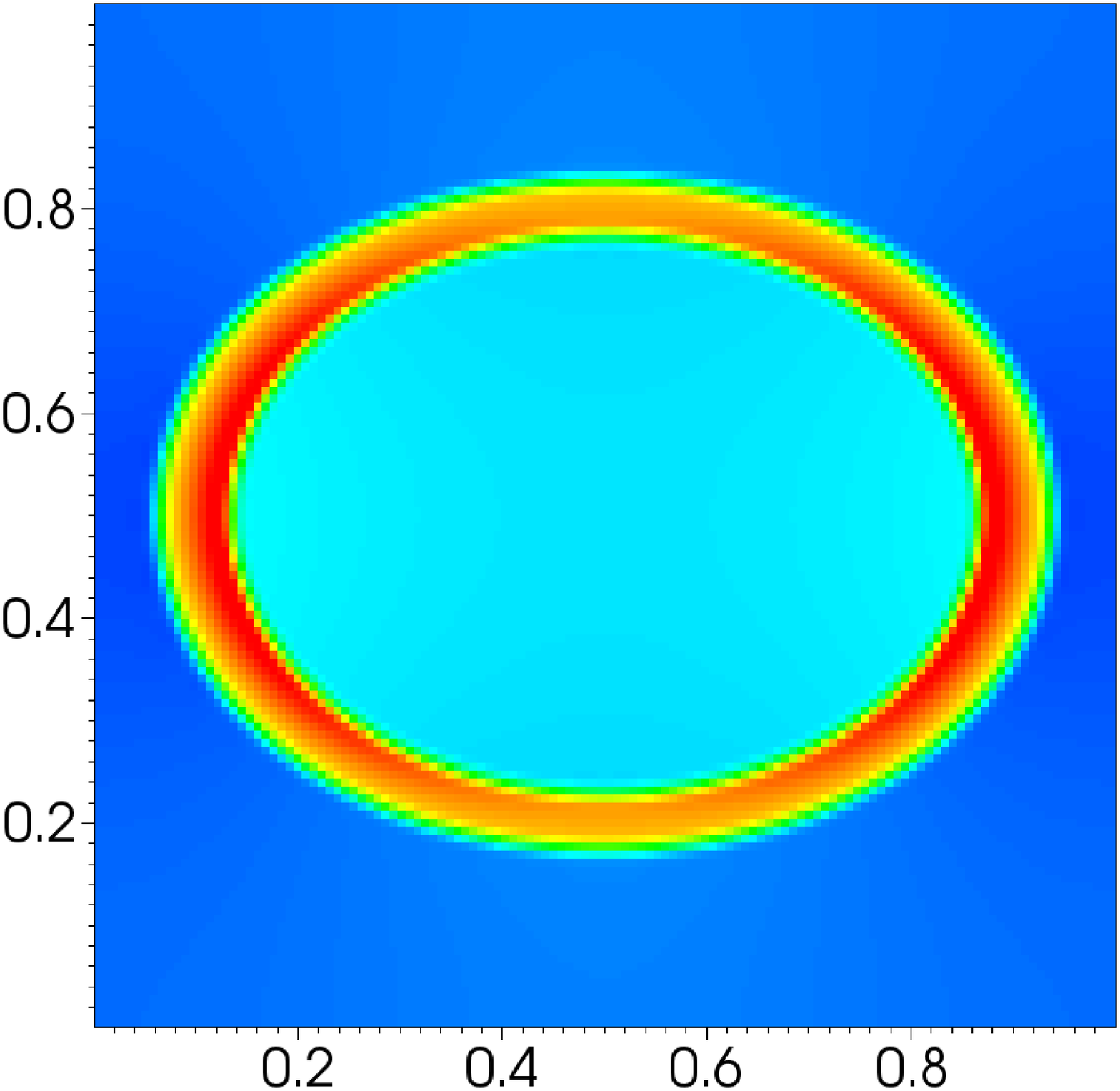}
    & \includegraphics[height=100pt]{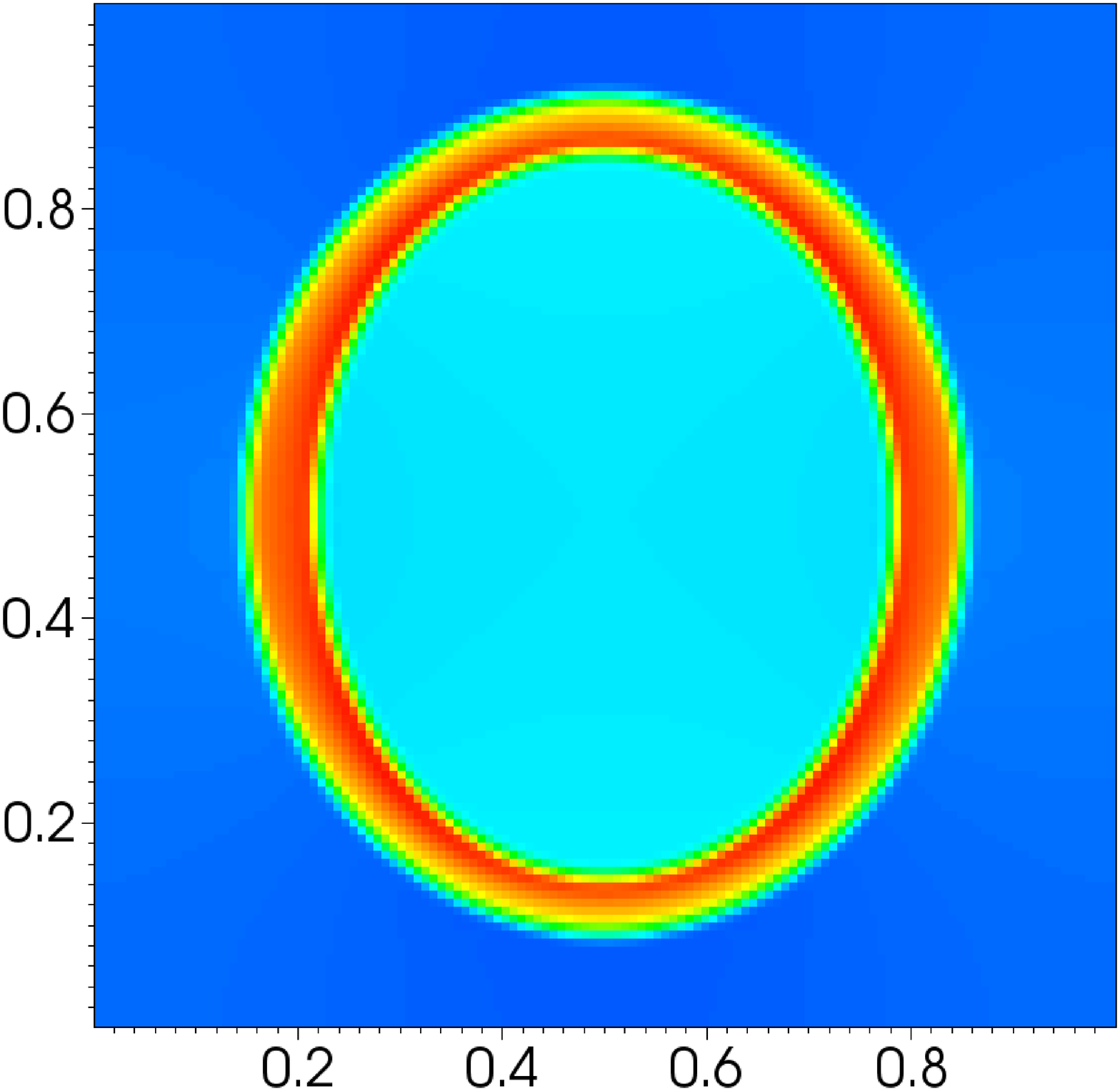} \\
    \vspace{-10pt} {\bf B.} & & & \\
    & \includegraphics[height=100pt]{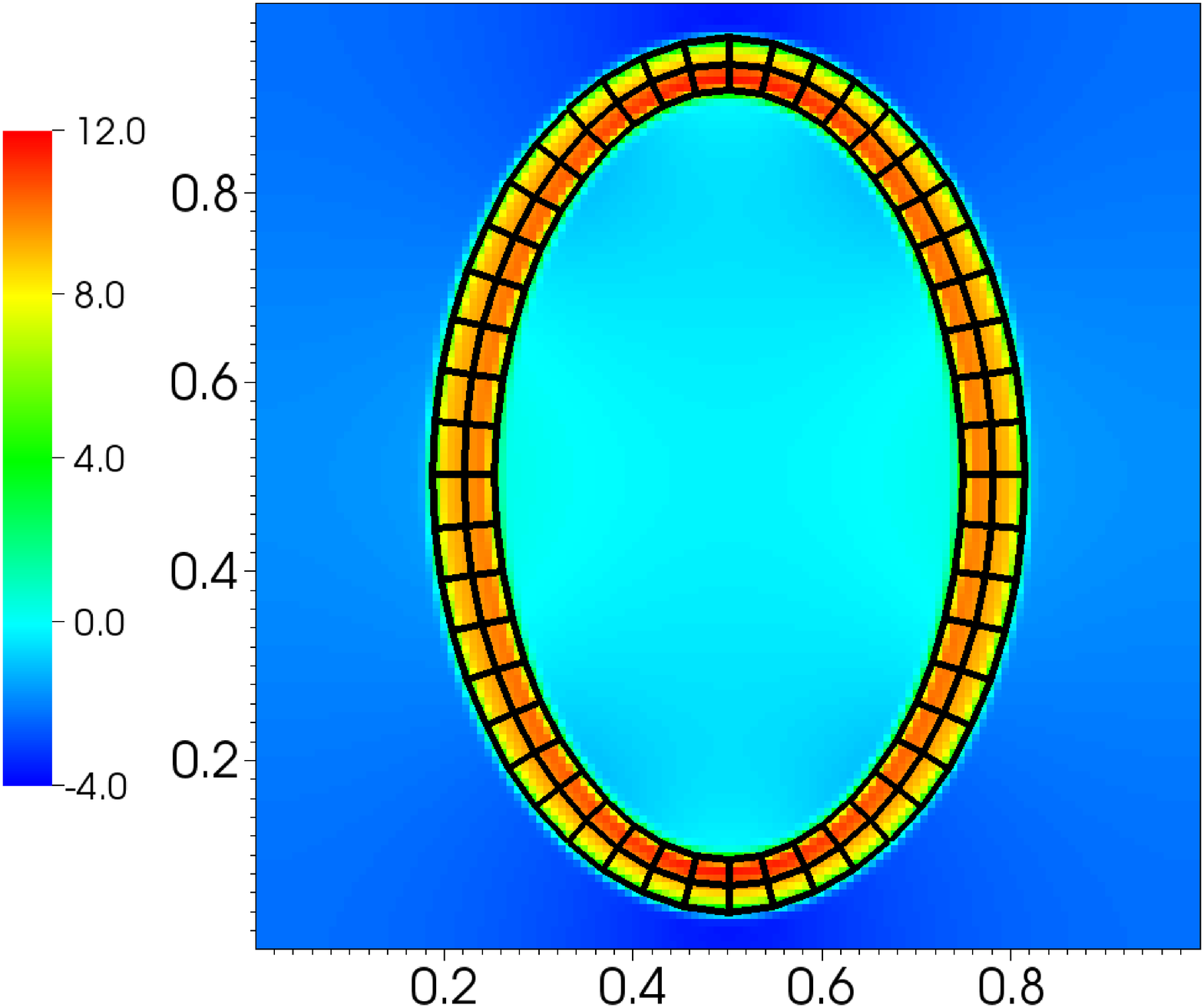}
    & \includegraphics[height=100pt]{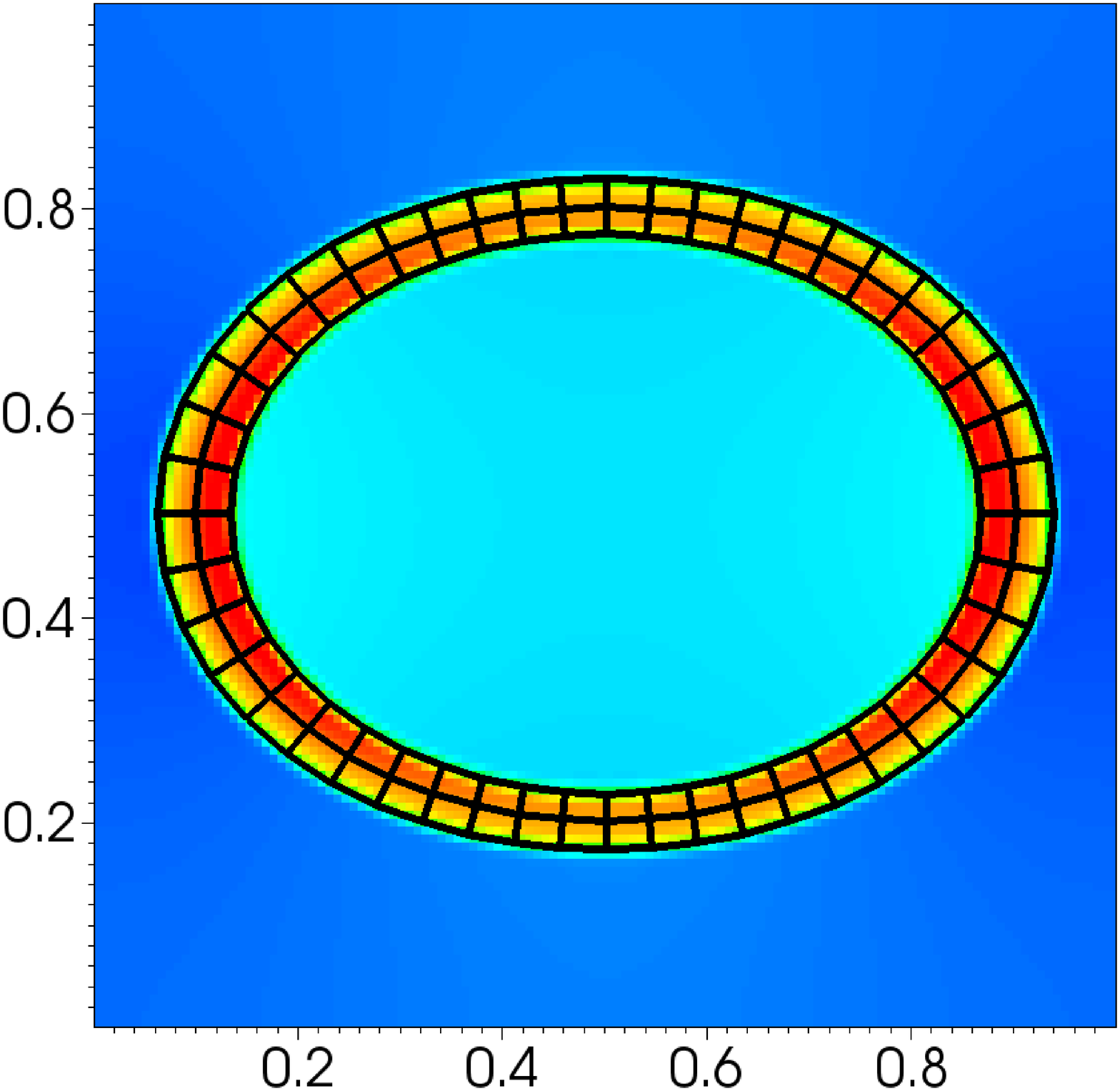}
    & \includegraphics[height=100pt]{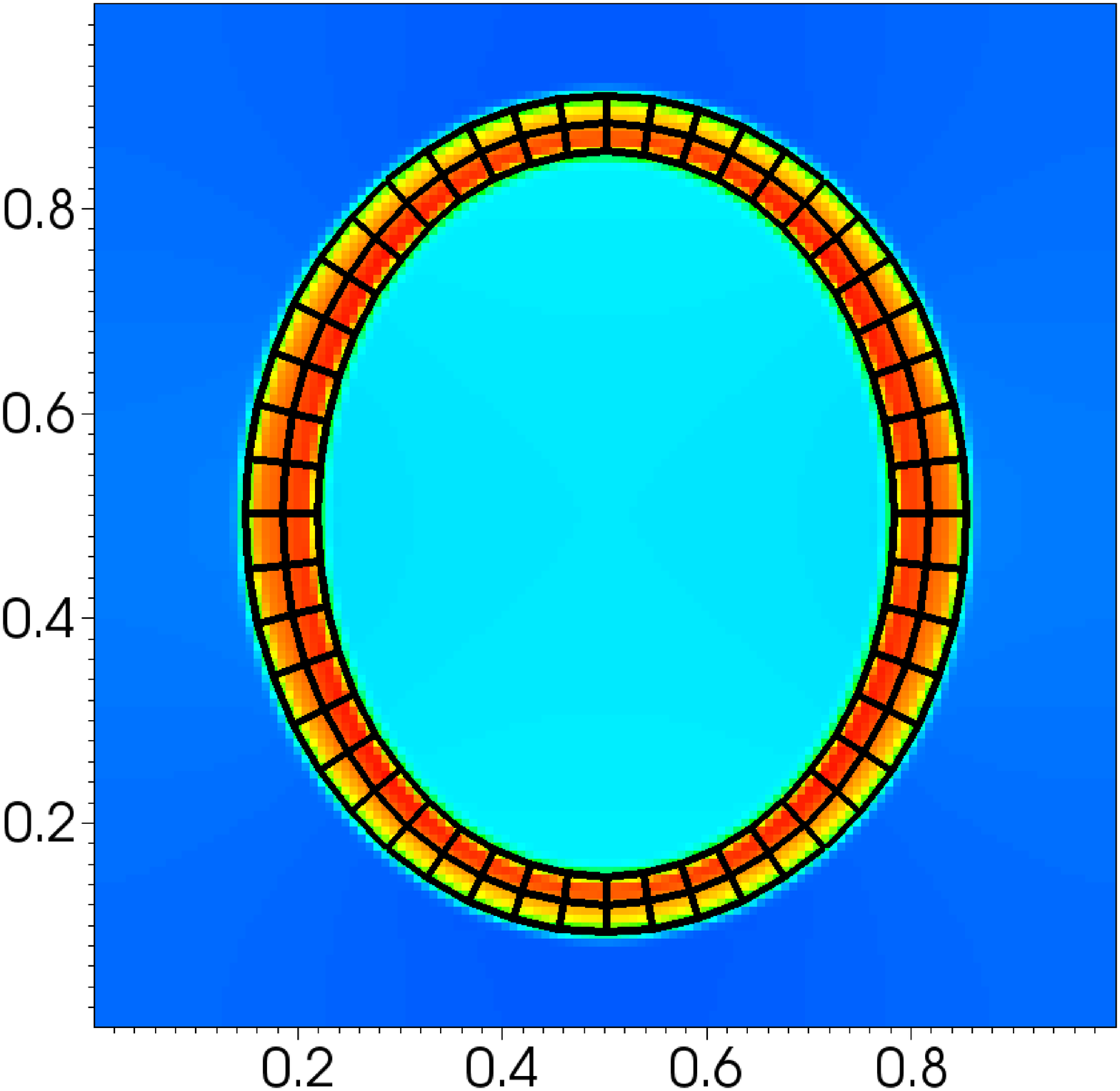} \\
    \vspace{-10pt} {\bf C.} & & & \\
    & \includegraphics[height=100pt]{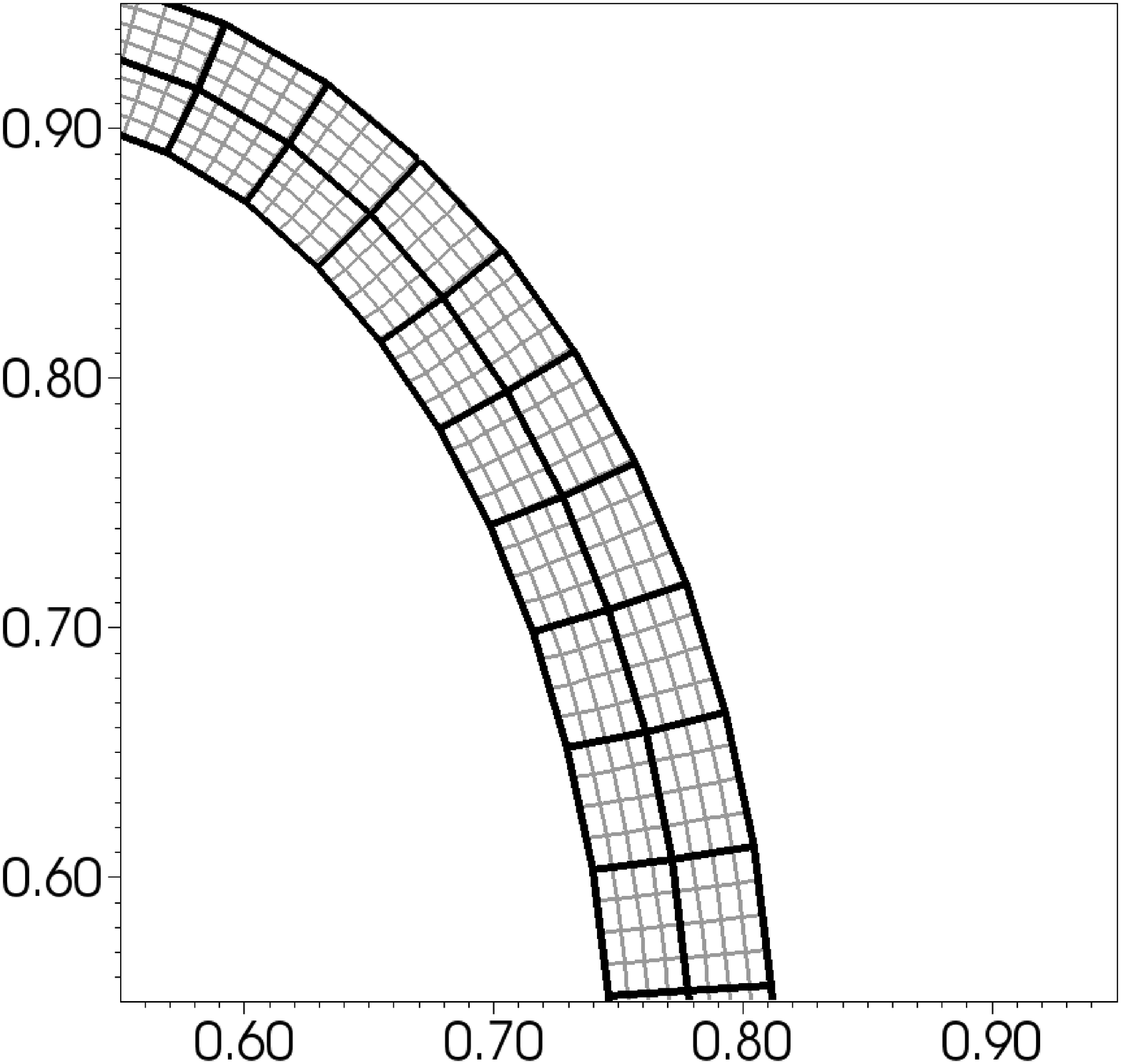}
    & \includegraphics[height=100pt]{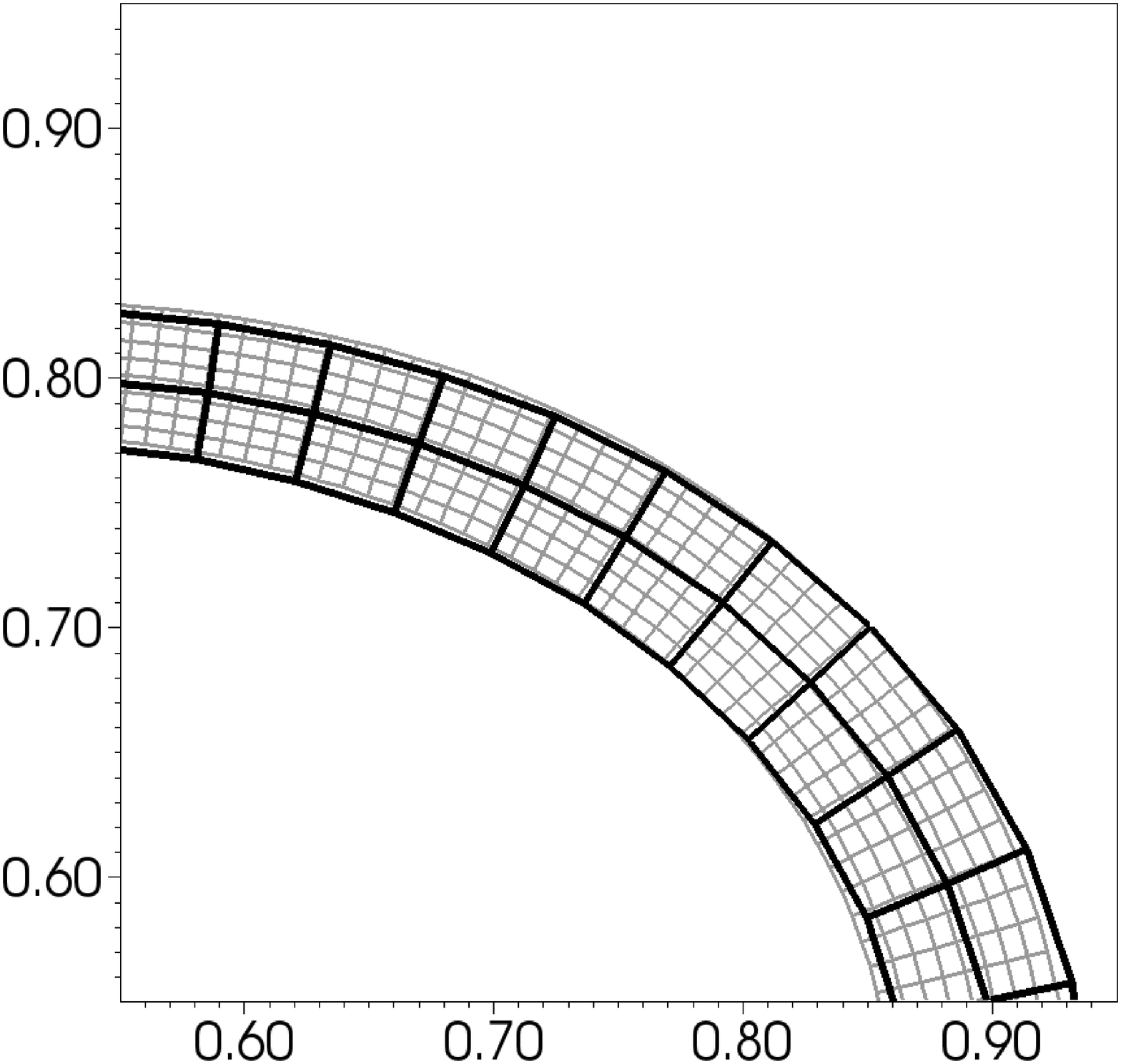}
    & \includegraphics[height=100pt]{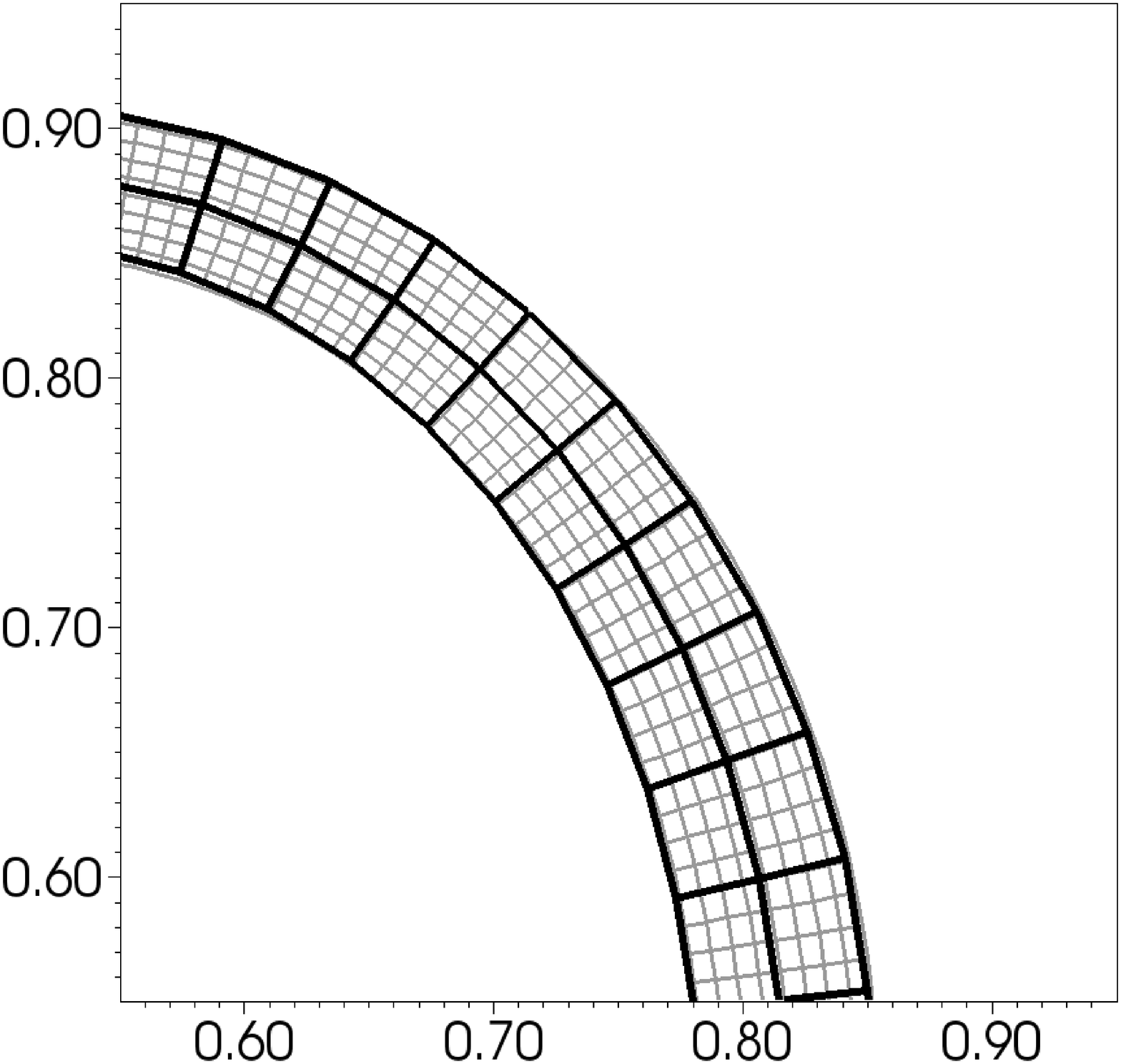} \\
  \end{tabular}
  \caption{
  	Similar to fig.~\ref{f:anisotropic_dynamic}, but here showing results obtained using the orthotropic shell model of sec.~\ref{s:orthotropic_shell} for $N=128$ and the partitioned (split) weak formulation over the time interval $0 \leq t \leq 1.25$.
  	As in fig.~\ref{f:anisotropic_dynamic}, the coarse and fine structural meshes yield very similar kinematics.
  }
  \label{f:orthotropic_dynamic}
\end{figure}

This section presents results from tests that use a thick elastic shell \cite{BEGriffith05-ib_accuracy,BEGriffith07-ibamr_paper,DBoffi08} to demonstrate the convergence properties of our method for different types of material models.
In these computations, the physical domain is $\Omega = [0,1] \times [0,1]$ with periodic boundary conditions, and, following Boffi et al.~\cite{DBoffi08}, the Lagrangian coordinate domain $U$ is defined using curvilinear coordinates $\vec{s} = (s_1,s_2) \in U$ instead of reference coordinates, with $U = [0 , 2\pi R] \times [0 , w]$ for $R = 0.25$ and $w = 0.0625$, and with periodic boundary conditions in the $s_1$ direction.
The coordinate mapping $\X:(U,t)\mapsto\Omega$ is given at time $t=0$ by
\begin{equation}
  \X(\vec{s},0) = \left( \cos(s_1/R)(R+s_2)+0.5, \sin(s_1/R)(R+\gamma+s_2)+0.5 \right).
\end{equation}
For $\gamma = 0$, the initial configuration of the shell is a circular annulus with inner radius $R$ and thickness $w$, which corresponds to an equilibrium configuration of the structure.
For $\gamma \neq 0$, the initial configuration is an elliptical annulus that is out of equilibrium.
In our tests, we use $\gamma=0$ for static problems and $\gamma=0.15$ for dynamic problems.
In either case, we discretize $\Omega$ using an $N$-by-$N$ Cartesian grid.
The Lagrangian discretization is constructed so that the nodes of the Lagrangian mesh are physically separated by a distance of approximately $\Mfac \dx$.
Specifically, we discretize $U$ using a $28 M$-by-$M$ mesh of bilinear ($Q^1$) elements, with $M = \frac{1}{\Mfac} \frac{N}{16}$.
Representative numerical results using $N=128$ are shown in figs.~\ref{f:anisotropic_dynamic} and \ref{f:orthotropic_dynamic}.

Although this is a relatively simple benchmark problem, the static version ($\gamma=0$) is one of the only test problems available for the IB method that we know of that permits a simple analytic solution.
Moreover, because certain choices of $\PPe$ allow the IB method to attain higher-order convergence rates, this test case allows us to verify that our implementation attains its designed order of accuracy.

\subsubsection{Anisotropic shell}
\label{s:anisotropic_shell}

\begin{figure}
  \centering
  \includegraphics[width=0.775\textwidth]{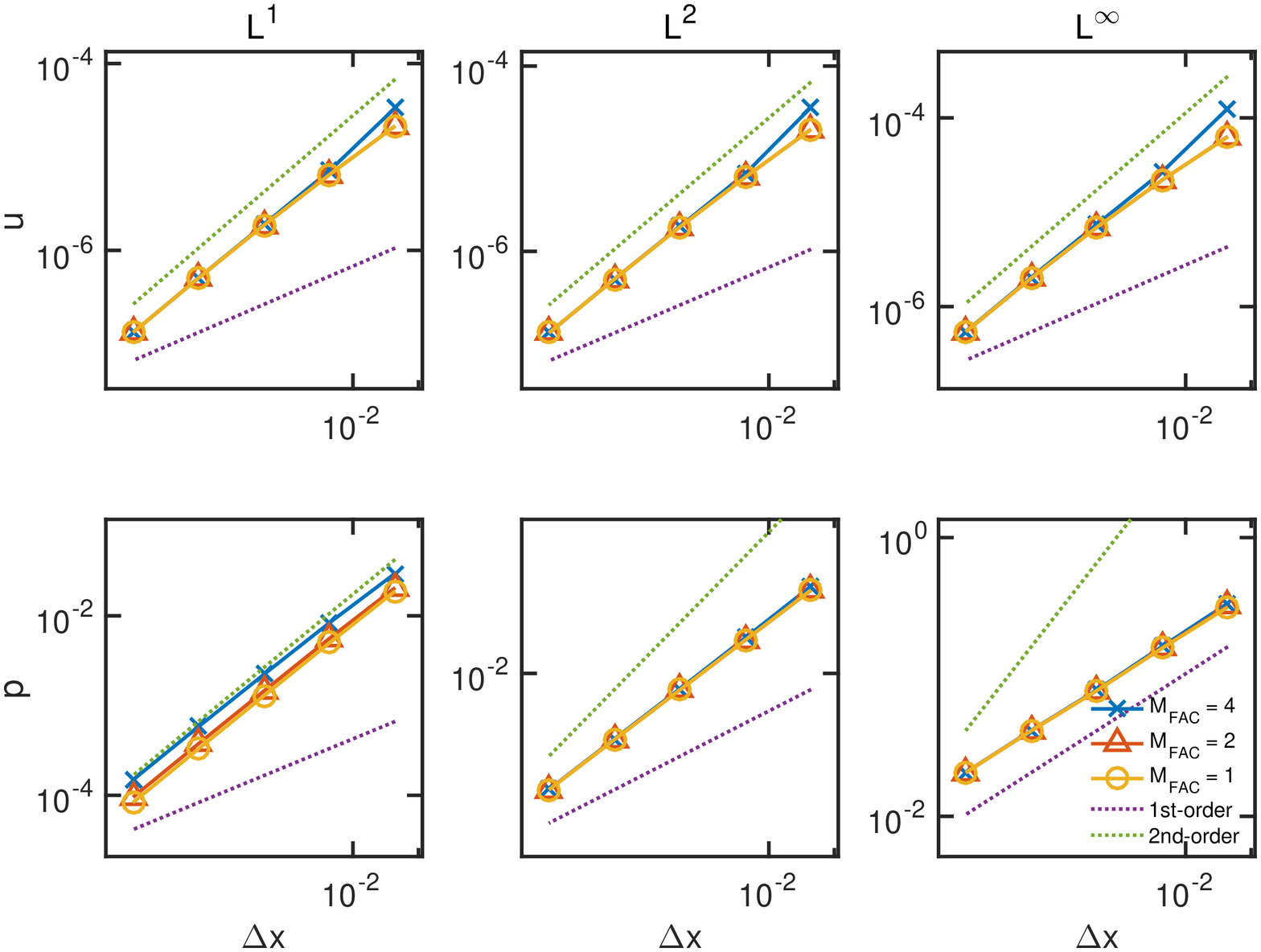}
  \caption{
  	Errors in $\u$ and $p$ in 
  	$L^1$, $L^2$, and $L^{\infty}$ norms for the static ($\gamma=0$) version of the idealized anisotropic shell model of sec.~\ref{s:anisotropic_shell}.
  	Reference lines with slopes of -1 and -2 are also shown.
  	The velocity converges at second-order accuracy in all norms, whereas the pressure converges at second-order in the $L^1$ norm, at first-order in the $L^\infty$ norm, and at order 1.5 in the $L^2$ norm.
  }
  \label{f:smooth_shell_static_error}
\end{figure}

\begin{figure}
  \centering
  \includegraphics[width=0.775\textwidth]{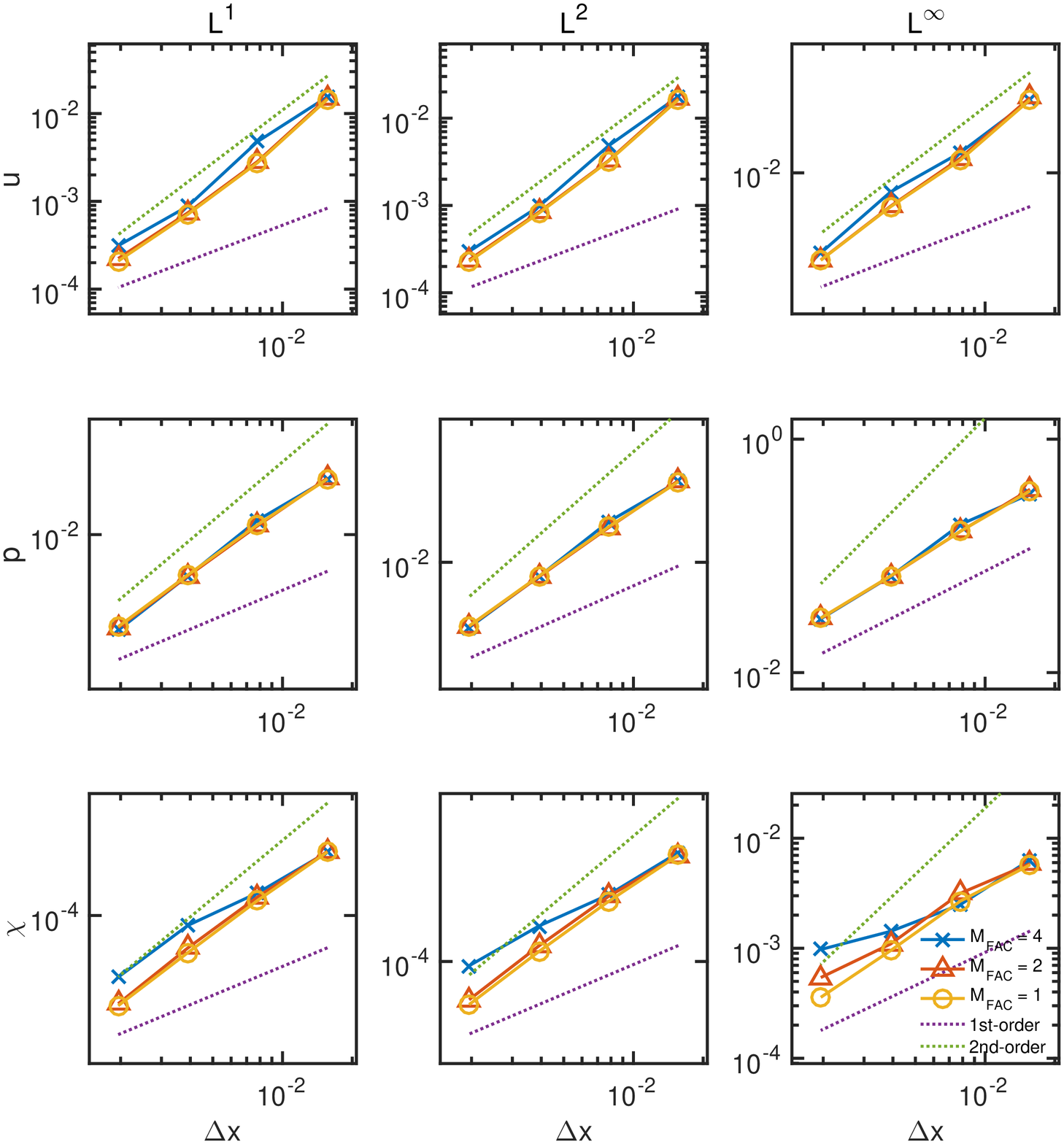}
  \caption{
  	Errors in $\u$, $p$, and $\X$ in 
  	$L^1$, $L^2$, and $L^{\infty}$ norms for the dynamic ($\gamma=0.15$) version of the idealized anisotropic shell model of sec.~\ref{s:anisotropic_shell}.
  	Reference lines with slopes of -1 and -2 are also shown.
  	The velocity converges at second-order accuracy in all norms, whereas the pressure converges at second-order in the $L^1$ norm, at first-order in the $L^\infty$ norm, and at order 1.5 in the $L^2$ norm.
  	The displacement converges at essentially second-order rates in the $L^1$ and $L^2$ norms, and at a slightly lower rate in the $L^\infty$ norm.
  }
  \label{f:smooth_shell_dynamic_error}
\end{figure}

We first consider an idealized anisotropic material model defined in terms of
\begin{equation}
  \We(\FF) = \frac{\mu^{\text{e}}}{2w} \left\| \D{\X}{s_1} \right\|^2 = \frac{\mu^{\text{e}}}{2w} \FF_{\alpha 1}\FF_{\alpha 1},
    \label{e:anisotropic_shell_model}
\end{equation}
so that
\begin{equation}
  \PPe = \D{\We}{\FF} = \frac{\mu^{\text{e}}}{w} \begin{pmatrix} \D{\chi_1}{s_1} & 0 \\ \D{\chi_2}{s_1} & 0 \end{pmatrix} = \frac{\mu^{\text{e}}}{w} \begin{pmatrix} \FF_{11} & 0 \\ \FF_{21} & 0 \end{pmatrix}.
\end{equation}
This model corresponds to an idealized 
elastic material composed of a continuous family of extension-resistant fibers that wrap the thick shell.
Because $U$ is periodic in the $s_1$ direction, $\PPe \, \N \equiv 0$ along $\p U$.
If we view the structure as a fiber-reinforced material, none of the fibers terminate along the boundary of the structure.
Because the transmission force vanishes in this case, the unified and partitioned weak formulations are identical.

Setting $\gamma = 0$, so that the structure is in equilibrium, and requiring $\int_\Omega p(\x,t) \, \Dx = 0$, it can be shown \cite{DBoffi08} that
\begin{equation}
  p(\x,t) = \begin{cases}
    p_0 + \frac{\mu^{\text{e}}}{R}                     &       r \leq R,   \\
    p_0 + \frac{\mu^{\text{e}}}{w} \frac{1}{R} (R+w-r) & R   < r \leq R+w, \\
    p_0                                                & R+w < r,
  \end{cases}
\end{equation}
with $r = \| \x - (0.5,0.5)\|$ and $p_0 = \frac{\pi\mu^{\text{e}}}{3w}\left(R^2 - \frac{(R+w)^3}{R}\right)$.
We set $\rho = 1$, $\mu = 1$, and $\mu^{\text{e}} = 1$, and we consider the time interval $0 \leq t \leq 3$.
Fig.~\ref{f:smooth_shell_static_error} summarizes the error data at time $t=3$ for $N = 64$, $128$, $256$, $512$, and $1024$, using $\Mfac = 1$, $2$, and $4$, with $\dt = 0.25 \dx$.
Second-order convergence rates are observed in the 
$L^1$, $L^2$, and $L^\infty$ norms for the velocity field.
Second-order convergence rates are also observed for the pressure in the 
$L^1$ norm; however, because the pressure field is $C^0$ but not $C^1$, only first-order convergence rates are observed for the pressure in the 
$L^\infty$ norm, and intermediate convergence rates (approximately order 1.5) are observed in the $L^2$ norm.

We also consider the case in which $\gamma = 0.15$, so that the initial configuration of the shell is not in equilibrium.
We set $\rho = 1$, $\mu = 0.01$, and $\mu^{\text{e}} = 1$, yielding a Reynolds number of approximately 50.
We consider the time interval $0 \leq t \leq 0.75$, which corresponds to approximately one damped oscillation of the shell.
Because an exact solution is not available in this case, we use a Richardson extrapolation approach, as described in detail in previous work \cite{BEGriffith05-ib_accuracy}.
Fig.~\ref{f:smooth_shell_dynamic_error} summarizes the error data at time $t = 0.75$ for $N = 64$, $128$, $256$, and $512$ and $\Mfac = 1$, $2$, and $4$, with $\dt = 0.25 \dx$.
Essentially second-order convergence rates are observed in the 
$L^1$, $L^2$, and $L^\infty$ norms for the velocity field.
The pressure converges at a second-order rate in the $L^1$ norm, at a first-order rate in the $L^\infty$ norm, and at an intermediate rate (approximately 1.5) in the $L^2$ norm.
Convergence rates for the deformation are somewhat less regular, with nearly second-order convergence rates being observed in the 
$L^1$ and $L^2$ norms and between first-~and second-order convergence rates observed in the 
$L^\infty$ norm.
The robustness of the convergence rate in the deformation can be improved by using higher-order structural elements.
Because the overall accuracy of the discretization is also limited by the Eulerian discretization and the form of the regularized kernel function, however, the use of higher-order elements does not in itself increase the overall order of accuracy of the method.

Notice that in both the static and dynamic test cases, virtually identical errors are attained for all of the values of $\Mfac$ considered.
This indicates that for these tests, the method is able to use relatively coarse structural meshes without appreciable loss in accuracy.
In particular, these results suggest that the scheme does not allow leaks at fluid-structure interfaces, even for Lagrangian meshes that are quite coarse compared to the Eulerian grid.

\subsubsection{Orthotropic shell}
\label{s:orthotropic_shell}

\begin{figure}
  \centering
  \includegraphics[width=0.775\textwidth]{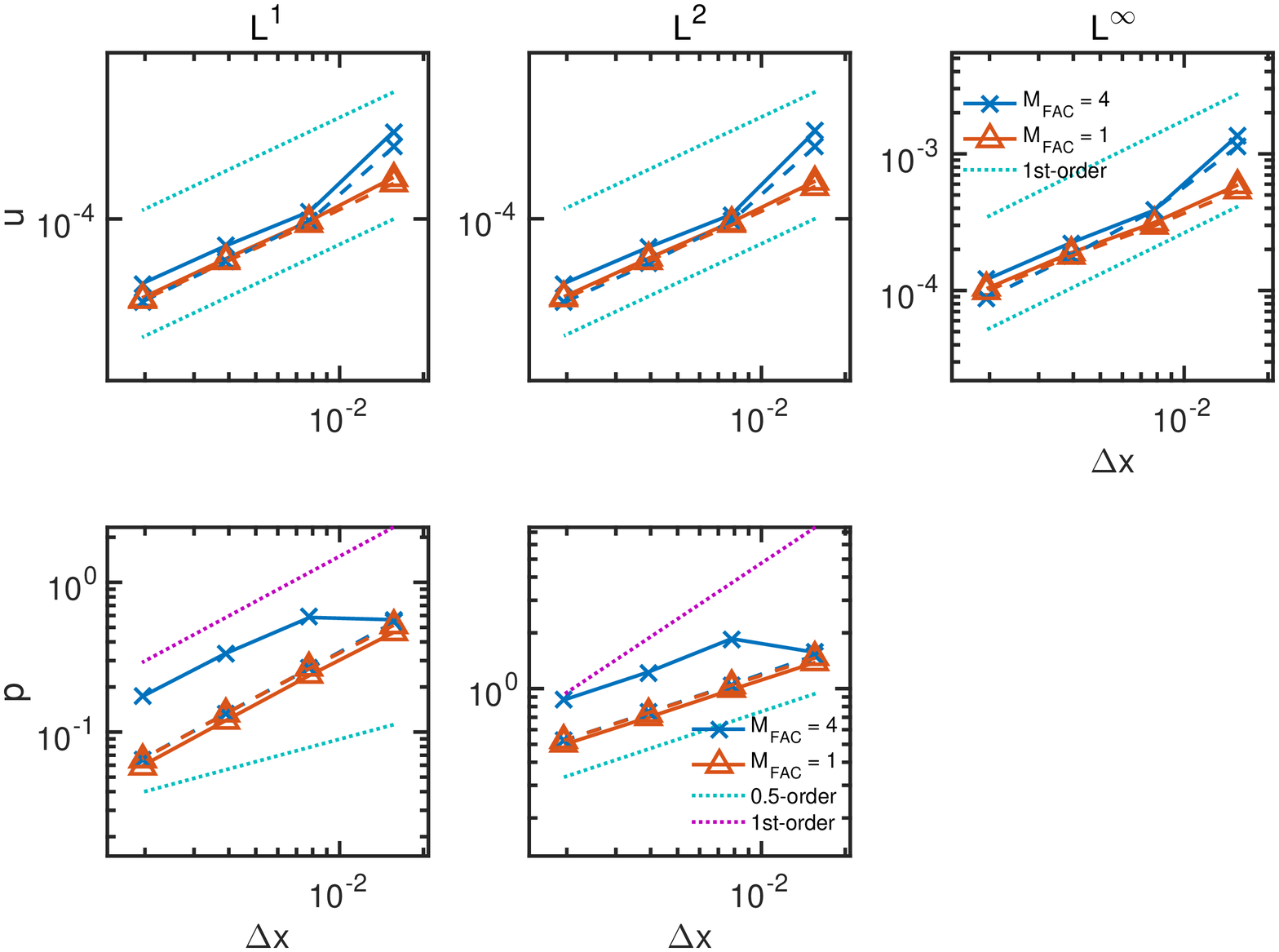}
  \caption{
  	Errors in $\u$ and $p$ in 
  	$L^1$, $L^2$, and $L^{\infty}$ norms for the static ($\gamma=0$) version of the orthotropic shell model of sec.~\ref{s:orthotropic_shell}.
  	Errors for the unified formulation appear as solid lines, and errors for the partitioned formulation appear as dashed lines.
  	Reference lines with slope -1 are provided for the $\u$ error data.
  	For $p$, reference lines with slopes -1 and -0.5 are provided for the $L^1$ and $L^2$ norm data, respectively.
  	Because this test includes discontinuities in the pressure at fluid-structure interfaces, the present method does not yield convergence in $p$ in the 
  	$L^{\infty}$ norm.
  	The partitioned formulation generally yields improved accuracy compared to the unified formulation, especially for the pressure for relatively coarse structural mesh spacings.
  }
  \label{f:sharp_shell_static_error}
\end{figure}

\begin{figure}
  \centering
  \includegraphics[width=0.775\textwidth]{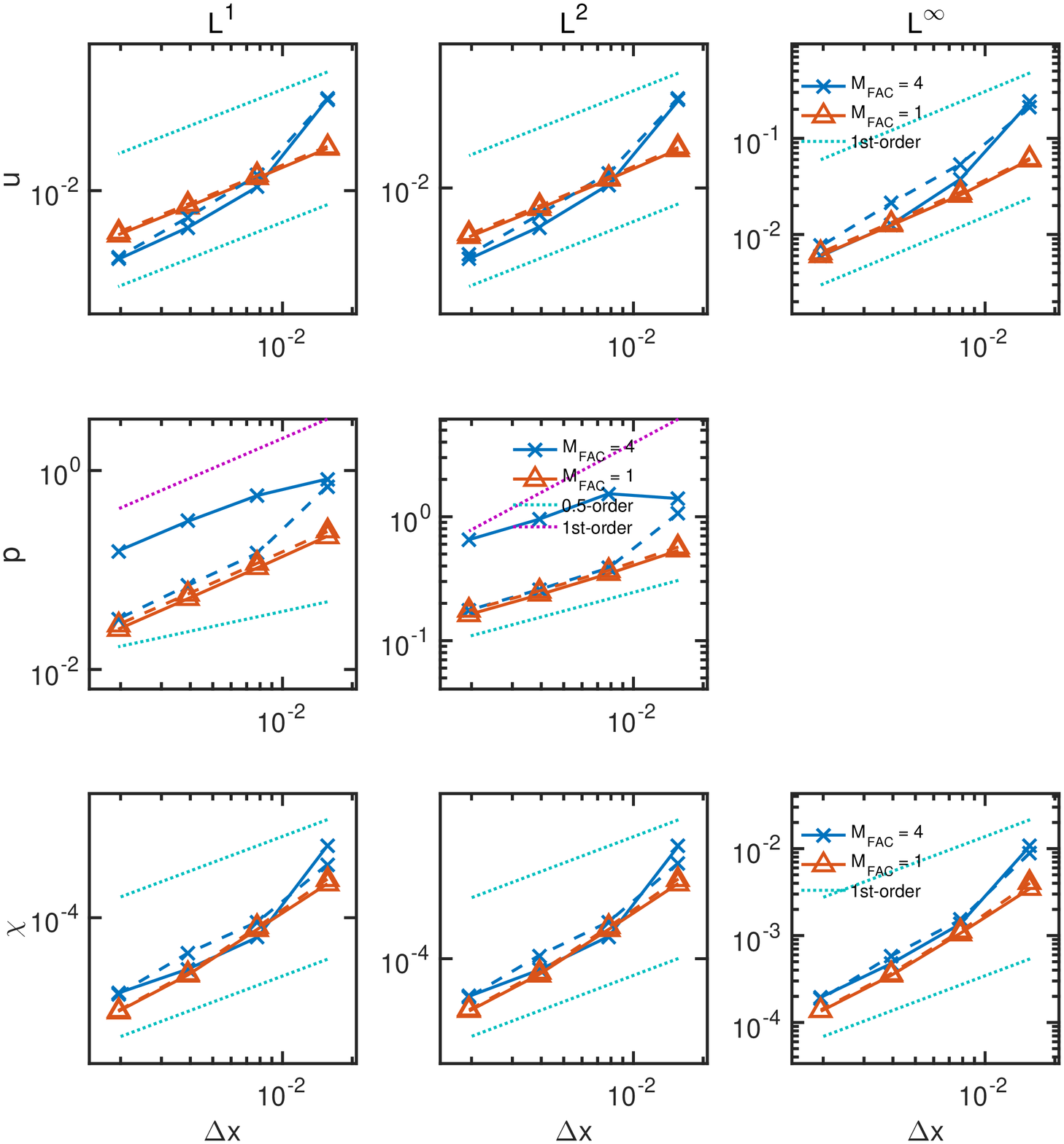}
  \caption{
  	Errors in $\u$, $p$, and $\X$ in 
  	$L^1$, $L^2$, and $L^{\infty}$ norms for the dynamic ($\gamma=0.15$) version of the orthotropic shell model of sec.~\ref{s:orthotropic_shell}.
  	Errors for the unified formulation appear as solid lines, and errors for the partitioned formulation appear as dashed lines.
  	Reference lines with slope -1 are provided for the $\u$ and $\X$ error data.
  	For $p$, reference lines with slopes -1 and -0.5 are provided for the $L^1$ and $L^2$ norm data, respectively.
  	The unified formulation generally yields modest improvements in the accuracy for $\u$ and $\X$, whereas the partitioned formulation generally yields improved accuracy for $p$, especially for relatively coarse structural mesh spacings.
  }
  \label{f:sharp_shell_dynamic_error}
\end{figure}

The second case that we consider uses a neo-Hookean material model,
\begin{equation}
  \We(\FF) = \frac{\mu^{\text{e}}}{2w} I_1(\CC), \label{e:orthotropic_shell_model}
\end{equation}
with $\CC = \FF^{\mathrm T} \FF$ and $I_1(\CC) = \mathrm{tr}(\CC)$, so that
\begin{equation}
  \PPe = \frac{\mu^{\text{e}}}{w} \FF.
\end{equation}
Because of the form of the mapping from curvilinear coordinates to initial coordinates, this material behaves as an orthotropic fiber-reinforced solid rather than as an isotropic material.
One way of viewing the elastic response is that the body is composed of two continuous families of fibers.
The first family of fibers wraps the elliptical shell circumferentially, and the second family is composed of radial fibers that are orthogonal to the circumferential fibers.
 Because one family of fibers terminates along the fluid-structure interfaces, there are singular force layers along $\p \X(U,t)$ that must be balanced by discontinuities in the pressure and viscous stress.
 Therefore, in this case the discretized unified and partitioned formulations yield different results.

Setting $\gamma = 0$, so that the structure is in equilibrium, and requiring $\int_{\Omega} p(\x,t) \, \Dx = 0$, it can be shown \cite{DBoffi08}
that
\begin{equation}
  p(\x,t) = \begin{cases}
    p_0 + \mu^{\text{e}} \left(\frac{1}{R} - \frac{1}{R+w}\right)                   &       r \leq R,   \\
    p_0 + \frac{\mu^{\text{e}}}{w} \left(\frac{1}{R} (R+w-r) + \frac{R}{R+w}\right) & R   < r \leq R+w, \\
    p_0                                                                             & R+w < r,
  \end{cases}
\end{equation}
with $r = \| \x - (0.5,0.5)\|$ and $p_0 = \frac{\pi\mu^{\text{e}}}{3w}\left(3 w R + R^2 - \frac{(R+w)^3}{R}\right)$.
We set $\rho = 1$, $\mu = 1$, and $\mu^{\text{e}} = 1$, and we consider the time interval $0 \leq t \leq 3$.
Fig.~\ref{f:sharp_shell_static_error} summarizes the error data at time $t=3$ for $N = 64$, $128$, $256$, $512$, and $1024$, using $\Mfac = 1$, $2$, and $4$, with $\dt = 0.25 \dx$.
First-order convergence rates are observed for $\u$ in all norms.
First-order convergence rates are also observed for $p$ in the $L^1$ norm.
Because $p$ possesses discontinuities at fluid-structure interfaces for this problem, however, the present method yields convergence rates of 0.5 in the $L^2$ norm and does not converge in the $L^{\infty}$ norm.

We also consider the case in which $\gamma = 0.15$, so that the initial configuration of the shell is not in equilibrium.
We set $\rho = 1$, $\mu = 0.01$, and $\mu^{\text{e}} = 1$, yielding a Reynolds number of approximately 100.
We consider the time interval $0 \leq t \leq 1.25$, which corresponds to approximately one damped oscillation of the shell.
Again, an exact solution is not available, and so convergence rates are estimated using Richardson extrapolation \cite{BEGriffith05-ib_accuracy}.
Fig.~\ref{f:sharp_shell_dynamic_error} summarize the error data at time $t = 0.75$ for $N = 64$, $128$, $256$, and $512$ and $\Mfac = 1$ and $4$, with $\dt = 0.25 \dx$.
Essentially first-order convergence rates are observed for $\u$ and $\X$ in all norms, whereas $p$ exhibits first-order convergence in only the $L^1$ norm.

For this problem, we find that the unified and partitioned formulations yield similar accuracy in most cases for $\u$ and $\X$.
By contrast, the partitioned formulation offers significantly better accuracy for the pressure for relatively coarse Lagrangian meshes.
This property appears also to result in improvements in volume conservation; see sec.~\ref{s:soft_elastic_disc}.

\subsection{Soft elastic disc in lid driven cavity}
\label{s:soft_elastic_disc}

\begin{figure}
	\centering
	\includegraphics[width=0.175\textwidth]{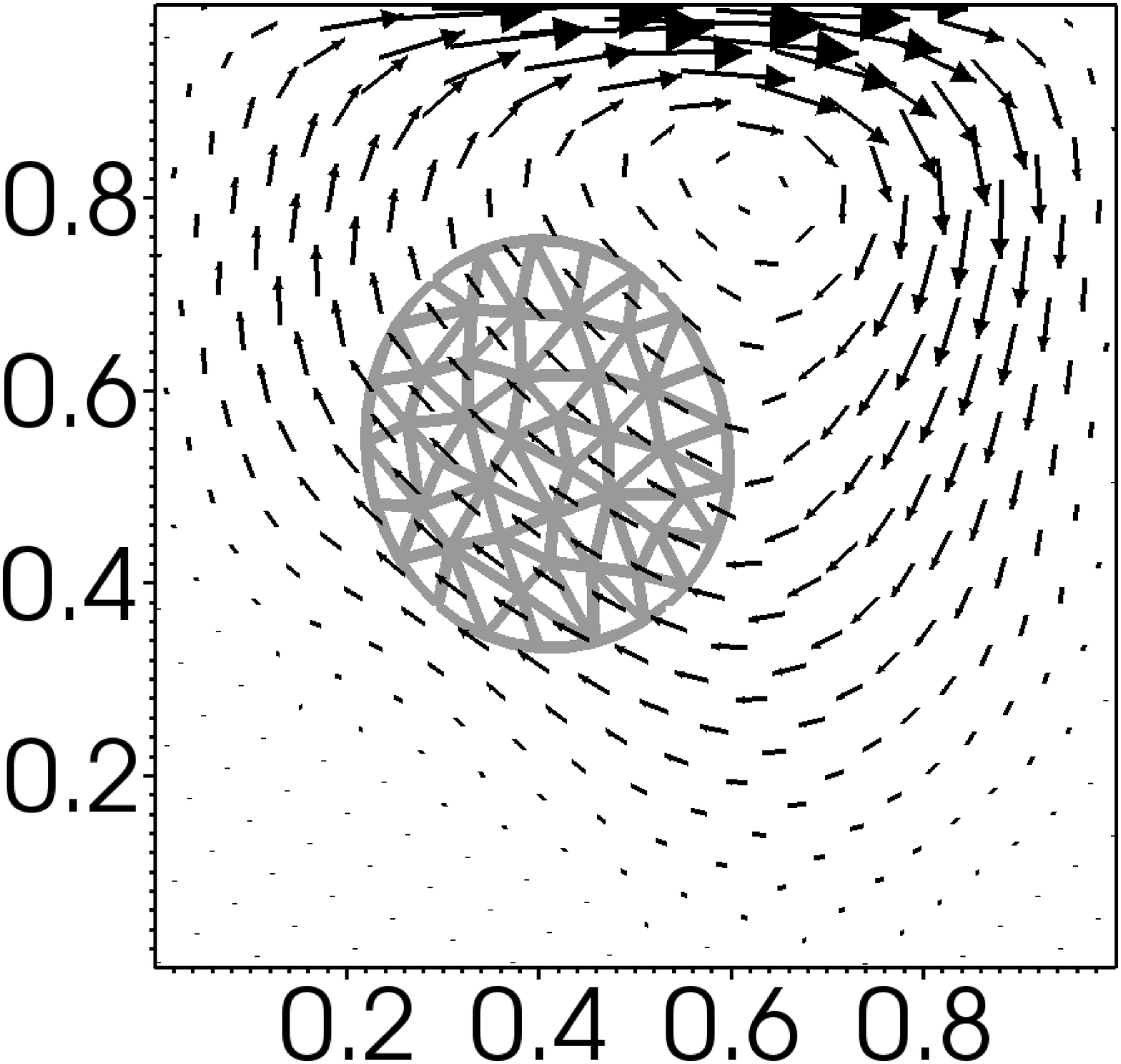} \
	\includegraphics[width=0.175\textwidth]{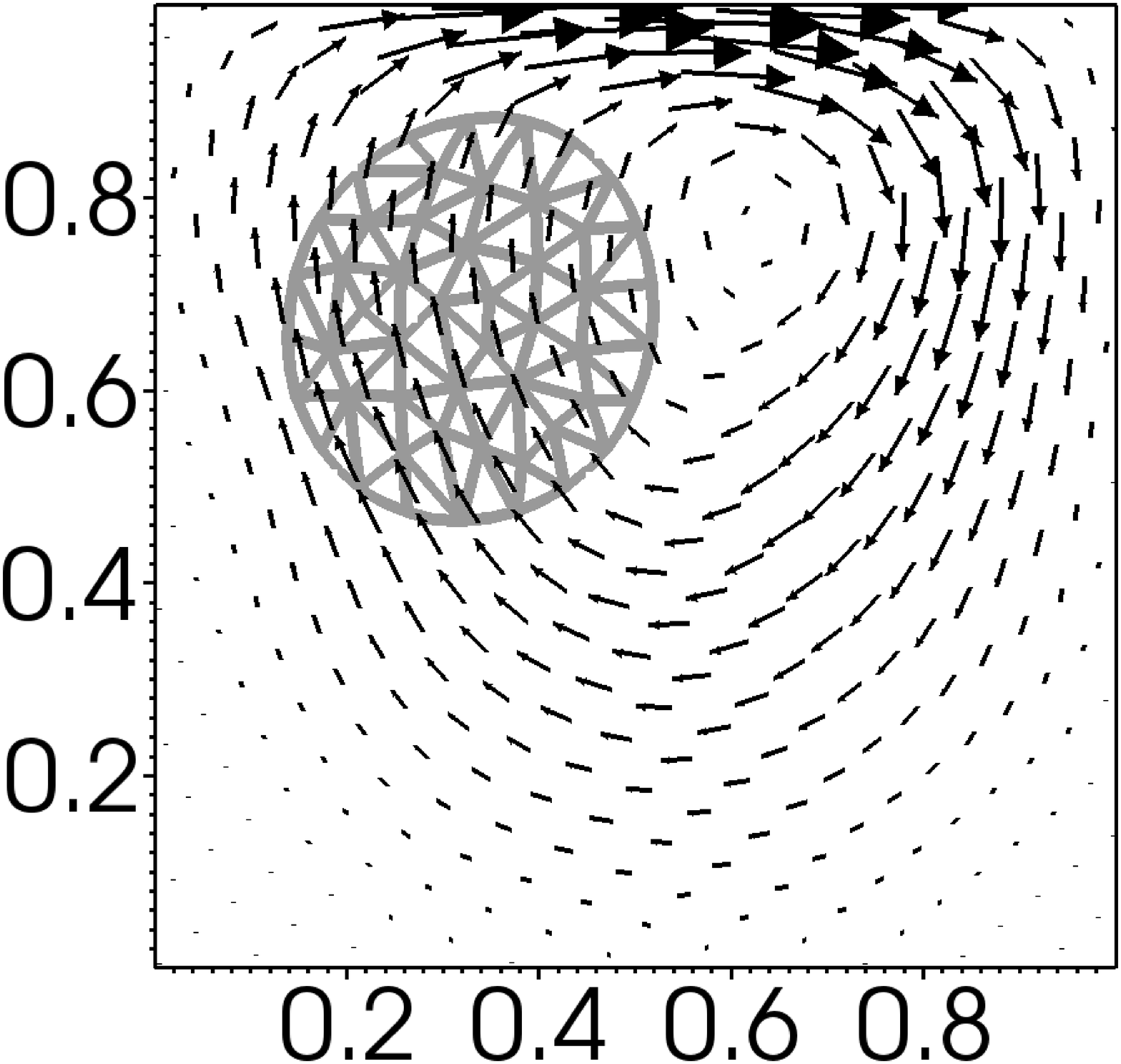} \
	\includegraphics[width=0.175\textwidth]{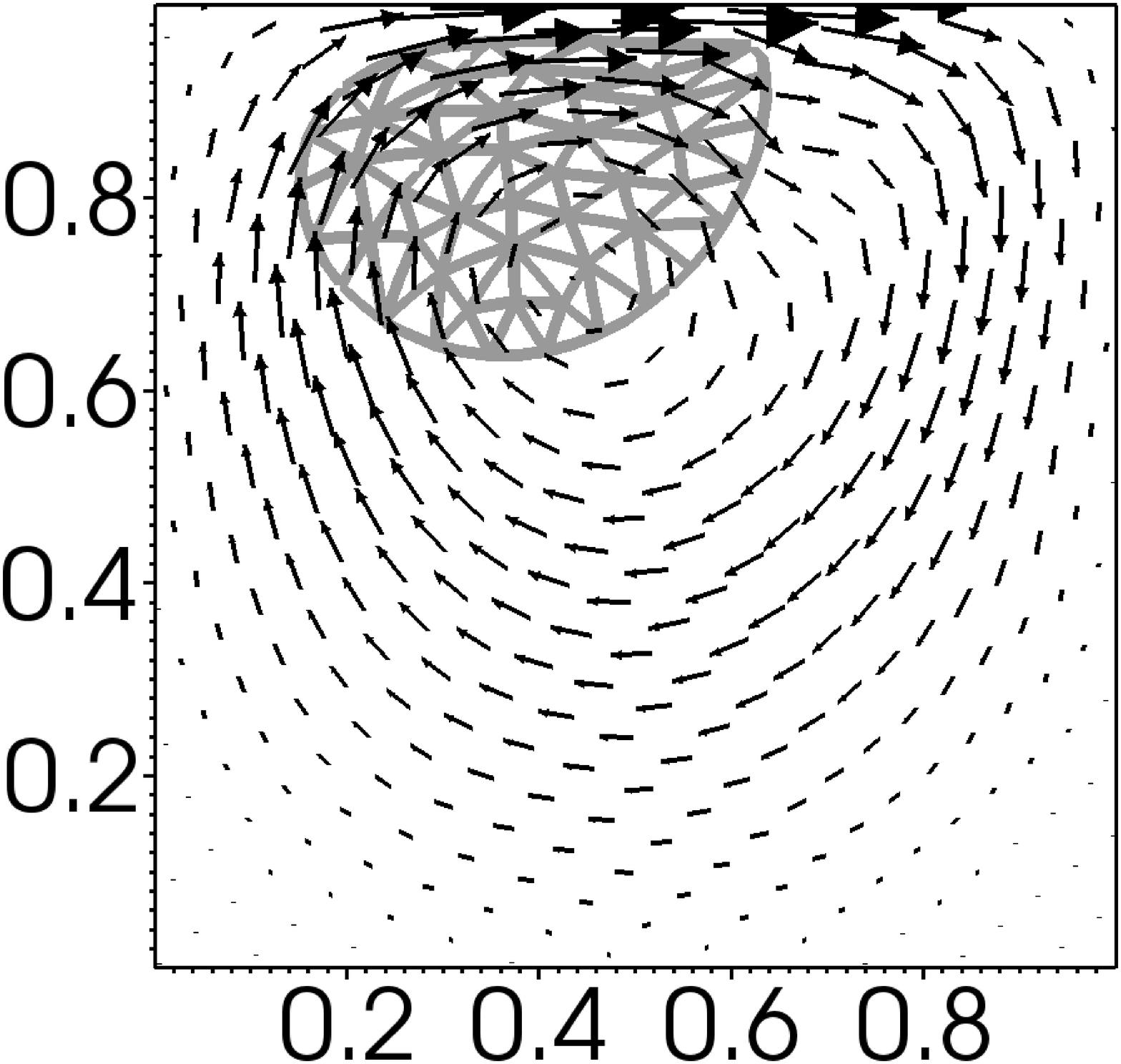} \
	\includegraphics[width=0.175\textwidth]{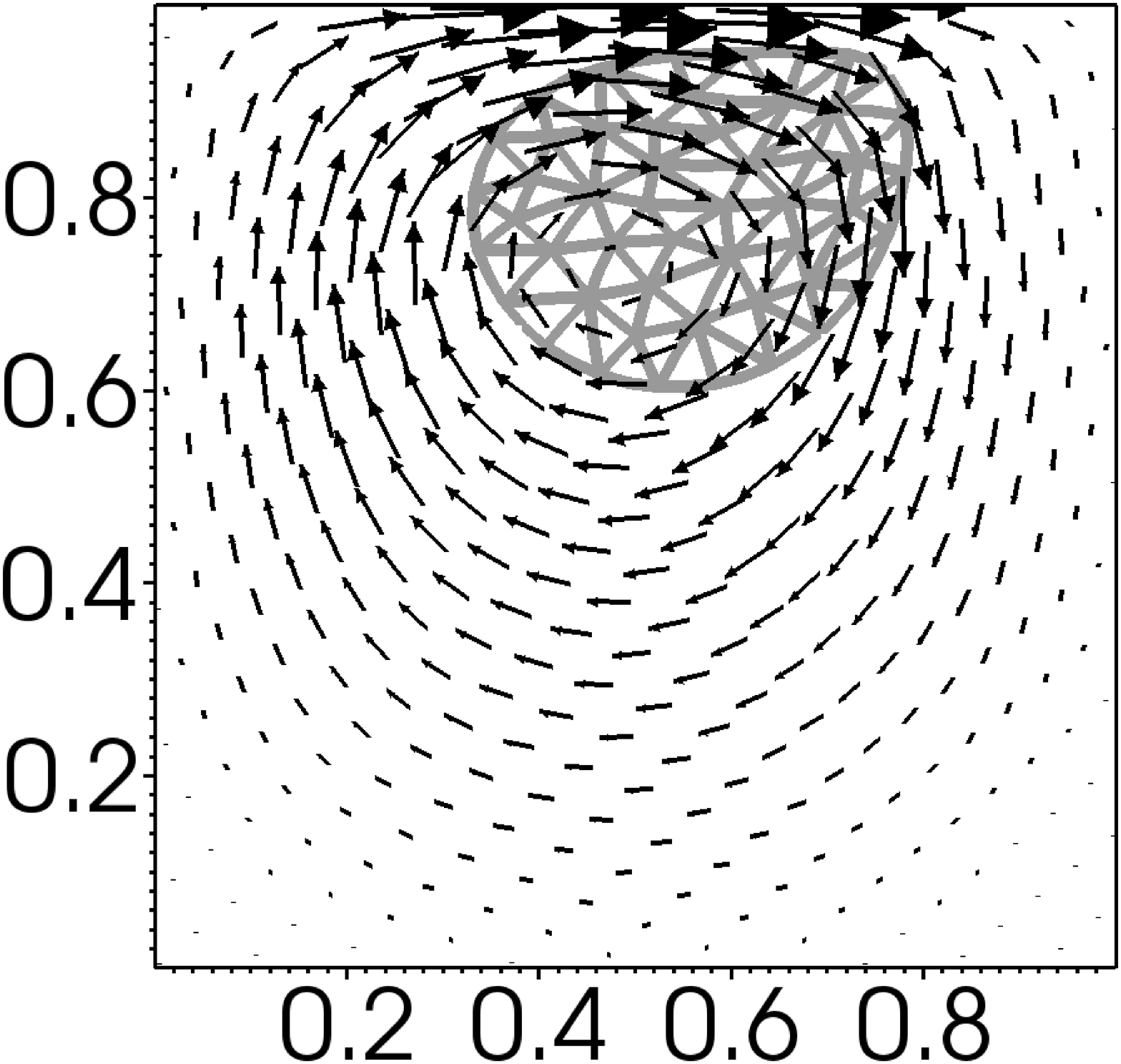} \
	\includegraphics[width=0.175\textwidth]{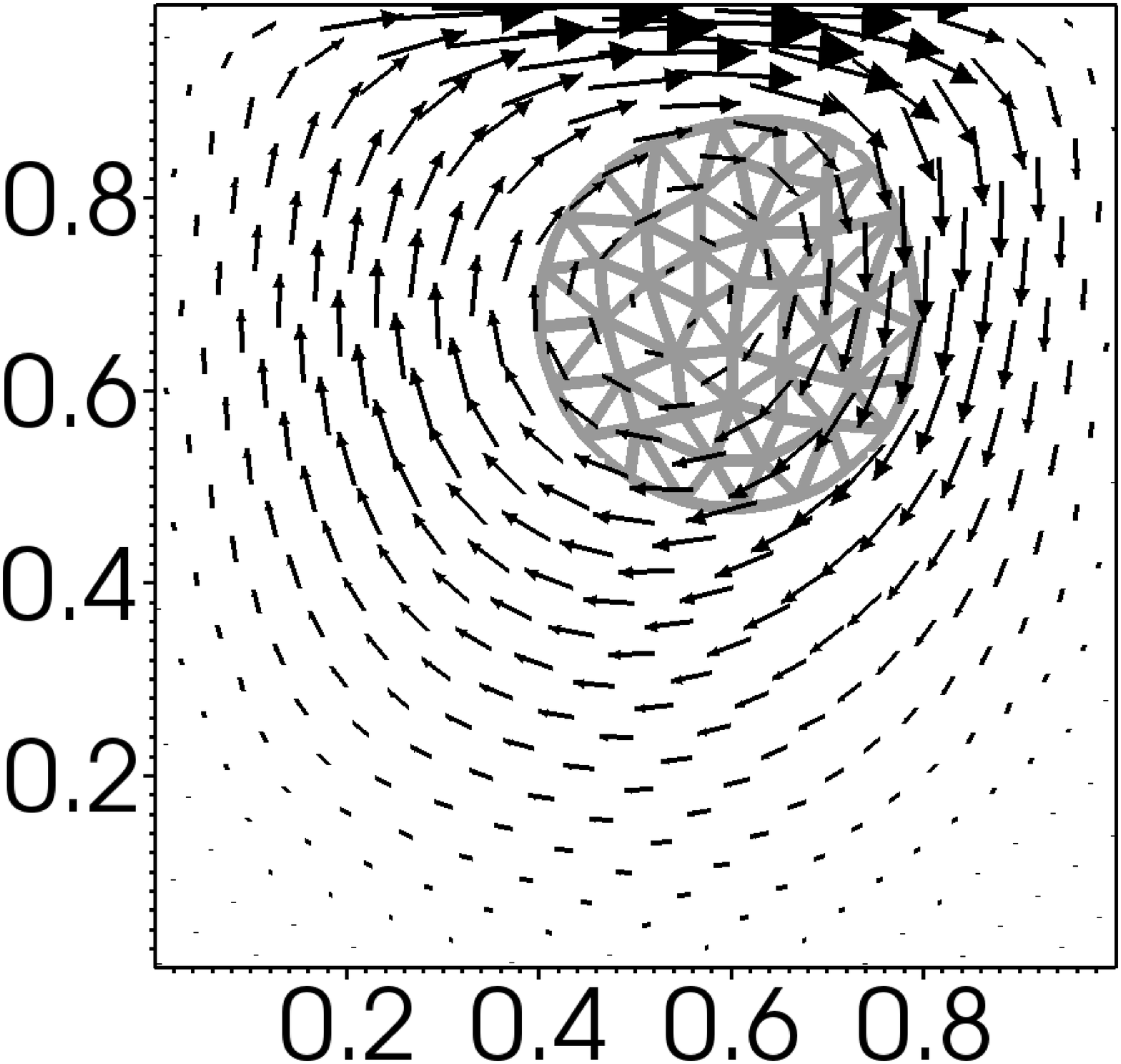}
	\caption{
		A soft neo-Hookean disc in a lid-driven cavity flow using the partitioned (split) weak formulation with $N=128$ and $\Mfac = 4$ over the time interval $3 \leq t \leq 7$.
	}
	\label{f:soft_disc}
\end{figure}

This section presents results from tests that use a soft elastic structure in a lid-driven cavity flow to demonstrate the volume-conservation properties of our method.

In these computations, the physical domain is $\Omega = [0,1] \times [0,1]$, and the immersed structure is a disc of radius $0.2$ that is initially centered about $\x = (0.6,0.5)$.
The boundary conditions imposed along $\p\Omega$ are $\u \equiv 0$ on the left ($x_1 = 0$), right ($x_1 = 1$), and bottom ($x_2 = 0$) boundaries of $\Omega$, and $\u \equiv (1,0)$ on the top ($x_2 = 1$) boundary of $\Omega$.
We use an isotropic neo-Hookean model,
\begin{equation}
  \PPe = \mu^\text{e} \, \FF - p_0 \, \FF^{-{\mathrm T}}, \label{e:soft_disc_stress}
\end{equation}
and we consider both $p_0 = 0$ and $p_0 = \mu^\text{e}$.
Because generally $\PPe \, \N \not \equiv 0$, the solution has discontinuities in the pressure and viscous stress at fluid-structure interfaces.
In such cases, we expect the IB method to yield no better than first-order convergence rates.\footnote{
	This flow also possesses well-known corner singularities that reduce the convergence rate of the incompressible flow solver.
	Although it is possible to devise numerical schemes that accurately treat the corner singularities present in the classical lid-driven cavity flow \cite{BotellaPeyret98}, we do not employ such a method in this work.
}
The flow induced by the driven lid brings the structure nearly into contact with the moving upper boundary of the domain; see fig.~\ref{f:soft_disc}.
This near contact is handled automatically by the IB formulation using a version of a modified kernel function approach introduced by Griffith et al.~\cite{BEGriffith09-heart_valves} and enhanced by Kallemov et al.~\cite{BKallemov16-RigidIBAMR}.
No additional specialized methods are required by the present scheme to handle this case.

\begin{figure}
	\centering
	\begin{tabular}{>{\centering\arraybackslash}m{0.1\textwidth} >{\centering\arraybackslash}m{0.8\textwidth} }
	{\bf A.} \vspace{-30pt} & \\ 
	$p_0 = 0$ & {\includegraphics[width=0.775\textwidth]{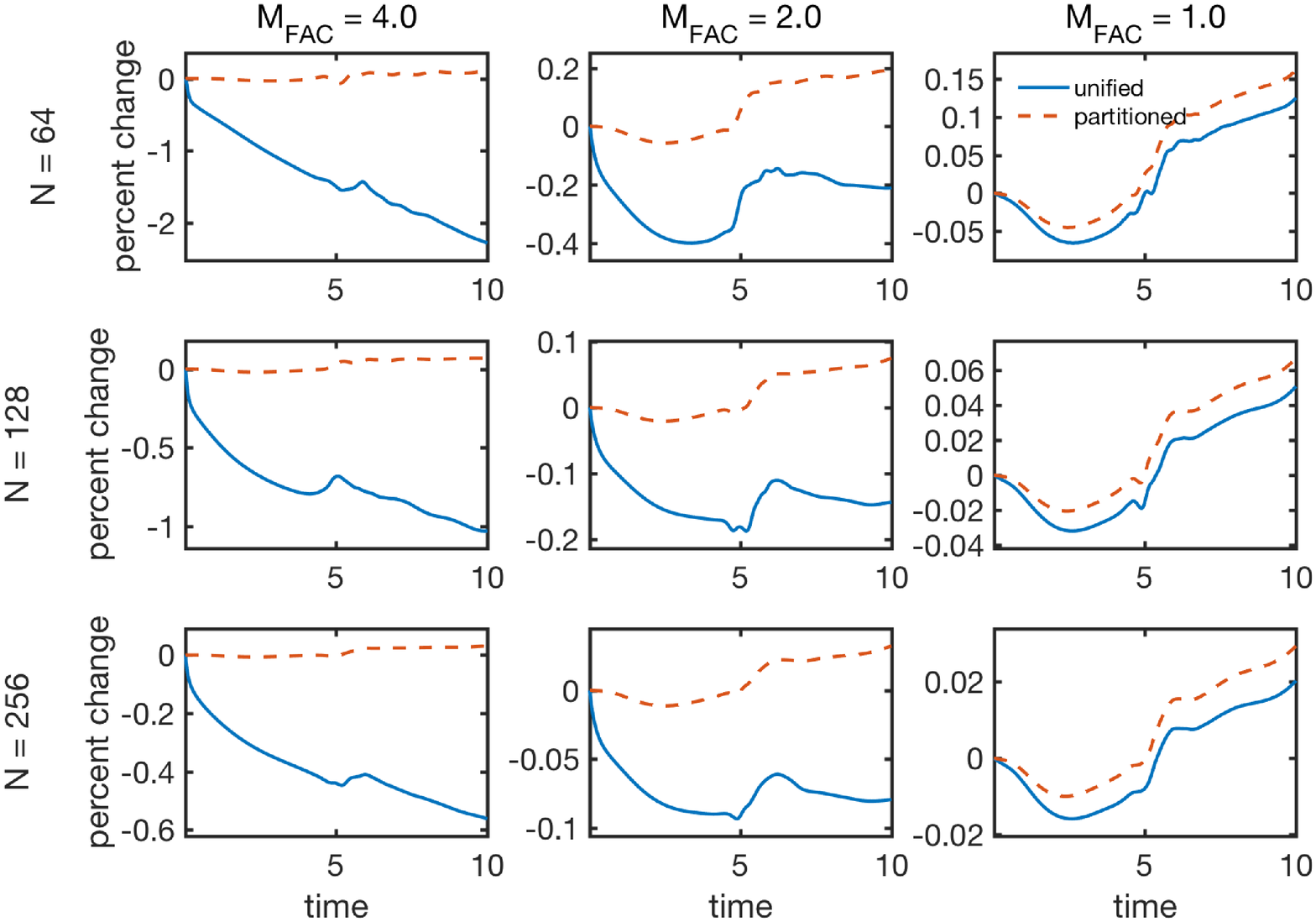}}	\vspace{-1.1\baselineskip} \\
	{\bf B.} \vspace{-30pt} & \\
	$p_0 = \mu^\text{e}$ & {\includegraphics[width=0.775\textwidth]{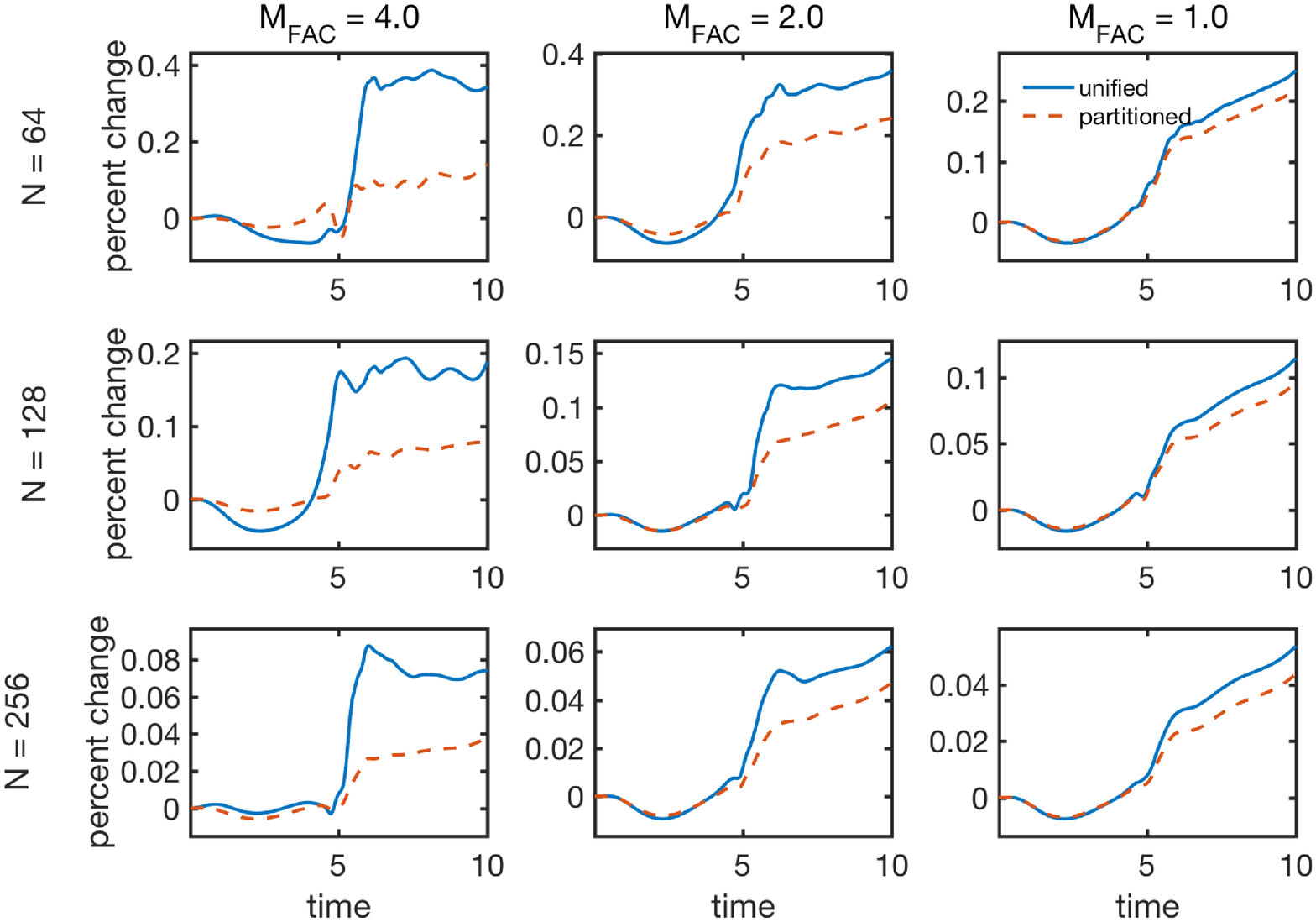}}	\vspace{-1.1\baselineskip} \\
  \end{tabular}
	\caption{
		Percent change in structure volume for the soft elastic disc benchmark of sec.~\ref{s:soft_elastic_disc} as a function of time using eq.~\eqref{e:soft_disc_stress} with {\bf A.}~$p_0 = 0$ and {\bf B.}~$p_0=\mu^\text{e}$.
		Results are shown for Cartesian grids of size $N = 64$, $128$, and $256$ with relative Lagrangian mesh spacings $\Mfac = 4$, $2$, and $1$.
		The amount of spurious volume change converges to zero at first order.
		At coarser relative mesh spacings, the partitioned (split) formulation yields substantially better volume (area in two spatial dimensions) conservation than the unified (unsplit) formulation.
		Volume conservation is substantially improved for the unified formulation by setting $p_0 = \mu^\text{e}$.
		The volume change is less than 0.4\% for all cases considered except for the unified formulation with $p_0 = 0$ and $\Mfac = 4$.
	}
	\label{f:soft_disc_volume}
\end{figure}

As in previous studies of this test case \cite{HZhao08,XSWang10}, we set $\mu = 0.01$ and $\rho = 1$.
We consider $\mu^{\text{e}} = 0.2$, which is a relatively ``soft'' case.
The initial velocity is $\u \equiv 0$, and the reference coordinates $\s$ are taken to be the initial coordinates, so that $\X(\s,0) \equiv \s$.
The physical domain is discretized using an $N \times N$ Cartesian grid.
The Lagrangian coordinate domain is discretized using an unstructured mesh of quadratic triangular ($P^2$) elements constructed so that the elements are approximately a factor of $\Mfac$ coarser than the background Eulerian grid, so that in the reference configuration, the nodes of the Lagrangian mesh are physically separated by a distance of approximately $\Mfac \dx$.
The time step size is $\dt = 0.125 \dx$.
We consider the time interval $0 \leq t \leq 10$, during which the disc is subjected to slightly more than one rotation within the cavity.
The structure becomes entrained in the shearing flow along the cavity lid from $t \approx 4$ until $t \approx 6$, and during this time is subjected to very large deformations.
Fig.~\ref{f:soft_disc_volume} shows the percent change in disc volume for different values of $N$ and $\Mfac$.
With $p_0 = 0$, the maximum volume change yielded by the unified formulation is approximately 2.3\% for $N=64$ and $\Mfac=4$, 0.4\% for $N=64$ and $\Mfac=2$, and 0.2\% for $N=64$ and $\Mfac=1$.
The split formulation yields substantially improved accuracy: the maximum volume change is approximately 0.12\% for $N=64$ and $\Mfac=4$, 0.2\% for $N=64$ and $\Mfac=2$, and 0.15\% for $N=64$ and $\Mfac=1$.
Setting $p_0 = \mu^\text{e}$ improves the accuracy of the unified formulation substantially, especially at coarser relative mesh spacings, whereas it results in slightly poorer volume conservation for the split formulation.
With $p_0 = \mu^\text{e}$, the maximum volume change is less than 0.4\% in all cases considered.
At smaller values of $\Mfac$, there is essentially no difference in the volume changes produced by the unified and partitioned formulations.
In all cases, the maximum volume change converges to zero at a first-order rate.

In general, the partitioned formulation appears to yield superior volume conservation for this problem, especially at relatively coarse Lagrangian mesh spacings. 
Moreover, the volume-conservation properties of the partitioned scheme seem to be relatively insensitive to the relative coarseness of the Lagrangian mesh.
Using either formulation, volume errors converge to zero at essentially a first-order rate.
These results compare very favorably to those obtained by the IFE method, which, even for relatively fine Lagrangian meshes, can yield volume losses of up to 20\% when applied to the same test without using a volume-conservation algorithm.
A modification of the IFE method to improve its volume conservation still yields volume losses of approximately 2.5\% for this test \cite{XSWang10}.

\subsection{Flow past a cylinder}
\label{s:flow_past_cylinder}

\begin{figure}
	\centering
    \begin{tabular}{lrlr}
    \vspace{-10pt} {\bf A.} & & {\bf B.} & \\
	& \includegraphics[height=125pt]{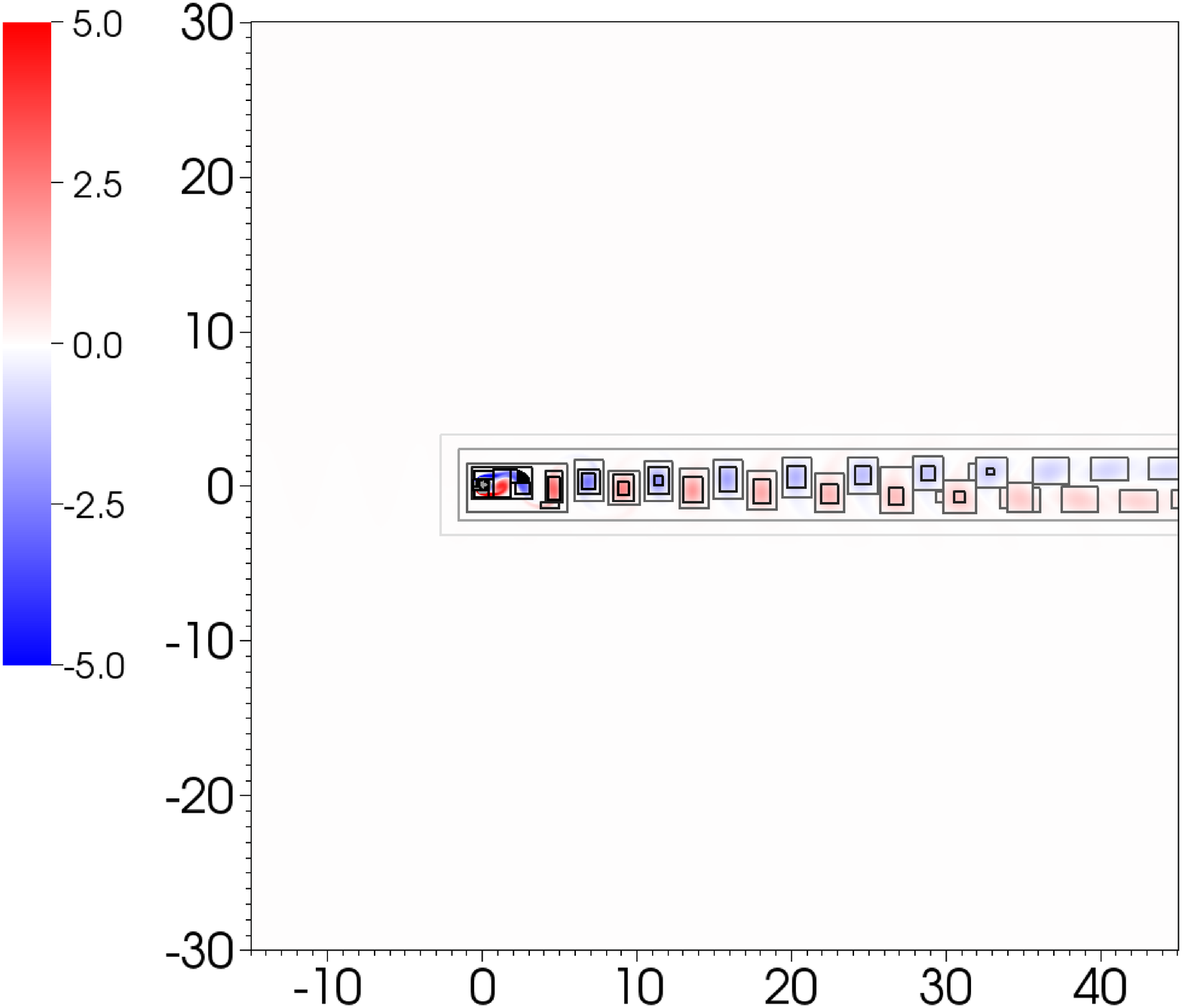} & & \includegraphics[height=125pt]{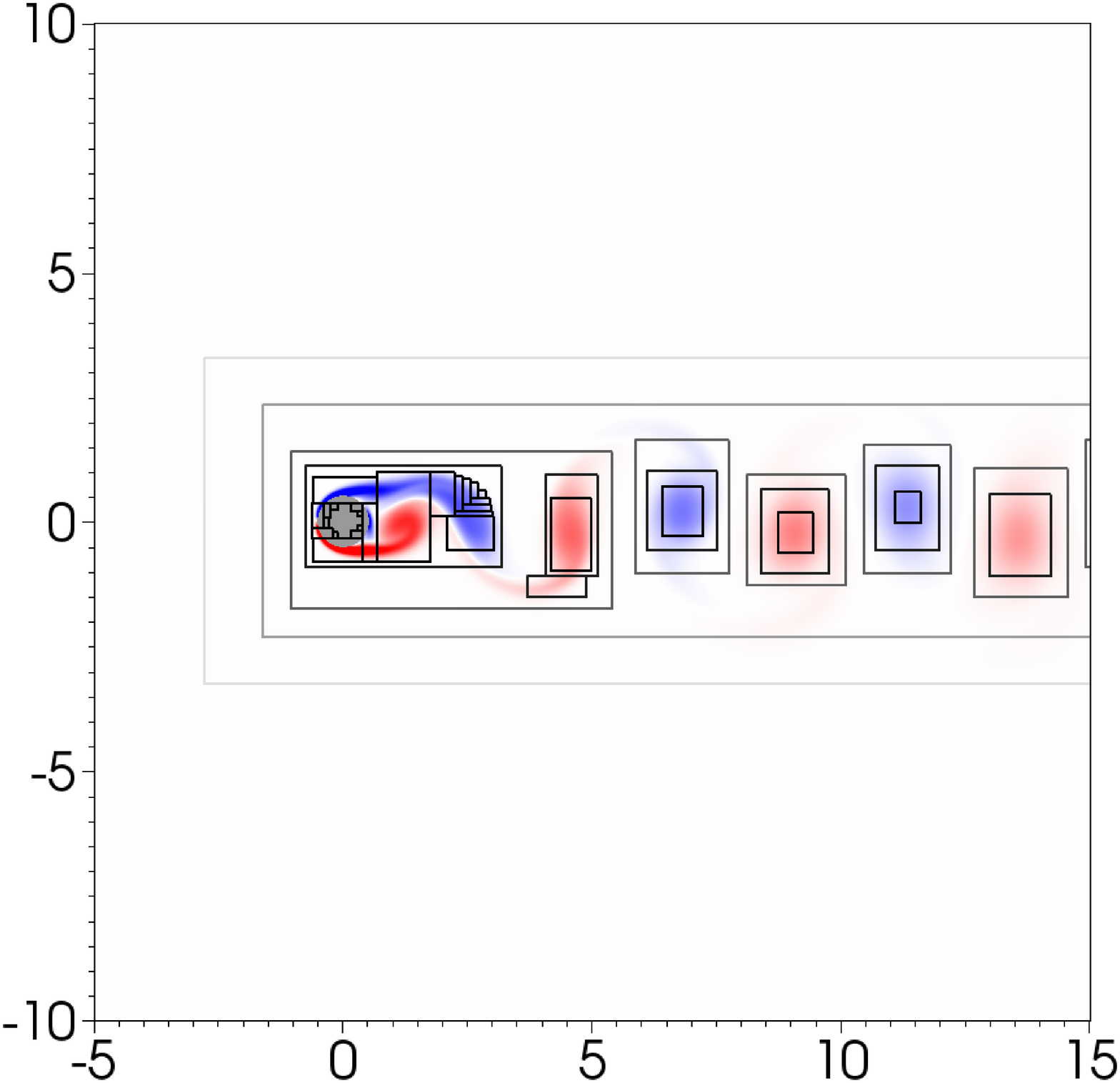}
	\end{tabular}
	\caption{
		Vortices shed from a stationary circular cylinder at $Re = 200$.
		This computation uses a six-level locally refined grid with a refinement ratio of two between levels and an $N \times N$ coarse grid with $N = 128$.
		Regions of local mesh refinement are indicated by rectangular boxes, with lighter grey boxes indicating coarser levels of refinement and black boxes indicating the finest regions of the locally refined grid.
		{\bf A.}~The full computational domain and {\bf B.}~a magnified view of the flow field near the cylinder.	
	}
	\label{f:flow_past_cylinder}
\end{figure}

This section presents results from tests using the widely used benchmark of viscous flow past a stationary circular cylinder.
In these computations, the physical domain is $\Omega = [-15,45]\times[-30,30]$, and the immersed structure is a thin circular interface of radius $0.5$ centered about $\x = (0,0)$.
This domain size and cylinder placement corresponds to `Case C' considered by Taira and Colonius \cite{KTaira07}.
Along the inflow boundary ($x_1 = -15$), we set a uniform inflow velocity, $\u \equiv (1,0)$.
Along the outflow boundary ($x_1 = 45$), we set the normal traction and tangential velocity to zero, whereas along the top ($x_2 = 30$) and bottom ($x_2 = -30$) boundaries, we set the normal velocity and tangential traction to zero.
The boundary condition treatment is the same as that described by Griffith \cite{BEGriffith09-efficient}.
We set $\rho = 1$ and $\mu = 0.005$.
Using the inflow velocity as the characteristic velocity $\u_\infty$ and the cylinder diameter $d$ as the characteristic length, the Reynolds number is $Re = \frac{\rho u_\infty d}{\mu} = 200$.
The computational domain is discretized using an adaptively refined Cartesian grid with six nested grid levels and a refinement ratio of two between levels.
The Cartesian grid spacing on the finest grid level is $\dx_\text{finest} = 2^{-5} \dx_\text{coarsest}$, with $\dx_\text{coarsest} = \frac{60}{N}$.
The cylinder is discretized using a mesh of one-dimensional quadratic ($P^2$) elements with a node spacing of approximately $\Mfac \dx_\text{finest}$.
The time step size is $\dt = 0.1 \dx_\text{finest}$, yielding a maximum CFL number of approximately $0.1$--$0.2$.
Rigidity constraints are approximately imposed using tether forces as in eq.~\eqref{e:rigidity}.
For each grid spacing considered, we use approximately the largest stable values of the penalty parameters $\kappa$ and $\eta$ as permitted by the time step size.
These values are empirically determined by a simple optimization procedure.
Representative results are shown in fig.~\ref{f:flow_past_cylinder}.

\begin{table}
\centering
	\begin{tabular}{l|lll}
		& $C_\text{L}$ & $C_\text{D}$ & $St$ \\ \hline
		Lai and Peskin \cite{LaiPeskin00}                 & ---        & ---              & $0.190$ \\
		Linnick and Fasel \cite{MNLinnick05}              & $\pm 0.69$ & $1.34 \pm 0.044$ & $0.197$ \\
		Liu et al. \cite{CLiu98}                          & $\pm 0.69$ & $1.31 \pm 0.049$ & $0.192$ \\
		Taira and Colonius \cite{KTaira07} (Case C)       & $\pm 0.68$ & $1.34 \pm 0.047$ & $0.195$ \\
		Roshko (experimental) \cite{CLiu98}               & ---        & ---              & $0.190$ \\
		Williamson (experimental) \cite{CLiu98}           & ---        & ---              & $0.197$	\\
		present (four-point kernel, $N=128$, $\Mfac = 2$) & $\pm 0.70$ & $1.36 \pm 0.046$ & $0.195$
	\end{tabular}
	\caption{
		Comparison of experimental and computational values of lift ($C_\text{L}$) and drag ($C_\text{D}$) coefficients and Strouhal numbers ($St$) for flow past a circular cylinder at $Re = 200$.
	}
	\label{t:cylinder_comparison}
\end{table}

\begin{figure}
	\centering
	\includegraphics[width=0.775\textwidth]{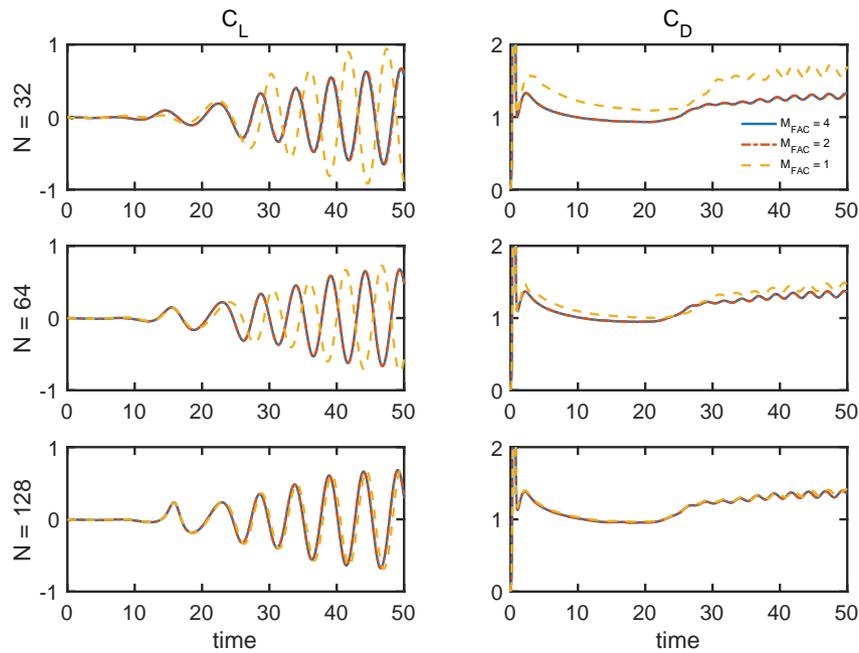}
	\caption{
		Lift ($C_\text{L}$) and drag ($C_\text{D}$) coefficients for flow past a stationary cylinder at $Re = 200$.
		The computational domain $\Omega$ is discretized using a six-level locally refined grid with a refinement ratio of two between levels and an $N \times N$ coarse grid.
		The spacing between the nodes of the immersed structure is approximately $\Mfac \dx_\text{finest}$.
		Under simultaneous Lagrangian and Eulerian grid refinement, the scheme converges to the same dynamics for all values of $\Mfac$.
		Notice, however, that the best accuracy is obtained for \emph{larger} values of $\Mfac$.
		Using a relatively \emph{finer} structural mesh spacing results in larger lift and drag coefficients and lower vortex shedding frequencies.
		Thus, smaller values of $\Mfac$ increase the \emph{effective numerical size} of the cylinder.
	}
	\label{f:lift_and_drag}
\end{figure}

\begin{figure}
	\centering
	\includegraphics[width=0.775\textwidth]{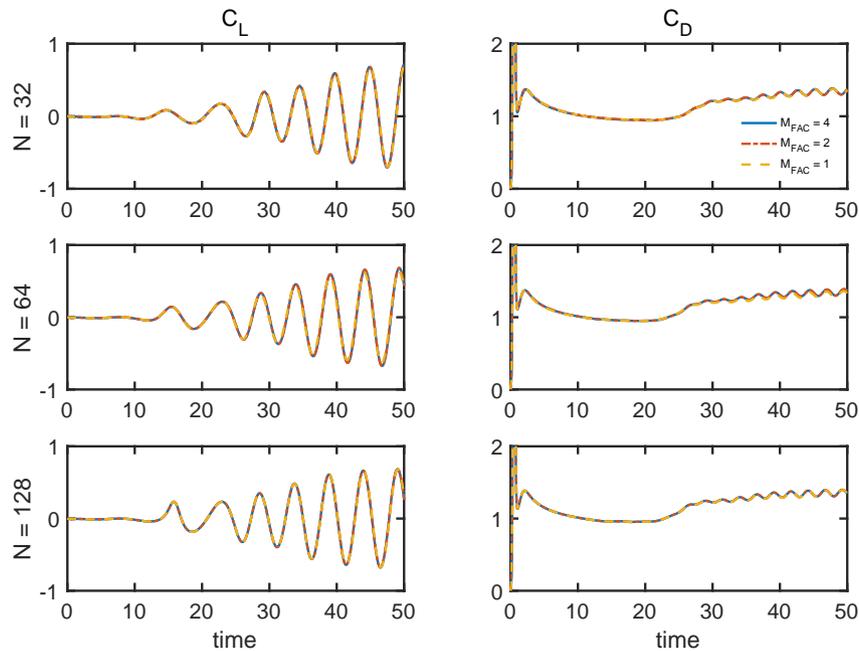}
	\caption{
		Similar to figs.~\ref{f:lift_and_drag} and \ref{f:lift_and_drag-IB_6}, but here using a three-point kernel function \cite{RomaPeskinBerger99}.
		Unlike the case of the four-~ and six-point kernels, comparable accuracy is obtained for all values of $\Mfac$ considered.
	}
	\label{f:lift_and_drag-IB_3}
\end{figure}

\begin{figure}
	\centering
	\includegraphics[width=0.775\textwidth]{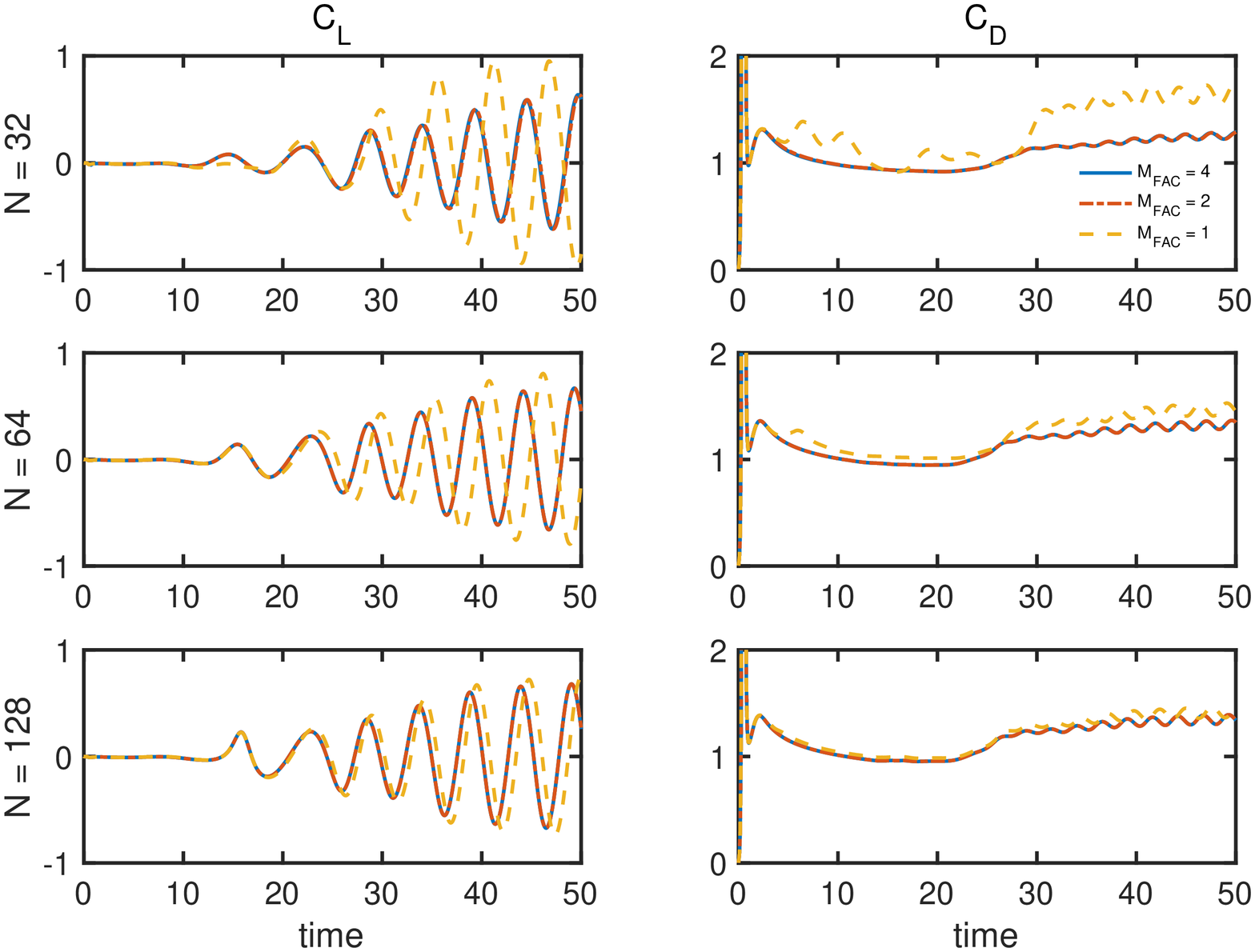}
	\caption{
		Similar to figs.~\ref{f:lift_and_drag} and \ref{f:lift_and_drag-IB_3}, but here using a six-point kernel function with higher continuity order \cite{YBao16-six_point_delta}.
		As with the four-point kernel, the best accuracy is obtained for \emph{larger} values of $\Mfac$.
		In this case, however, using relatively fine structural mesh spacings yields erratic results on coarser Eulerian grids (e.g.~$\Mfac=1$ for $N=32$).
	}
	\label{f:lift_and_drag-IB_6}
\end{figure}

Table \ref{t:cylinder_comparison} compares results obtained using the present method with the standard four-point kernel \cite{Peskin02}, $N = 128$, and $\Mfac = 2$ with corresponding experimental and computational results from prior studies.
We compute the lift coefficient, $C_\text{L} = F_y / (\rho u^2_\infty d/2)$, drag coefficient, $C_\text{D} = F_x / (\rho u^2_\infty d/2)$, and Strouhal number, $St = f_\text{s} d / u_\infty$, in which $\F = (F_x,F_y)$ is the net force on the cylinder and $f_\text{s}$ is the shedding frequency.
The present results are seen to be in excellent quantitative agreement with these earlier results.

Fig.~\ref{f:lift_and_drag} shows the lift 
and drag 
coefficients as functions of time for $N = 32$, $64$, and $128$ and for $\Mfac = 1$, $2$, and $4$ using the four-point regularized kernel function \cite{Peskin02}.
Similar results are shown in fig.~\ref{f:lift_and_drag-IB_3} for the similar three-point kernel \cite{RomaPeskinBerger99}, whereas fig.~\ref{f:lift_and_drag-IB_6} shows results for a recently developed six-point kernel with higher continuity order \cite{YBao16-six_point_delta}.
By construction, the structural meshes obtained for a fixed value of $N$ are nested versions of each other (i.e.~they can be viewed as obtained via Lagrangian mesh refinement).
Under simultaneous Lagrangian and Eulerian grid refinement, the scheme converges to the same dynamics for all values of $\Mfac$.
For the four-~and six-point kernels, however, the best accuracy is obtained for \emph{larger} values of $\Mfac$.
In particular, by using relatively \emph{coarser} structural meshes, the scheme yields more accurate forces ($C_\text{L}$ and $C_\text{D}$) and vortex shedding dynamics (characterized by, e.g., $St$).
In fact, for the coarsest Eulerian discretization considered ($N=32$), the six-point kernel yields erratic results for the finest structural discretization ($\Mfac = 1$) that do not occur for the coarser structural discretizations ($\Mfac = 2$ and $4$).
In contrast, with the three-point kernel, comparable results are obtained for all values of $\Mfac$ considered here.
Although not shown here, results obtained using a two-point piecewise-linear kernel are similar to those obtained using the three-point kernel.

A striking feature of these results is that the use of Lagrangian mesh refinement alone generally \emph{lowers} the accuracy of the computation for a fixed Eulerian grid.
A complete theoretical explanation for this behavior is lacking at present. 

\subsection{Idealized model of left ventricular mechanics}
\label{s:idealized_left_ventricle}

\begin{figure}
	\centering
	\small
	\begin{tabular}{lclclc}
	\vspace{-10pt} {\bf A.} & & {\bf B.} & & {\bf C.} & \\
	& \includegraphics[width=0.2375\textwidth]{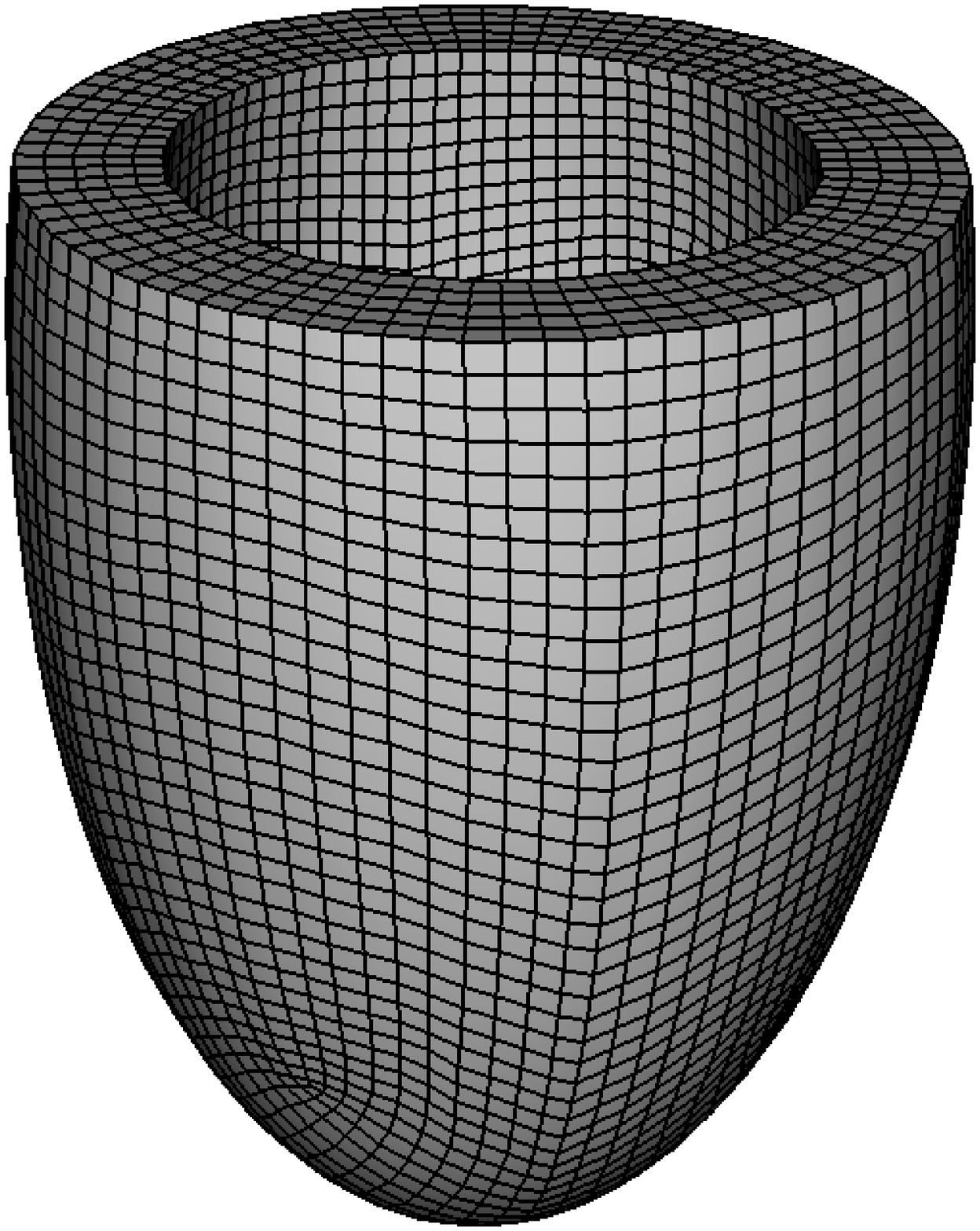} &	
	& \includegraphics[width=0.2375\textwidth]{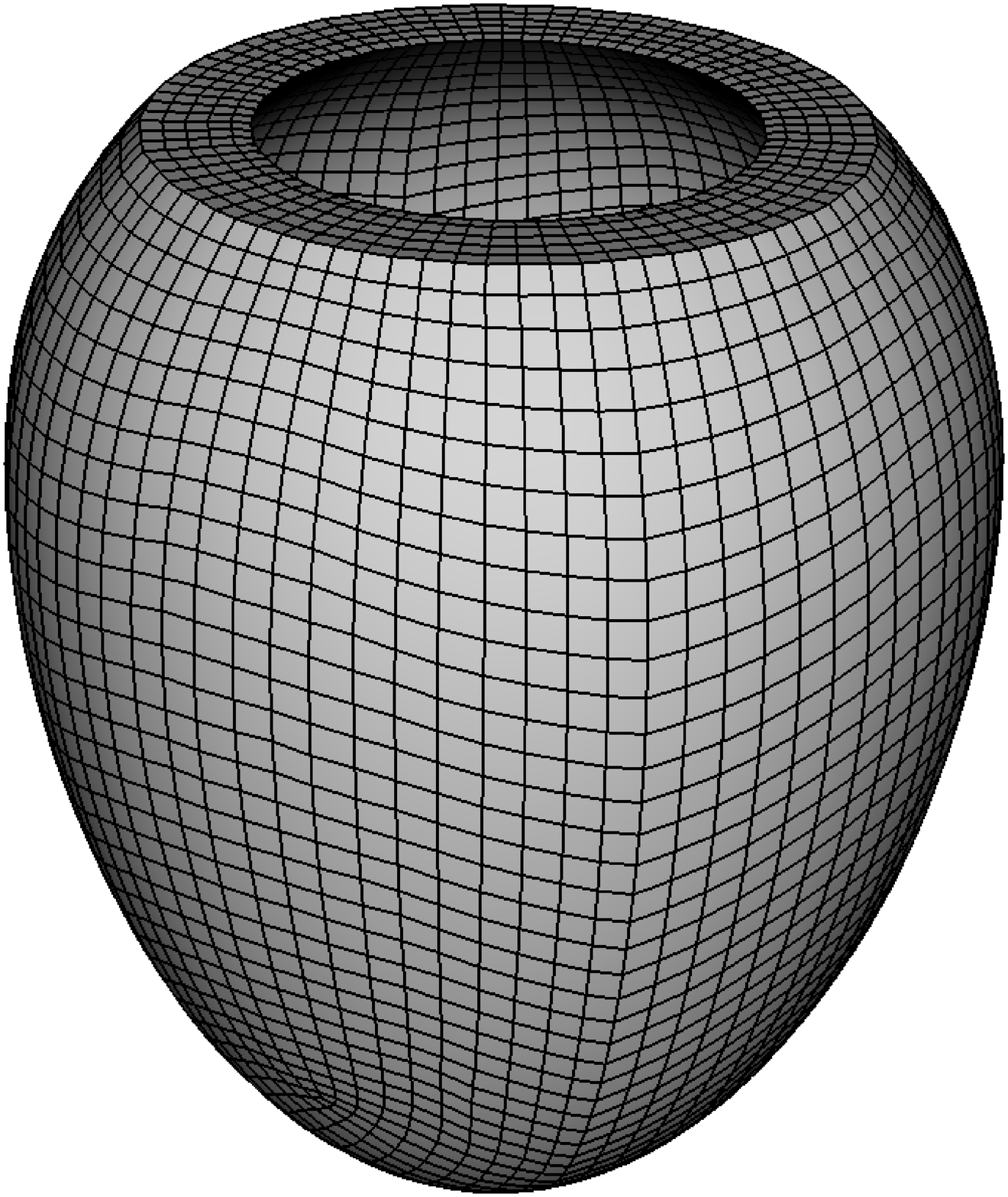} &	
	& \includegraphics[width=0.2375\textwidth]{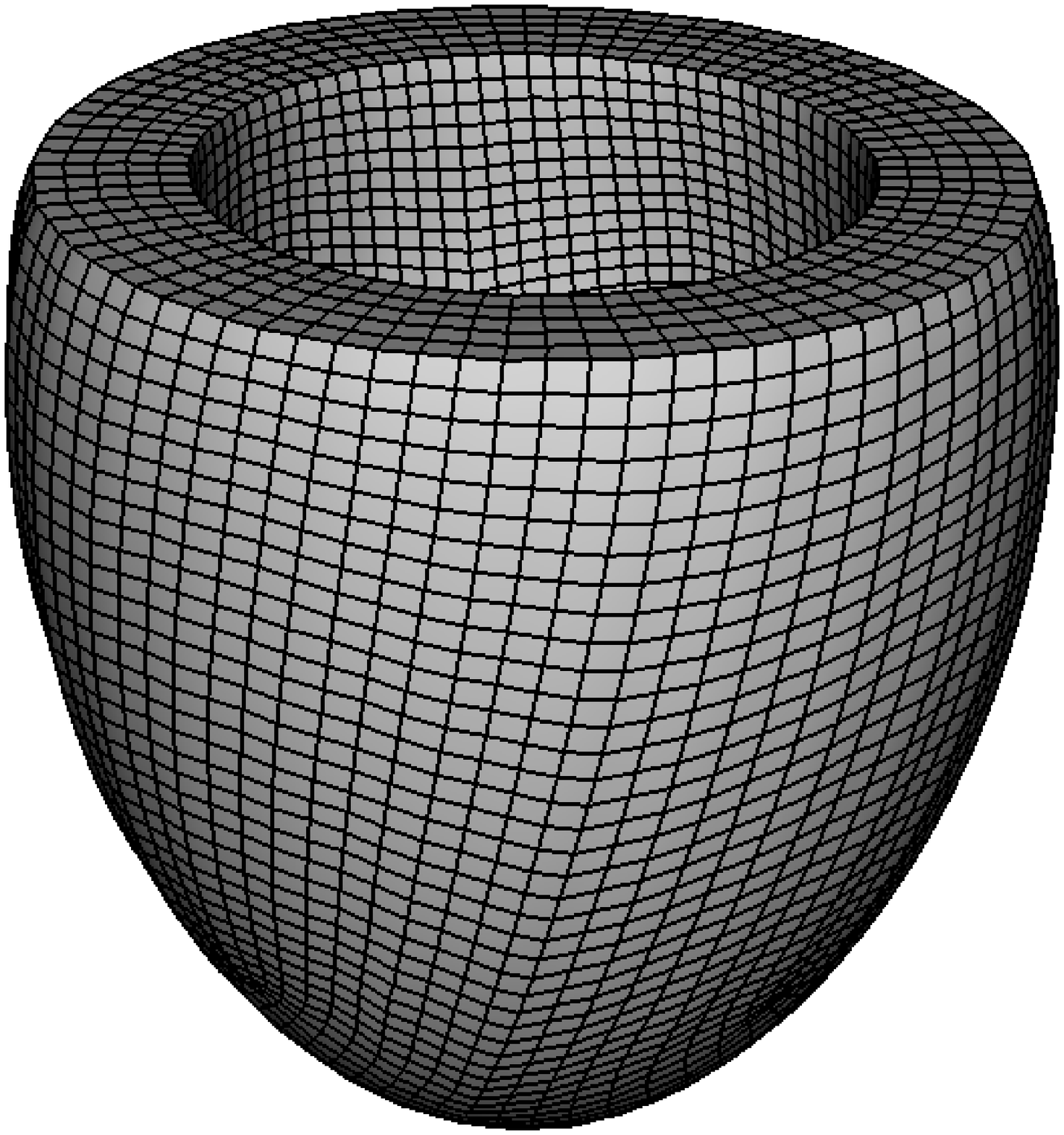} \\
	\vspace{-10pt} {\bf D.} & & {\bf E.} & & {\bf F.} & \\
	& \includegraphics[width=0.2375\textwidth]{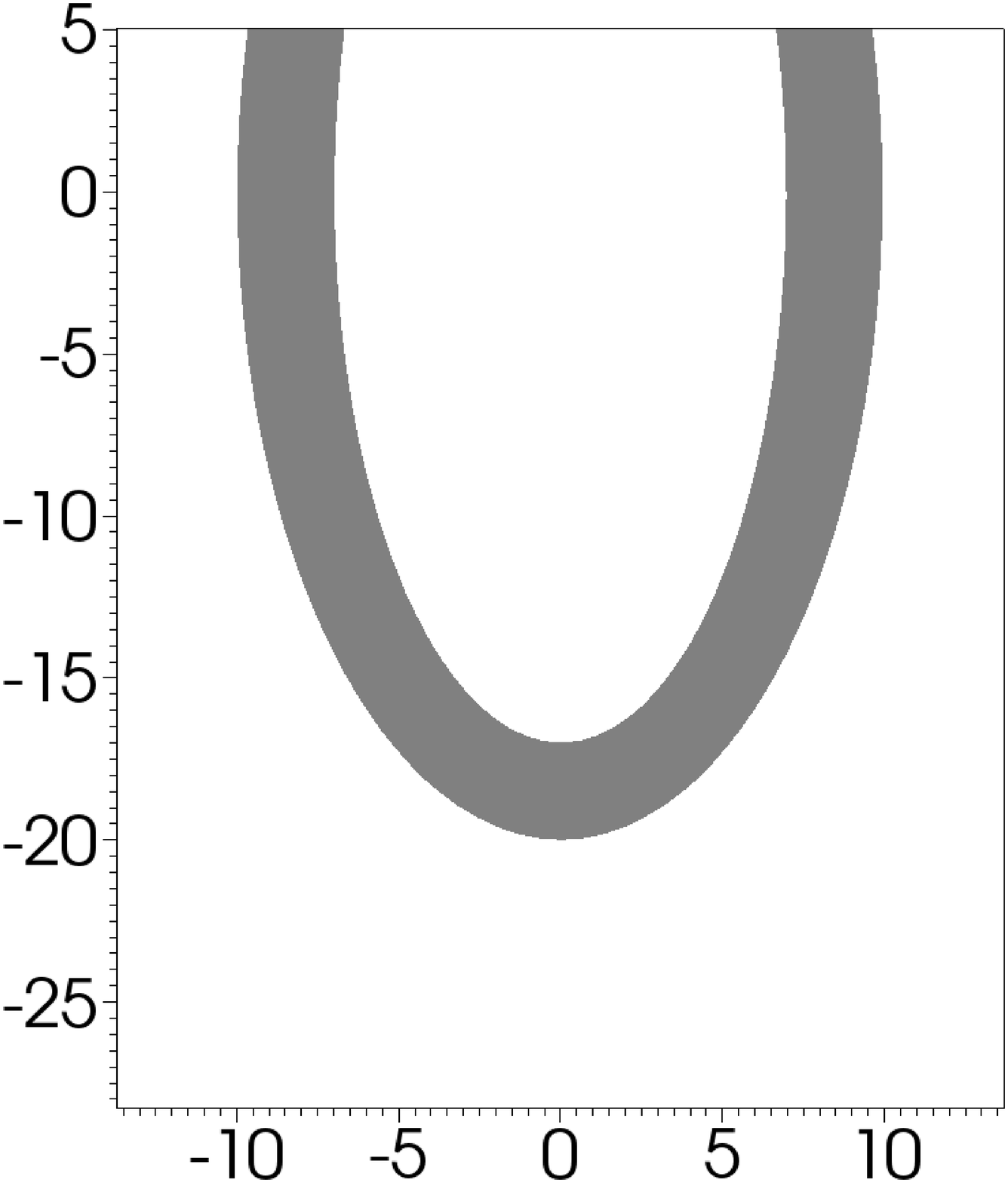} &	
	& \includegraphics[width=0.2375\textwidth]{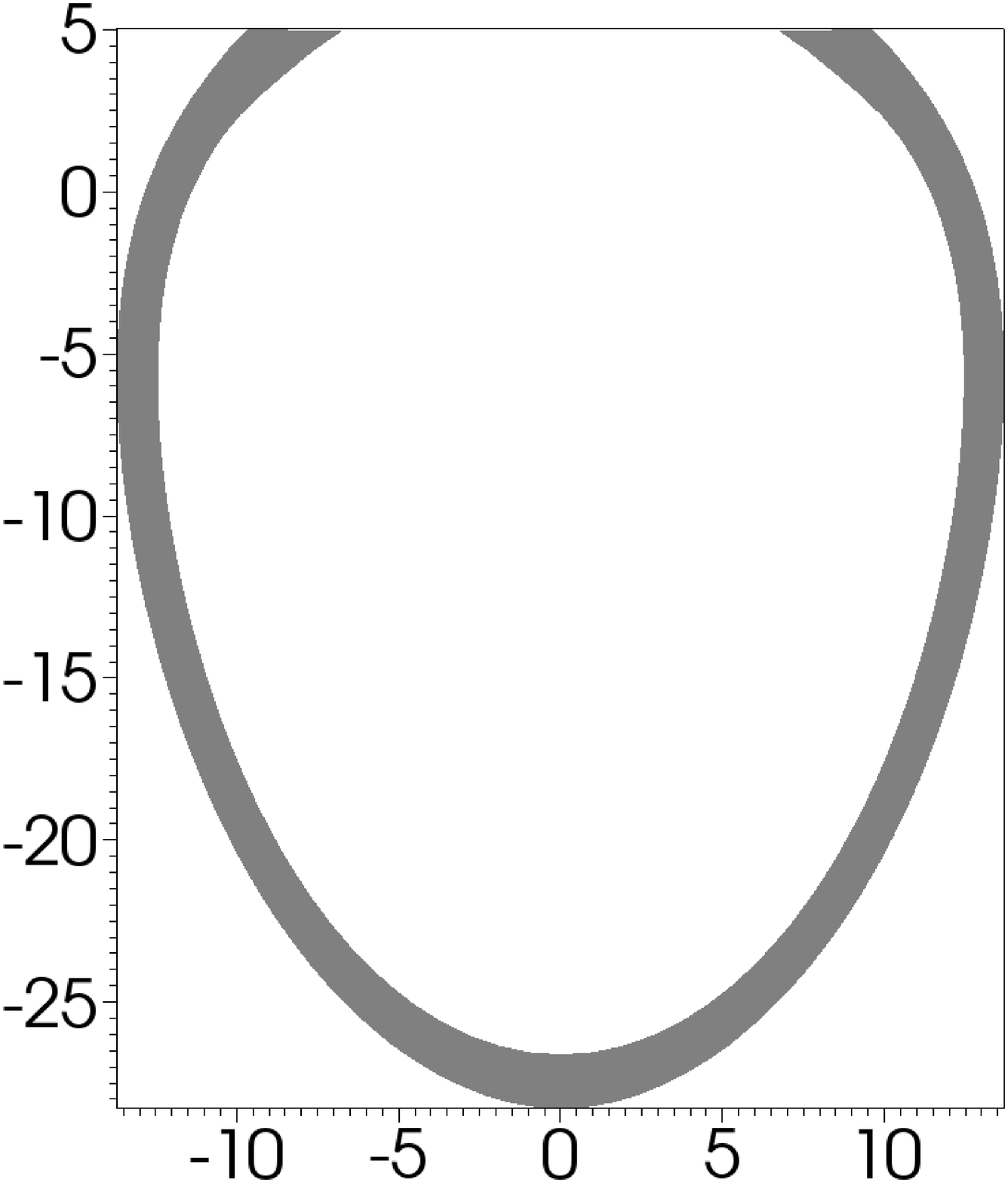} &	
	& \includegraphics[width=0.2375\textwidth]{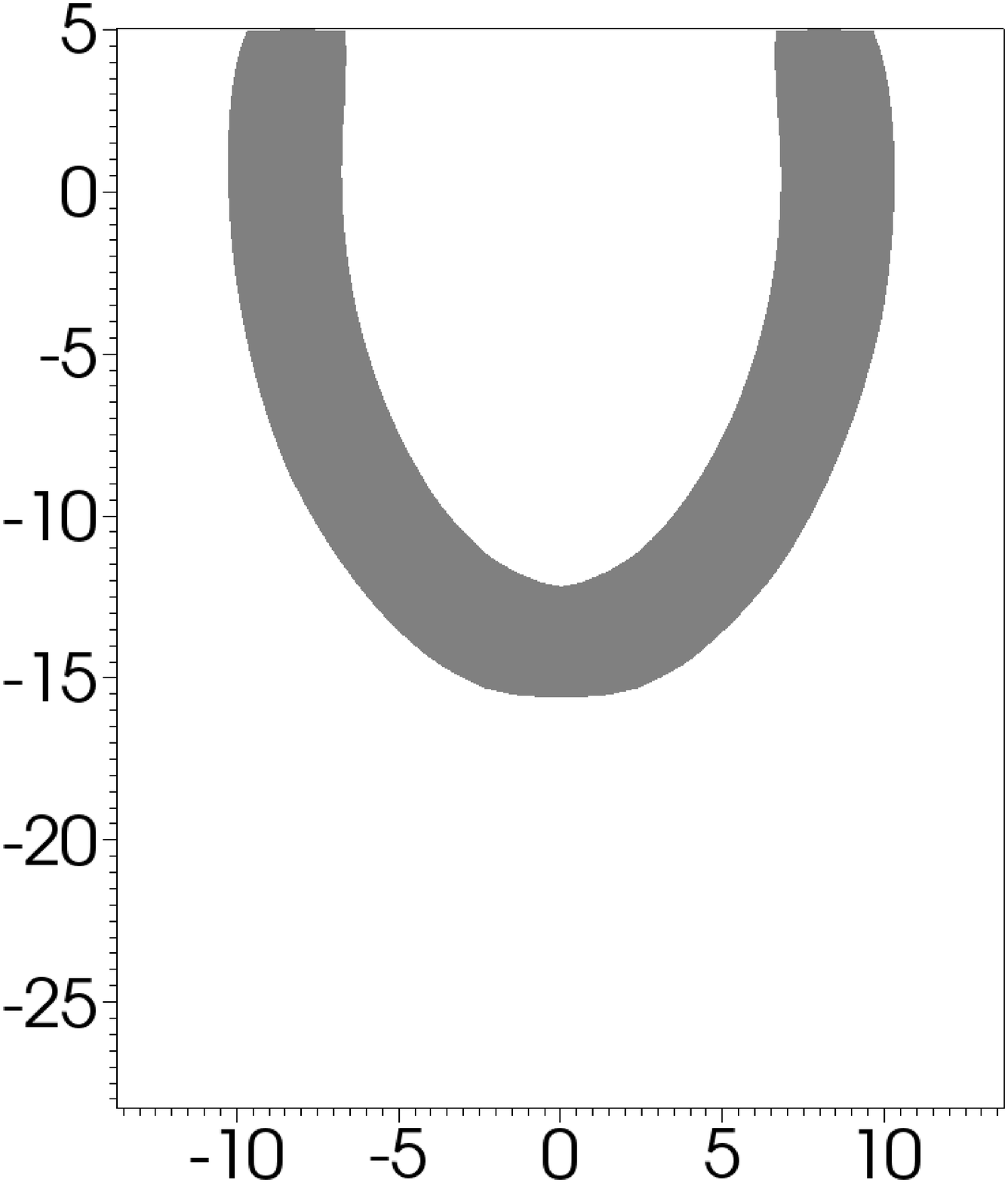}	
	\end{tabular}
	\caption{
		Inflation and active contraction of an idealized model of the left ventricle of the heart.
		{\bf A.}~The reference configuration of the left ventricle model is a truncated ellipsoid.
		{\bf B.}~Passive inflation of an isotropic version of the left ventricular model.
		{\bf C.}~Active contraction of a transversely isotropic left ventricular model that includes transmural fiber rotation.
		Notice the torsion induced by the active contractile forces.
		{\bf D--F.}~Slices of the configurations shown in panels {\bf A}--{\bf C}.
	}
	\label{f:idealized_LV}
\end{figure}

\begin{figure}
	\centering
	\small
	\begin{tabular}{lr}	
	\vspace{-10pt} {\bf A.} \quad \mbox{} & \\
	& \includegraphics[scale=0.3]{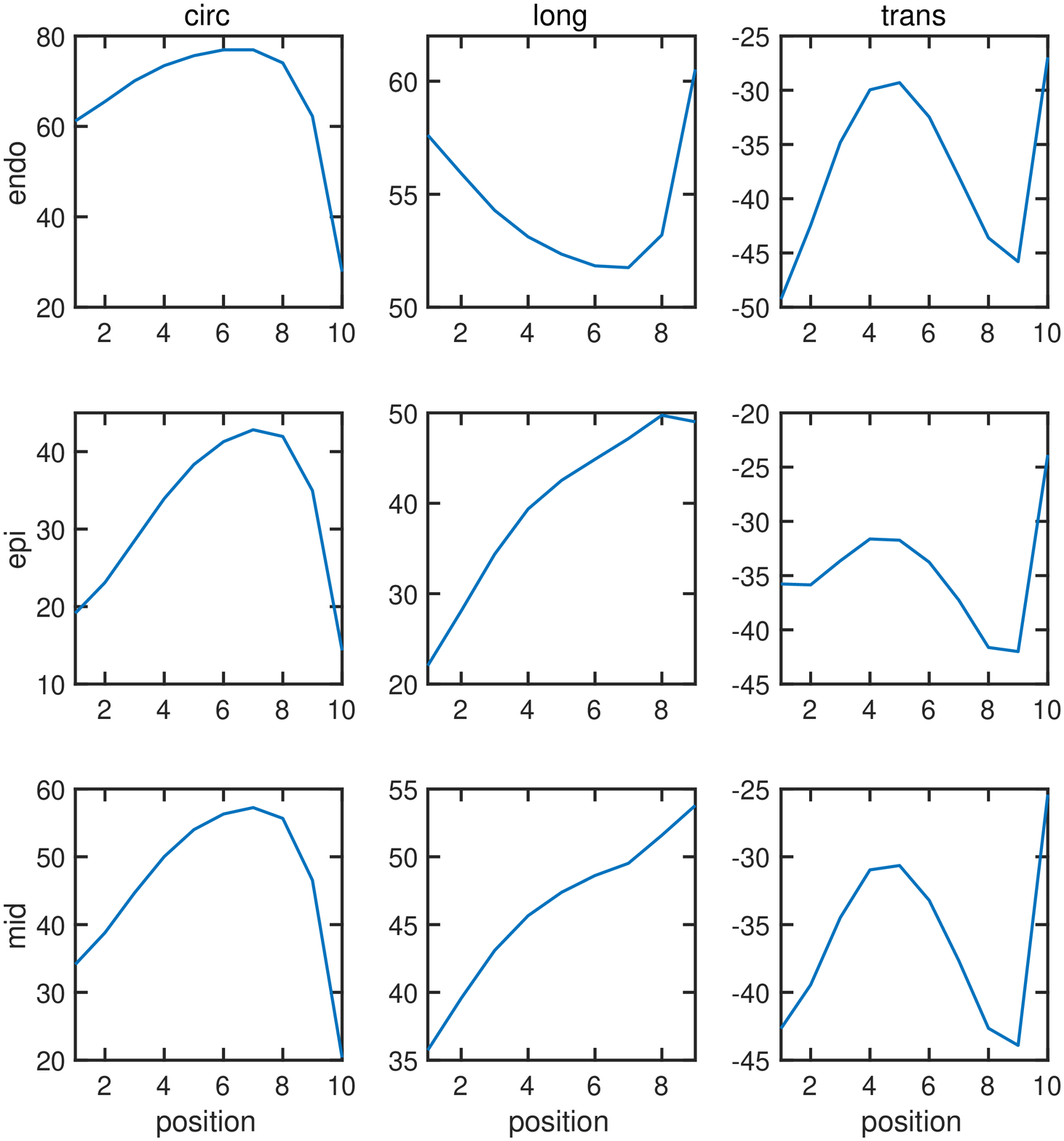} \\
	& \\
	\vspace{-10pt} {\bf B.} \quad \mbox{} & \\
	& \includegraphics[scale=0.3]{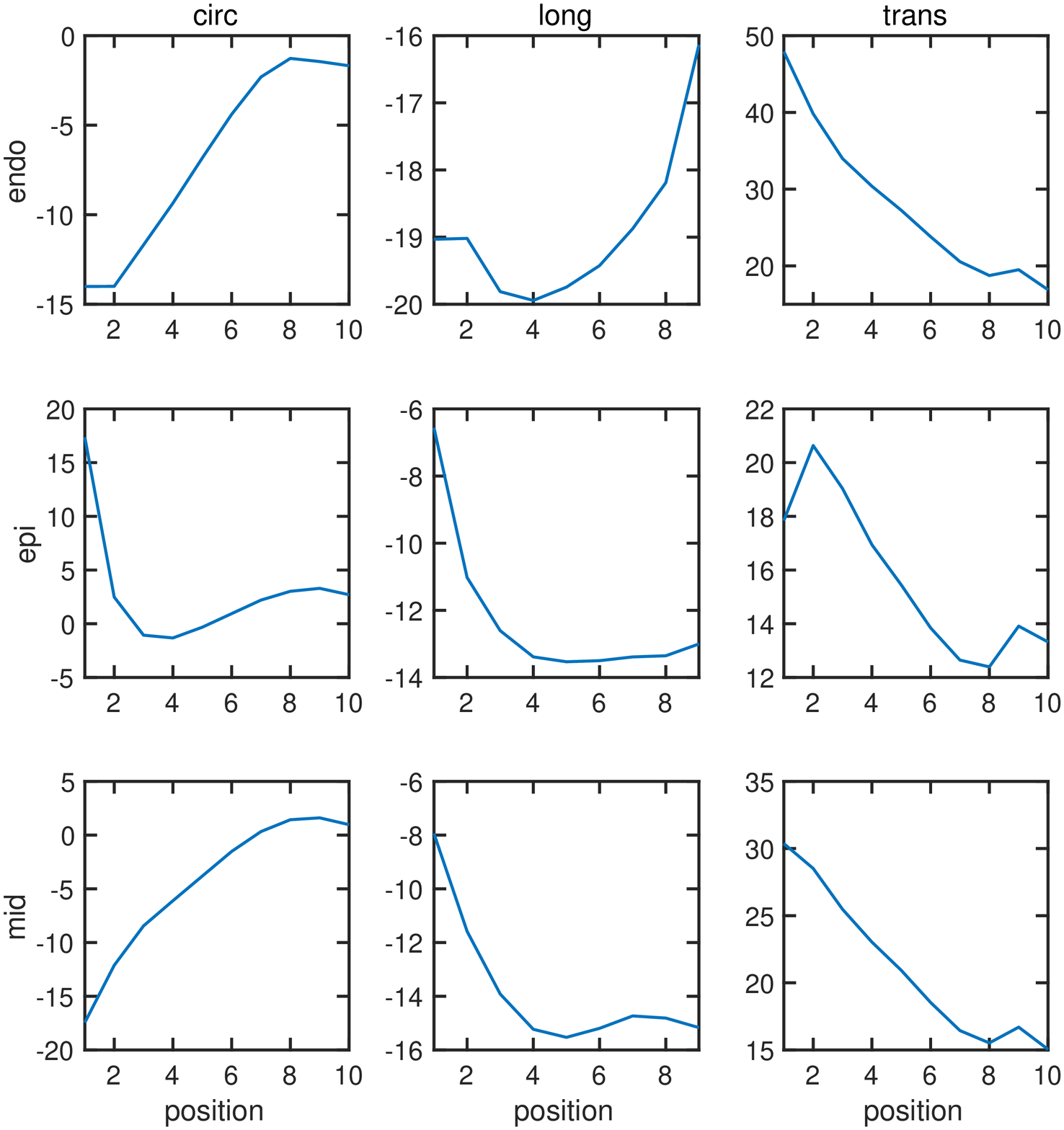}
	\end{tabular}
	\caption{
		Circumferential, longitudinal, and transverse strains along the endocardium, epicardium, and midwall in the idealized left ventricular model for ({\bf A}) passive inflation and ({\bf B}) active contraction.
		Strains are plotted at selected points running from apex to base.
		See Land et al.~\cite{SLand15-cardiac_verification} for additional details on the positions of the sample points.
		The computed strains are in excellent agreement (within $1\%$) with consensus results obtained by other structural mechanics codes using the same model \cite{SLand15-cardiac_verification}.
		}
	\label{f:LV_strains}
\end{figure}

To demonstrate the capabilities of the methodology to treat more complex geometries and structural models, this section briefly presents results from tests using an idealized model of the left ventricle of the heart, as used in a previous benchmarking study of cardiac mechanics codes by Land et al.~\cite{SLand15-cardiac_verification}.

For these tests, the undeformed geometry of the idealized left ventricle is described as a truncated ellipsoid.
Using parametric coordinates $(t,u,v)$, the ventricular geometry is defined by
\begin{equation}
	\s(t,u,v) = \left((7+3t) \sin u \cos v , (7+3t) \sin u \sin v , (17+3t) \cos u\right),
\end{equation}
with length measured in millimeters.
The endocardial surface is defined by $t = 0$, $u \in [-\pi,-\arccos \frac{5}{17}]$, and $v \in [-\pi,\pi]$, and the epicardial surface is defined by $t = 1$, $u \in [-\pi,-\arccos \frac{5}{20}]$, and $v \in [-\pi,\pi]$.
The model ventricle is truncated at the base, which is taken to correspond to $z = 5~\text{mm}$.
See fig.~\ref{f:idealized_LV}(A).
The passive elasticity of the heart wall is described using a transversely isotropic hyperelastic constitutive law by Guccione et al.~\cite{JMGuccione95}, which is defined with respect to a local fiber-aligned coordinate system via
\begin{equation}
	\We = \frac{C}{2}\left(e^Q - 1\right),
\end{equation}
with
\begin{equation}
	Q = b_\text{f} E_{11}^2 + b_\text{t} \left(E_{22}^2 + E_{33}^2 + E_{23}^2 + E_{32}^2\right)	+ b_\text{fs} \left(E_{12}^2 + E_{21}^2 + E_{13}^2 + E_{31}^2\right),
\end{equation}
in which $\EE = \half(\FF^\text{T} \FF - \II) = \left(E_{ij}\right)$ is the Green-Lagrange strain tensor and $C$, $b_\text{f}$, $b_\text{t}$, $b_\text{fs}$ are material parameters.
The structure is discretized using trilinear ($Q^1$) hexahedral elements.
The computational domain $\Omega$ is a $4~\text{cm} \times 4~\text{cm} \times 4~\text{cm}$ region.
Except where noted, $\Omega$ is discretized using a uniform $64 \times 64 \times 64$ Cartesian grid.
The model is run to steady state under steady loading conditions, corresponding to a quasi-static analysis, as in earlier work \cite{HGao14-iblv_diastole}.

\subsubsection{Passive inflation of an isotropic left ventricle model}

Figs.~\ref{f:idealized_LV}(B,E) and \ref{f:LV_strains}(A) show results of passively inflating an isotropic version of the left ventricle model.
We use $C = 10~\text{kPa}$ and $b_\text{f} = b_\text{t} = b_\text{fs} = 1$ along with a uniform pressure load of $10~\text{kPa}$.
This corresponds to `Problem 2' from Land et al.~\cite{SLand15-cardiac_verification}.
Notice that our formulation permits very large structural deformations.
The computed longitudinal, circumferential, and radial strains in the inflated configuration are in excellent agreement (within $1\%$) with values obtained by other structural mechanics codes using the same model \cite{SLand15-cardiac_verification}.
This test problem was judged to be relatively easy, and all of the methods considered in the benchmarking study yielded virtually identical results.

\subsubsection{Active contraction of a fiber-reinforced left ventricle model}

\begin{figure}
	\centering
	\small
	\begin{tabular}{lclc}
	\vspace{-10pt} {\bf A.} & & {\bf B.} \\
	& \includegraphics[width=0.4\textwidth]{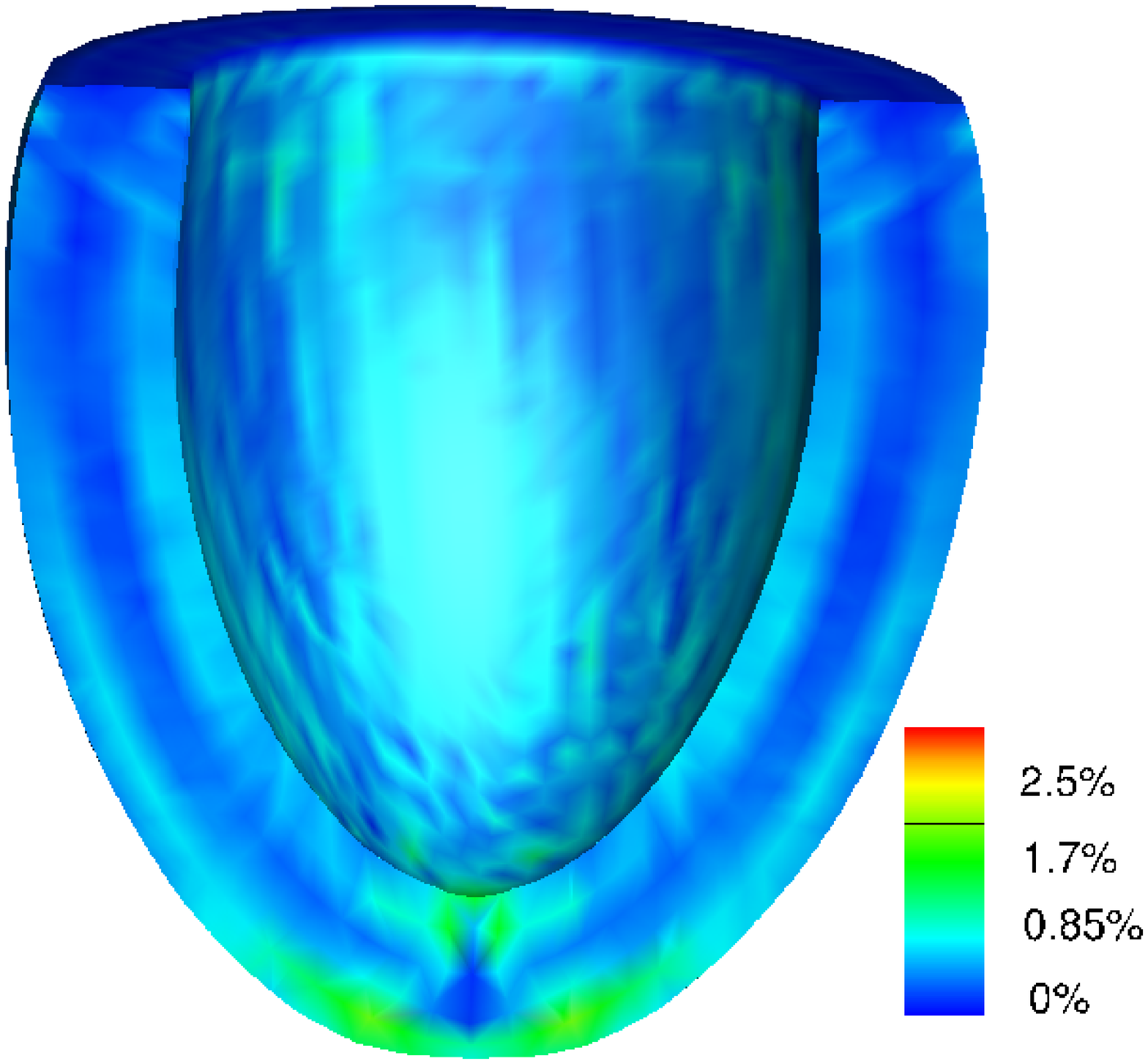} &
	& \includegraphics[width=0.4\textwidth]{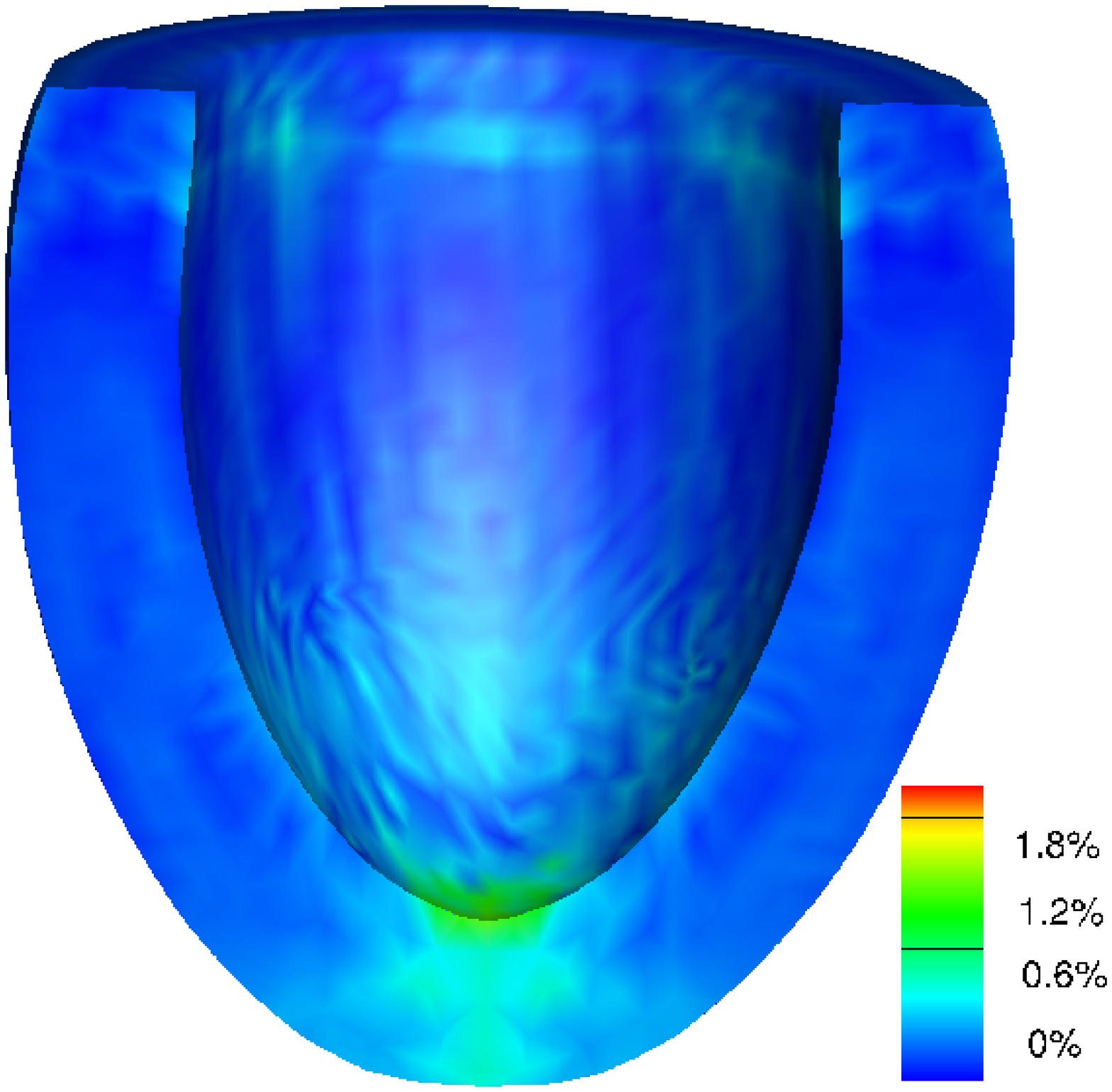} \\
	\end{tabular}
	\caption{
		Difference in the displacement of the actively contracting fiber-reinforced left ventricle model using different $N \times N \times N$ Cartesian grids.
		{\bf A.}~Distribution of difference in the structural displacement obtained using $N = 48$ and $N = 64$.
		{\bf B.}~Distribution of difference in the structural displacement obtained using $N = 64$ and $N = 96$.
	}
	\label{f:LV_diff}
\end{figure}

\begin{figure}
	\centering
	\includegraphics[width=0.4\textwidth]{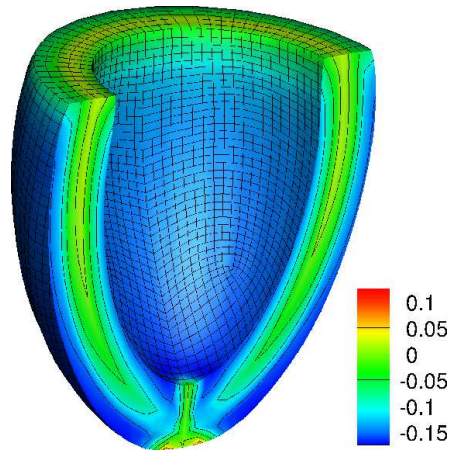} 
	\caption{
		Fiber strain distribution in the actively contracting fiber-reinforced left ventricle model obtained using a $64 \times 64 \times 64$ Cartesian grid.
	}
	\label{f:LV_active_strains}	
\end{figure}

Figs.~\ref{f:idealized_LV}(C,F) and \ref{f:LV_strains}(B) show the results from a simulation of active contraction using this model.
In this case, we use a fiber-reinforced material model, in which the fiber direction in the reference configuration $\E_\text{f}$ is given by
\begin{equation}
	\E_\text{f} = \sin \alpha \left. \D{\s}{u} \middle/ \left\|\D{\s}{u}\right\| \right. + \cos \alpha \left. \D{\s}{v} \middle/ \left\|\D{\s}{v}\right\| \right., \quad \alpha = 90 - 180 t,
\end{equation}
yielding a 180$^\circ$ transmural fiber rotation with a linear relationship between fiber angle and tissue depth.
In addition to passive elasticity, the test also includes active contractile forces, which are accounted for in the first Piola-Kirchhoff stress via
\begin{equation}
	\PPe = \D{W}{\FF} + T_\text{a} \, \FF \, \E_\text{f} \otimes \E_\text{f},
\end{equation}
with $T_\text{a}$ corresponding to the active stress acting in the fiber direction.
We use $C = 2~\text{kPa}$, $b_\text{f} = 8$, $b_\text{t} = 2$, $b_\text{fs} = 4$, and $T_\text{a} = 60~\text{kPa}$.
A uniform pressure load of $15~\text{kPa}$ is also applied along the endocardial surface.
This corresponds to `Problem 3' from Land et al.~\cite{SLand15-cardiac_verification}.
No special treatment is used for the fiber singularity at $u = -\pi$.
To demonstrate convergence, the differences in the structural displacement obtained using uniform $N \times N \times N$ grids for $N = 48$, $64$, and $96$ are shown in fig.~\ref{f:LV_diff}.
The maximum difference in the displacement between $N=48$ and $N = 64$ is 2.5\%; this difference is 1.8\% between $N = 64$ and $N=96$, indicating approximately first-order convergence.
This test problem was judged to be relatively difficult when compared to Problem 2 \cite{SLand15-cardiac_verification}, but as in the case of passive deformations, the actively contracting model is in excellent agreement (within $1\%$) with grid-converged consensus results obtained by other structural codes using the same idealized left ventricular model\cite{SLand15-cardiac_verification}.
The distribution of fiber strain  within the model is shown in fig.~\ref{f:LV_active_strains}.
It is interesting to see that although the ventricular wall is mostly compressed, as indicated by the negative values of the fiber strain, there are parts of the wall that are very slightly stretched, primarily in the midwall.
This is because the fibers rotate across the wall, and the middle layer experiences circumferential stretch even though the whole left ventricle is substantially compressed in the longitudinal direction.
The comparatively small stretches near the apex reflect the different fiber distribution in that region.

\section{Discussion and Conclusions}
\label{s:conclusions}

This paper describes a version of the IB method for fluid-structure interaction that uses either a structured or unstructured FE discretization of the immersed structure while retaining a Cartesian grid finite difference scheme for the Eulerian equations.
This method has already proved useful in a variety of biological and biomedical applications \cite{BEGriffith13-FMD2013, HGao14-iblv_diastole, HGao14-iblv_mi, HGao14-MV, WWChen16-coupled_lv, VFlamini16-aortic_root, SKJones16-bristles, APHoover17-active_jellyfish, HGaoXX-coupled_mv_lv, AHasanXX-imaged_based_aortic_root, WKouXX-continuum_eso}, but its performance has not previously been examined using standard benchmark cases such as those that are the focus of this work.
This paper also demonstrates the suitability of this method to simulate cardiac mechanics using an idealized model of the left ventricle of the heart \cite{SLand15-cardiac_verification}.

A feature of this method is that it uses standard discretization technology for both the Lagrangian and Eulerian equations.
In practice, it should be feasible to use this approach to couple existing structural analysis and incompressible flow codes by passing forces from the structural code to the fluid solver, and by passing velocities from the fluid solver back to the structural code.
Although we focus on cases in which the immersed structure is hyperelastic, our numerical scheme does not rely on the availability of a hyperelastic energy functional, and the present approach is also demonstrated to treat immersed rigid bodies and structures with active contractile stresses.

A key contribution of this paper is that it introduces an approach to Lagrangian-Eulerian interaction that overcomes a longstanding limitation of the IB method, namely the requirement that the Lagrangian mesh must be relatively fine compared to the background Eulerian grid to avoid leaks at fluid-structure interfaces.
This restriction results in structural meshes that are excessively dense, potentially leading to reduced efficiency and even reduced accuracy.
Numerical examples demonstrate that our scheme permits the use of Lagrangian meshes that are at least four times as coarse as the background Eulerian grid without leak, and we expect that there are cases in which even coarser structural meshes would yield good accuracy.
Even with very coarse structural discretizations, the present methodology can yield improved volume conservation when compared to the IFE method \cite{XSWang10}.
Of course, fine Lagrangian meshes are still needed in practice to resolve fine-scale geometrical features, and an approach such as Lagrangian adaptive mesh refinement may be crucial for accurately treating extremely large structural deformations.
Nonetheless, the present approach enables the effective use of spatial discretizations that are tailored to the requirements of the structural model rather than dictated by the background Eulerian discretization.

Our method is based on a formulation of the continuous IB equations introduced by Boffi et al.~\cite{DBoffi08}.
We consider two different statements of these equations that each use a weak formulation of the Lagrangian equations of nonlinear elasticity.
One of these formulations treats the internal and boundary stresses within the immersed elastic structure using a single, unified, volumetric elastic force density.
This form of the equations of motion is essentially that used in the IFE method \cite{LZhang04,WKLiu06,LTZhang07,XSWang10,XSWang12-semi_implicit_IFE} as well as other similar approaches \cite{Peskin-energy-functions,DDevendran12,DBoffi08,Heltai12,SRoy15-benchmarking_ibfe}.
The other formulation considered in this work partitions these stresses into a internal elastic force density supported on the interior of the structure, and a transmission elastic force density supported on fluid-structure interfaces; this approach does not appear to have been used previously to develop a numerical scheme.
Both formulations are demonstrated to yield similar convergence rates, but the partitioned formulation is seen to yield higher accuracy for Lagrangian meshes that are relatively coarse in comparison to the background Eulerian grid, especially in terms of volume conservation.

A limitation of the partitioned formulation is that it does not satisfy a discrete power identity that implies that energy is conserved during Lagrangian-Eulerian interaction.
Such a power identity is satisfied by the unified formulation, and may be necessary to obtain an unconditionally stable implicit time-stepping scheme \cite{Newren07}.
Developing a symmetric partitioned formulation will likely require the introduction of additional boundary degrees of freedom, so that volumetric operators (i.e.~$\cS$ and $\cJ$) couple the volumetric structural variables to the background grid, and corresponding surface operators (e.g.~$\cS^{\p U}$ and $\cJ^{\p U}$) couple the surface degrees of freedom to the Eulerian variables.
Lagrange multipliers or penalty methods could be used to ensure that the surface discretization moves with the immersed body.

For immersed \emph{rigid} bodies, tests suggest that relatively coarse structural discretizations can yield superior accuracy when compared to relatively fine discretizations.
For the standard test case of viscous flow past a circular cylinder, we demonstrate that relatively fine structural discretizations can produce spurious drag, suggesting that the effective numerical size of the immersed body is determined in part by the spacing of the control points.
Although not studied here in detail, there appears to be a relatively sharp transition, and the results obtained by the method for structural meshes that are \emph{coarser} than a minimum spacing appear to be relatively insensitive to the choice of Lagrangian meshwidth.
The threshold spacing appears to depend on the choice of regularized kernel function.
For instance, when using a three-point kernel function, the method yields good results for relative mesh spacing values of $\Mfac = 1$, whereas this same relative structural mesh spacing yields poor results with four-~and six-point kernels.

In a constrained formulation, using a mesh of control points that is too fine with respect to the background grid will result in an ill-conditioned system of equations, and for a sufficiently fine mesh of control points, the constrained formulation becomes singular.
It is interesting to note that in our formulation, we obtain high accuracy even with extremely dense meshes of \emph{interaction} points so long as the \emph{control} points are sufficiently far apart from each other.
Projecting the standard IB velocity field $\U^\text{IB}$ onto the FE shape functions \emph{filters} velocity fluctuations at length scales that are smaller than structural mesh width.
For relatively coarse structural meshes, this additional filtering appears to improve the accuracy of the IB method, dramatically reducing errors in the lift and drag forces.
Such errors have been previously described to be a fundamental aspect of the IB method, but we believe that our results show that such numerical artifacts are not intrinsic to this methodology, but rather result from the use of overly dense structural meshes.
The present work provides a framework for coupling relatively course structural models without leak at fluid-structure interfaces.

In closing, we remark that the partitioned formulation developed in this work could be useful in constructing higher-order versions of the IB method.
For instance, this formulation is well suited for developing a hybrid approach in which the IB method is used to treat the volumetric internal force density, and another method is used to treat the singular transmission force density.
Because the transmission force density of the partitioned formulation is defined on a closed surface, it is feasible to treat it with higher-order accuracy by a version of the immersed interface method \cite{LeVequeLi97,LiLai01,LeeLeVeque03,SXu06,SXu06b,SXu08,SXu11-3d_iim}.
Such an extension of this method could yield a fully second-order accurate generalization of the IB method for cases, like those considered in this work, in which the immersed structure is of codimension zero with respect to the fluid.

\section*{Acknowledgements}

This work was sponsored in part by an award from the Royal Society of Edinburgh.
B.E.G.~acknowledges research support from the National Science Foundation (NSF awards ACI 1047734/1460334, ACI 1450327, and DMS 1016554/1460368) and the National Institutes of Health (NIH award HL117063).
He also gratefully acknowledges discussions of this work with Charles Peskin and David McQueen.
X.Y.L.~acknowledges research support from the Engineering and Physical Sciences Research Council (UK EPSRC awards EP/I029990 and EP/N014642/1) and from the Leverhulme Trust (RF-2015-510).
We also acknowledge Hao Gao (U.~Glasgow) and Viatcheslav Gurev (IBM T.~J.~Watson Research Center) for their assistance with the cardiac mechanics benchmark.

\appendix

\section{Temporal discretization}
\label{s:temporal_discretization}

\subsection{Basic time-stepping scheme}
\label{s:basic_time_stepping_scheme}

Let $\X^{n}$, $\u^{n}$, and $p^{n-\half}$ denote the approximations to the values of $\X$ and $\u$ at time $t^{n}$, and to the value of $p$ at time $t^{n-\half}$, respectively.
We advance the solution values forward in time by the increment $\dt$ as follows.
First, we determine a preliminary approximation to the deformed structure configuration at time $t^{n+\half}$ via
\begin{equation}
  \frac{{\X}^{n+\half} - \X^{n}}{\dt/2} = \cJ\left(\X^{n}\right) \u^{n}.
\end{equation}
Then, we solve
\begin{align}
  \rho \left(\frac{\u^{n+1} - \u^{n}}{\dt} + \vec{A}^{n+\half}\right) &= -\grad_h p^{n+\half} + \mu \grad_h^2 \left(\frac{\u^{n+1} + \u^{n}}{2}\right) + \vec{f}^{n+\half} + \vec{t}^{n+\half}, \label{e:staggered1} \\
  \grad_h \cdot \u^{n+1} &= 0, \label{e:staggered2} \\
  \f^{n+\half} &= \cS\left(\X^{n+\half}\right) \, \F\left(\X^{n+\half}\right), \label{e:staggered3} \\
  \t^{n+\half} &= \cS^{\p U}\left(\X^{n+\half}\right) \, \T\left(\X^{n+\half}\right), \label{e:staggered4} \\
  \frac{\X^{n+1} - \X^{n}}{\dt} &= \cJ\left(\X^{n+\half}\right) \, \left(\frac{\u^{n+1}+\u^{n}}{2}\right), \label{e:staggered5}
\end{align}
for $\X^{n+1}$, $\u^{n+1}$, and $p^{n+\half}$, in which $\vec{A}^{n+\half} = \frac{3}{2} \u^{n} \cdot \grad_h \u^{n} - \frac{1}{2} \u^{n-1} \cdot \grad_h \u^{n-1}$ is computed via a PPM-type approximation to the nonlinear advection term \cite{BEGriffith09-efficient,BEGriffith12-ib_volume_conservation}.
The time-stepping scheme for the unified formulation is analogous.
Notice that solving eqs.~\eqref{e:staggered1}--\eqref{e:staggered5} for $\X^{n+1}$, $\u^{n+1}$, and $p^{n+\half}$ requires the solution of a Crank-Nicolson-type discretization of the time-dependent incompressible Stokes equations.
We solve this system of equations via the flexible GMRES (FGMRES) algorithm \cite{Saad93}, using $\vec{u}^{n}$ and $p^{n-\half}$ as initial approximations to $\vec{u}^{n+1}$ and $p^{n+\half}$, and using a pressure-free projection method with inexact multigrid subdomain solvers as a preconditioner \cite{BEGriffith09-efficient}.

\subsection{Initial time step}
\label{e:initial_time_step}

Because time step-lagged values of $\u$ and $p$ are used by the time-stepping scheme of sec.~\ref{s:basic_time_stepping_scheme}, we cannot use that scheme for the initial time step.
Instead, we use a two-step predictor-corrector method:
First, we solve
\begin{align}
  \rho \left(\frac{\tilde{\u}^{n+1} - \u^{n}}{\dt} + \vec{A}^{n}\right) &= -\grad_h \tilde{p}^{n+\half} + \mu \grad_h^2 \left(\frac{\tilde{\u}^{n+1} + \u^{n}}{2}\right) + \vec{f}^{n} + \vec{t}^{n}, \\
  \grad_h \cdot \tilde{\u}^{n+1} &= 0, \\
  \f^{n} &= \cS\left(\X^{n}\right) \, \F\left(\X^{n}\right), \\
  \t^{n} &= \cS^{\p U}\left(\X^{n}\right) \, \T\left(\X^{n}\right), \\
  \frac{\tilde{\X}^{n+1} - \X^{n}}{\dt} &= \cJ\left(\X^{n}\right) \, \u^{n},
\end{align}
for $\tilde{\X}^{n+1}$, $\tilde{\u}^{n+1}$, and $\tilde{p}^{n+\half}$, with $\vec{A}^{n} = \u^{n} \cdot \grad_h \u^{n}$.
Because we do not have an initial value for the pressure, we use $p=0$ as an initial guess for $p^{n+\half}$.
Second, we set $\X^{n+\half} = \frac{\tilde{\X}^{n+1} + \X^{n}}{2}$ and solve eqs.~\eqref{e:staggered1}--\eqref{e:staggered5} for $\X^{n+1}$, $\u^{n+1}$, and $p^{n+\half}$, except that we use $\vec{A}^{n+\half} = \u^{n+\half} \cdot \grad_h \u^{n+\half}$ with $\u^{n+\half} = \half\left(\tilde{\u}^{n+1}+\u^{n}\right)$.

\bibliographystyle{wileyj}
\bibliography{ibfem}

\end{document}